\documentclass[a4paper,10pt]{amsart}
\usepackage{microtype}
\usepackage[utf8]{inputenc}
\usepackage{amsfonts, amsmath, amssymb, amsthm, mathtools}
\usepackage{xcolor}
\definecolor{darkblue}{rgb}{0,0,0.6}
\usepackage{tikz}
\usetikzlibrary{matrix,arrows,trees}
\usepackage{tikz-cd}
\usepackage[ocgcolorlinks,colorlinks=true, citecolor=darkblue, filecolor=darkblue, linkcolor=darkblue, urlcolor=darkblue]{hyperref}
\usepackage[nameinlink, capitalise]{cleveref}

\crefname{rem}{Remark}{Remarks}

\DeclareRobustCommand{\SkipTocEntry}[5]{}

\usepackage{eucal}                  

\hypersetup{bookmarksopen=true}

\newcommand*\cocolon{%
	\nobreak
	\mskip6mu plus1mu
	\mathpunct{}%
	\nonscript
	\mkern-\thinmuskip
	{:}%
	\mskip2mu
	\relax
}

\numberwithin{equation}{section}

\swapnumbers
\theoremstyle{definition}
\newtheorem{defn}[equation]{Definition}
\newtheorem{ex}[equation]{Example}

\newtheorem{rem}[equation]{Remark}
\newtheorem{theorem}[equation]{Theorem}
\newtheorem{prop}[equation]{Proposition}
\newtheorem{lem}[equation]{Lemma}
\newtheorem{cor}[equation]{Corollary}

\renewcommand{\epsilon}{\varepsilon}
\renewcommand{\theta}{\vartheta}
\renewcommand{\phi}{\varphi}

\newcommand{\bbE}{\mathbb{E}}

\newcommand{\bbN}{\mathbb{N}}

\newcommand{\bbZ}{\mathbb{Z}}

\newcommand{\cA}{\mathcal{A}}
\newcommand{\cB}{\mathcal{B}}
\newcommand{\cC}{\mathcal{C}}
\newcommand{\cD}{\mathcal{D}}

\newcommand{\cQ}{\mathcal{Q}}

\newcommand{\cU}{\mathcal{U}}
\newcommand{\cV}{\mathcal{V}}

\newcommand{\cX}{\mathcal{X}}
\newcommand{\cY}{\mathcal{Y}}
\newcommand{\cZ}{\mathcal{Z}}

\newcommand{\abs}[1]{\lvert #1 \rvert}

\DeclareMathOperator*{\colim}{colim}
\DeclareMathOperator{\cofib}{cofib}
\DeclareMathOperator{\fib}{fib}
\DeclareMathOperator{\Fun}{Fun}

\DeclareMathOperator{\Hom}{Hom}
\DeclareMathOperator{\id}{id}
\DeclareMathOperator{\inc}{inc}

\DeclareMathOperator{\kernel}{ker}
\newcommand{\Ldec}{\mathrm{L}}
\newcommand{\Rdec}{\mathrm{R}}

\DeclareMathOperator{\pr}{pr}

\newcommand{\st}{\mathrm{st}}

\newcommand{\Ab}{\mathrm{Ab}}
\newcommand{\catinf}{\mathrm{Cat}}

\newcommand{\catst}{\catinf^{\st}}

\newcommand{\Funex}[1]{\Fun^{\mathrm{ex}}_{#1}}
\newcommand{\op}{\mathrm{op}}

\newcommand{\psh}{\mathcal{P}}
\newcommand{\Set}{\mathrm{Set}}
\newcommand{\Sp}{\mathrm{Sp}}
\newcommand{\Spc}{\mathrm{An}}

\newcounter{commentcounter}

\newcommand{\bdd}{\mathrm{b}}

\newcommand{\egr}{\mathrm{eg}}
\newcommand{\exct}{\mathrm{ex}}

\newcommand{\fg}{\mathrm{fg}}

\newcommand{\ing}{\mathrm{in}}
\newcommand{\lex}{\mathrm{lex}}

\newcommand{\mx}{\mathrm{max}}

\newcommand{\rezk}{\mathrm{Rezk}}

\newcommand{\s}{\mathrm{s}}

\newcommand{\tf}{\mathrm{tf}}
\newcommand{\tors}{\mathrm{tors}}

\newcommand{\W}{\mathrm{W}}

\newcommand{\cc}{\mathrm{cc}}
\newcommand{\derb}{\cD^{\bdd}}
\newcommand{\exact}{\mathrm{Exact}}

\newcommand{\idem}[1]{{#1}^\natural}
\newcommand{\widem}[1]{{#1}^\flat}
\newcommand{\Idemrel}[2]{#2^{\natural{#1}}}
\newcommand{\LFib}{\mathrm{LFib}}
\newcommand{\Mod}{\mathrm{Mod}}

\newcommand{\Perf}{\mathrm{Perf}}

\newcommand{\quotex}[2]{({#1}/{#2})_\exct}
\newcommand{\quotw}[2]{({#1}/{#2})_\W}

\newcommand{\wald}{\mathrm{Wald}}

\DeclareMathOperator{\Ar}{Ar}

\DeclareMathOperator{\con}{con}
\DeclareMathOperator{\dec}{dec}
\DeclareMathOperator{\K}{K}

\DeclareMathOperator{\N}{N}

\DeclareMathOperator{\stab}{SW}
\DeclareMathOperator{\St}{St}
\DeclareMathOperator{\Tw}{Tw}
\DeclareMathOperator{\Un}{Un}

\DeclareFontFamily{U}{min}{}
\DeclareFontShape{U}{min}{m}{n}{<-> dmjhira}{}
\newcommand{\yo}{\text{\usefont{U}{min}{m}{n}\symbol{'110}}}

\title[Localisation theorems for exact categories]{Localisation theorems for the connective K-theory of exact categories}

\author{Christoph Winges}
\address{Fakult\"at f\"ur Mathematik, Universit\"at Regensburg, 93040 Regensburg, Germany}
\email{christoph.winges@ur.de}

\date{}
\keywords{algebraic K-theory, localisation theorem, theorem of the heart}
\subjclass{18F25; 18E10, 18N55, 19D10}

\begin{document}

\begin{abstract}
 We prove a localisation theorem for the K-theory of filtering subcategories of exact $\infty$-categories which subsumes the localisation theorem for stable $\infty$-categories, Quillen's localisation theorem for abelian categories, and Schlichting's localisation theorem for s-filtering subcategories. 
\end{abstract}

\maketitle
\tableofcontents

\section{Introduction}

Through the work of Blumberg, Gepner and Tabuada \cite{bgt:alg-k} and Barwick \cite{barwick:algK}, algebraic K-theory has been understood to be the universal space-valued additive invariant, either as a functor on stable $\infty$-categories, or more generally on Waldhausen $\infty$-categories.
One of its key structural properties, which is not immediate from the universal property, is the fact that K-theory is localising in the sense that certain non-split cofibre sequences $\cA \to \cB \to \cC$ of input categories induce fibre sequences $\K(\cA) \to \K(\cB) \to \K(\cC)$ in algebraic K-theory.

For stable $\infty$-categories, this is encoded in the notion of a Verdier-localising invariant \cite[Theorem~6.1]{hls:localisation}, which is rooted in Waldhausen's fibration theorem \cite[Theorem~1.6.4]{waldhausen:algKspaces}.
The latter also has an analogue in Barwick's setup \cite[Theorem~9.24]{barwick:algK}, at the additional expense of having to deal with labelled Waldhausen $\infty$-categories.

As a middle ground between stable and Waldhausen $\infty$-categories, Barwick introduced in \cite{barwick:heart} the class of exact $\infty$-categories, which encompasses both stable $\infty$-categories and Quillen's exact categories \cite{quillen:k-theory}.
Due to the Gillet--Waldhausen theorem (see \cite[Theorem~1.11.7]{TT} and \cite[Theorem~1.7]{sw:exact-stable}), algebraic K-theory also enjoys a universal property as a functor on exact $\infty$-categories \cite[Theorem~6.7]{sw:exact-stable}.

However, even in the setting of ordinary exact categories, there are several seemingly disconnected localisation theorems, among them Quillen's localisation theorem for abelian categories \cite[Theorem~5]{quillen:k-theory} and Schlichting's localisation theorem for s-filtering subcategories \cite[Theorem~2.1]{schlichting:delooping}.
The goal of this article is to prove a general localisation theorem for exact $\infty$-categories which unifies all of these statements (with the exception of Waldhausen's fibration theorem, which doesn't concern exact categories).

There is a general proof strategy for localisation theorems in algebraic K-theory using Waldhausen's generic fibration theorem \cite[Corollary~1.5.7]{waldhausen:algKspaces} (see \cite[Theorem~8.11]{barwick:algK} for an $\infty$-categorical version).
This strategy is implicit in Waldhausen's proof of the fibration theorem, but also in work of Staffeldt \cite{staffeldt:fundamental}, and has been further clarified in \cite{hls:localisation} based on unpublished work of Clausen.
A convenient formulation of this strategy for the purpose at hand can be found in \cref{sec:generic-fibration}, and uses the observation that there exist in fact two models for the generic cofibre, depending on the preferred choice of d\'ecalage for simplicial objects.

Given an inclusion $j \colon \cU \to \cC$ of a suitable full subcategory, the explicit form of the two generic cofibres suggests that, if the relative K-theoretic term can be identified as the K-theory of some category, it should be closely related to the localisation of $\cC$ at the cofibres of ingressive morphisms $u \rightarrowtail x$ with $u \in \cU$, or to the localisation of $\cC$ at the fibres of egressive morphisms $y \twoheadrightarrow w$ with $w \in \cU$. The latter makes most sense if $\cC$ is assumed to be exact.
However, these localisations have a better chance of modelling the cofibre of $j$ among Waldhausen respectively coWaldhausen $\infty$-categories.

\cref{sec:quotients} singles out a condition which guarantees that such a localisation models the desired cofibre.
Specifically, we introduce the following version of Schlichting's notion of a filtering subcategory \cite{schlichting:delooping}.

\begin{defn}\label[defn]{def:filtering-intro}
 Let $\cC$ be an exact $\infty$-category and let $\cU \subseteq \cC$ be an extension-closed subcategory.
 
 Then $\cU$ is \emph{left filtering} if the following holds for every $u \in \cU$ and every $x \in \cC$:
 \begin{itemize}
     \item For every $\alpha \in \pi_k\Hom_\cC(u,x)$, $k \geq 0$, there exists an egressive morphism $p \colon x \twoheadrightarrow x'$ whose fibre lies in $\cU$ and such that $p_*x = 0 \in \pi_k\Hom_\cC(u,x')$.
 \end{itemize}
 Moreover, $\cU$ is \emph{strongly left filtering} if $\cU$ is a left filtering subcategory and the following condition holds:
 \begin{itemize}
     \item For every ingressive morphism $i \colon x \rightarrowtail y$ and every egressive morphism $q \colon y \twoheadrightarrow \overline{y}$ with fibre in $\cU$, there exist an egressive morphism $p \colon x \twoheadrightarrow x'$ with fibre in $\cU$, an ingressive morphism $j \colon x' \rightarrowtail y'$, a morphism $g \colon \overline{y} \to y'$ such that $gq$ is an egressive morphism with fibre in $\cU$, and a commutative square
     \[\begin{tikzcd}
        x\ar[r, tail, "i"]\ar[d, two heads, "p"'] & y\ar[d, two heads, "gq"] \\
        x'\ar[r, tail, "j"] & y'
     \end{tikzcd}\]
 \end{itemize}
\end{defn}

In \cref{thm:W-quotient}, we describe the quotient (in Waldhausen $\infty$-categories) of an exact $\infty$-category by a left filtering subcategory as a Dwyer--Kan localisation.
\cref{sec:localisation} then employs the general proof strategy to establish the folllowing.

\begin{theorem}\label[theorem]{thm:localisation-intro}
    Let $\cC$ be an exact $\infty$-category and let $\cU \subseteq \cC$ be an idempotent complete, extension-closed subcategory.
    If $\cU$ is strongly left filtering,
    then the induced sequence
    \[ \K(\cU) \to \K(\cC) \to \K(\cC/\cU) \]
    is a cofibre sequence of connective spectra, where $\cC/\cU$ denotes the cofibre in the $\infty$-category of Waldhausen $\infty$-categories.
\end{theorem}

\begin{rem}
 \ \begin{enumerate}
   \item \cref{thm:localisation} in the body of the text provides a generalisation of this statement which does not assume $\cU$ to be idempotent complete.
   \item The idea of identifying the cofibre as the K-theory of a Waldhausen category is not entirely new: it is implicit in the proof of the theorem of the heart \cite{barwick:heart}, and has also been employed by Sarazola \cite{sarazola:cotorsion} in the case of ordinary exact categories.
 \end{enumerate}
\end{rem}

In \cref{sec:applications}, we will apply \cref{thm:localisation-intro} to reprove the localisation theorem for stable $\infty$-categories, Quillen's localisatiom theorem for abelian categories, Schlichting's localisation theorem for s-filtering subcategories and Barwick's theorem of the heart \cite[Theorem~6.1]{barwick:heart}.
Even though it does not require any of the machinery developed in this paper, we also take the opportunity to explain how Sarazola's localisation theorem for cotorsion pairs follows from the localisation theorem for stable $\infty$-categories.

\addtocontents{toc}{\SkipTocEntry}
\hypertarget{target:conventions}{\subsection*{Conventions}}
\begin{enumerate}
    \item From this point forward, the term \emph{category} means \emph{($\infty$,1)-category}.
     We will say \emph{ordinary category} to emphasise that a category has mapping sets.
     In particular, $\catinf$ denotes the category of small categories.
     We will have no need to consider either the ordinary category or the 2-category of small ordinary categories.
    \item The category of anima/spaces/$\infty$-groupoids is denoted by $\Spc$, and its objects will be referred to as anima.
    \item Since we are primarily interested in additive categories, we adopt the following non-standard convention: $\Hom_\cC(x,y)$ will denote the mapping anima in an arbitrary category $\cC$; if $\cC$ is additive, $\hom_\cC(x,y)$ will denote the \emph{connective} spectra delooping the mapping anima in $\cC$.
    \item The fully faithful functor $\Spc \to \catinf$ admits both a right adjoint, the \emph{groupoid core} $\iota \colon \catinf \to \Spc$, and a left adjoint $\abs{-} \colon \catinf \to \Spc$, which we call \emph{realisation}.
    Note that the realisation $\abs{\cC}$ of a category is a Dwyer--Kan localisation at the collection of all morphisms in $\cC$.
    \item The category of small Waldhausen categories and exact functors is denoted by $\wald$, and the category of small exact categories and exact functors is denoted by $\exact$. We denote by $\idem{\cC}$ the idempotent completion of an exact category $\cC$, and by $\widem{\cC}$ the weak idempotent completion.
    Further explanations can be found in \cref{sec:gabriel-quillen}.
    \item For any Waldhausen category $\cC$, we denote by $S_\bullet(\cC)$ the evaluation of the Segal--Waldhausen construction on $\cC$.
    \item An \emph{additive invariant} on Waldhausen categories is an ``additive theory'' in Barwick's sense \cite[Definition~7.5]{barwick:algK}. 
    One way of phrasing this is to say that an additive invariant $H \colon \wald \to \cX$, where $\cX$ may be any category with finite products, is a functor which
    \begin{enumerate}
        \item preserves finite products;
        \item is group-like in the sense that the shear map $H(\cC \times \cC) \to H(\cC \times \cC)$ is an equivalence;
        \item has the property that the map $H(S_2(\cC)) \to H(\cC) \times H(\cC)$ induced by the face maps $d_0$ and $d_2$ of the $S_\bullet$-construction is an equivalence.
    \end{enumerate}
    See also \cite[Theorem~7.4]{barwick:algK}.
    
    Since $\wald$ is semiadditive, the first assumption implies that such a functor canonically refines to a functor $\wald \to \mathrm{CMon}(\Spc)$ valued in $\bbE_\infty$-monoids.
    The second assumption then asserts that this refinement factors over the full subcategory of group-like objects, so any additive invariant canonically refines to a functor valued in connective spectra by May's recognition principle.
    \item We denote the \emph{connective} algebraic K-theory functor by
    \[ \K \colon \wald \to \Sp_{\geq 0}. \]
    By \cite[Theorem~7.9 \& Corollary~7.14]{barwick:algK} this functor is the initial group-like additive functor with a natural transformation $\Sigma^\infty \iota \Rightarrow \K$.
    It satisfies $\Omega^\infty \K(-) \simeq \Omega\abs{\iota S_\bullet(-)}$, where $\abs{-} := \colim_{\Delta^\op} (-)$ denotes the realisation of a simplicial anima.
    In particular, the functor
    \[ \K_0 := \pi_0\K \colon \wald \to \Ab \]
    is not invariant under idempotent completion.
\end{enumerate}

\addtocontents{toc}{\SkipTocEntry}
\subsection*{Acknowledgements}
I am grateful to Clark Barwick for numerous very helpful discussions relating to the contents of this article.
Victor Saunier provided helpful comments on an earlier draft of this article.
Moreover, I would like to thank the anonymous referee for a very diligent and helpful report.

The author was supported by the CRC 1085 ``Higher Invariants" funded by the Deutsche Forschungsgemeinschaft (DFG).

\section{\texorpdfstring{The Gabriel--Quillen embedding}{The Gabriel-Quillen embedding}}\label{sec:gabriel-quillen}

In this section, we recall some necessary terminology and collect a number of results about the presentable version of the Gabriel--Quillen embedding for exact categories, details of which can be found in \cite{nw:gabriel-quillen}.

Recall that a \emph{Waldhausen category} is a category $\cC$ together with a wide subcategory $\cC_\ing$ of \emph{ingressive} morphisms, denoted by  feathered arrows $\rightarrowtail$, such that the following holds:
\begin{enumerate}
 \item $\cC$ is pointed;
 \item for every object $x \in \cC$, the unique map $0 \to x$ is ingressive;
 \item for every span $z \leftarrow x \rightarrowtail y$ in $\cC$, the pushout $z \sqcup_x y$ exists, and the structure morphism $z \to z \sqcup_x y$ is ingressive.
\end{enumerate}
If $i \colon x \rightarrowtail y$ is ingressive, we say that $\cofib(i) := y \sqcup_x 0$ is the \emph{quotient} of $y$ by $x$ and call the sequence
\[ x \overset{i}{\rightarrowtail} y \to \cofib(i) \]
an \emph{exact} sequence.
 
An \emph{exact} functor $f \colon (\cC,\cC_\ing) \to (\cD,\cD_\ing)$ of Waldhausen categories is a functor $f \colon \cC \to \cD$ between the underlying categories which preserves zero objects and ingressive morphisms and sends pushouts along ingressives to pushouts.
 
The collection of small Waldhausen categories and exact functors assemble into a category $\wald$.
The category $\exact$ of \emph{exact categories} is the full subcategory of $\wald$ spanned by those Waldhausen categories satisfying the following properties:
\begin{enumerate}
 \item $\cC$ is additive;
 \item the class of morphisms which are cofibres of ingressives is closed under pullback;
 \item every ingressive is a fibre of its cofibre.
\end{enumerate}
In that case, the collection of \emph{egressive morphisms}
\[ \cC_\egr := \{ p \in \cC \mid p \text{ is the cofibre of an ingressive morphism} \} \]
forms a subcategory of $\cC$.
Morphisms in $\cC_\egr$ will typically be denoted by double-headed arrows $\twoheadrightarrow$.
By \cite[Corollary~4.8.1 \& 4.9]{barwick:heart}, this description is equivalent to the more traditional definition of exact categories in terms of Quillen's axioms (see \cite{barwick:heart} or \cref{sec:diagram-lemmas} for additional details).

\begin{ex}\label[ex]{ex:prestable-exact}
    Every prestable category becomes an exact category by declaring all morphisms to be ingressive, see \cref{cor:prestable-exact}.
\end{ex}

\begin{defn}
    Let $\cC$ be an exact category.
    A full subcategory $\cU \subseteq \cC$ is \emph{extension-closed} if it is pointed and for every exact sequence $u \rightarrowtail y \twoheadrightarrow w$ with $u,w \in \cU$, the object $y$ also lies in $\cU$.
\end{defn}

An extension-closed subcategory $\cU$ of an exact category $\cC$ inherits an exact structure in which a morphism is ingressive precisely if it is ingressive in $\cC$ and its cofibre lies in $\cU$ \cite[Lemma~2.9]{sw:exact-stable}.
The Gabriel--Quillen embedding asserts that every exact category arises as an extension-closed subcategory of a prestable category.

\begin{defn}
    Let $\cC$ be a small Waldhausen category.
    Define the category $\psh_\lex(\cC)$ of \emph{left exact presheaves} to be the full subcategory of $\psh(\cC)$ spanned by those presheaves $X$ satisfying the following:
    \begin{enumerate}
        \item the presheaf $X \colon \cC^\op \to \Spc$ is \emph{additive}, ie it sends finite coproducts in $\cC$ to products in $\Spc$;
        \item if $i \colon x \rightarrowtail y$ is ingressive in $\cC$, then $X(\cofib(i)) \to X(y) \xrightarrow{i^*} X(x)$ is a fibre sequence.
    \end{enumerate}
\end{defn}

\begin{rem}
 If the underlying category is additive, the category of additive presheaves is equivalent to the category of additive $\Sp_{\geq 0}$-valued presheaves via composition with the functor $\Omega^\infty \colon \Sp_{\geq 0} \to \Spc$.

 Since $\Omega^\infty$ preserves and detects limits, we will regard $\psh_\lex(\cC)$ as a full subcategory of the category of $\Sp_{\geq 0}$-valued presheaves whenever $\cC$ is exact.
\end{rem}

Since every representable presheaf is left exact, we interpret (the $\Sp_{\geq 0}$-valued refinement of) the Yoneda embedding as a fully faithful functor
\[ \yo \colon \cC \to \psh_\lex(\cC). \]
Recall that a \emph{Grothendieck prestable category} is a prestable and presentable category in which filtered colimits commute with finite limits.

\begin{theorem}[{Gabriel--Quillen embedding, \cite[Theorem~2.10]{nw:gabriel-quillen}}]
 Let $\cC$ be a small Waldhausen category. The following are equivalent:
 \begin{enumerate}
  \item\label{it:gabriel-quillen1} $\cC$ is exact.
  \item \label{it:gabriel-quillen2} The category of left exact presheaves $\psh_\lex(\cC)$ is Grothendieck prestable and the Yoneda embedding $\yo \colon \cC \to \psh_\lex(\cC)$ identifies $\cC$ with an extension-closed subcategory.
  \item \label{it:gabriel-quillen3} There exist a small prestable category $\cD$ and a fully faithful functor $\cC \to \cD$ which identifies $\cC$ with an extension-closed subcategory of $\cD$.
 \end{enumerate}
\end{theorem}

If $\cC$ is exact, the Grothendieck prestable category $\psh_\lex(\cC)$ enjoys the following universal property.

\begin{prop}[{\cite[Theorem~2.7]{nw:gabriel-quillen}}]\label[prop]{prop:psh-lex-univ-prop}
    Let $\cC$ be a small exact category and let $\cD$ be a pointed, cocomplete category.
    Then restriction along the Yoneda embedding induces an equivalence
    \[ \yo^* \colon \Fun^\Ldec(\psh_\lex(\cC),\cD) \xrightarrow{\sim} \Fun^\exct(\cC,\cD), \]
    where the left hand side denotes the category of colimit-preserving functors and the right hand side denotes the category of exact functors with respect to the maximal Waldhausen structure on $\cD$.
\end{prop}

From the universal property of the presentable Gabriel--Quillen embedding, it is easy to deduce the existence of Klemenc's stable hull of an exact category.

\begin{cor}[{\cite[Theorem~1.2]{klemenc:stablehull}, \cite[Proposition~2.22]{nw:gabriel-quillen}}]
    The inclusion functor $\catst \to \exact$ admits a left adjoint
    \[ \derb \colon \exact \to \catst. \]
\end{cor}

Concretely, $\derb(\cC)$ is the Spanier--Whitehead stabilisation of the full prestable subcategory of $\psh_\lex(\cC)$ generated by the essential image of the Yoneda embedding.
In particular, $\cC$ is contained as an extension-closed subcategory in $\derb(\cC)$.

\begin{defn}
    An extension-closed subcategory $\cU$ of an exact category $\cC$ is \emph{left special} if the following holds: for every egressive morphism $p \colon x \twoheadrightarrow u$ with $x \in \cC$ and $u \in \cU$, there exists a morphism $f \colon v \to x$ with $v \in \cU$ such that the composite $pf \colon v \to u$ is an egressive morphism in $\cU$.

    An extension-closed subcategory $\cU \subseteq \cC$ is \emph{right special} if $\cU^\op$ is left special in $\cC^\op$.
\end{defn}

\begin{ex}
    An exact category $\cC$ is left special in $\psh_\lex(\cC)$ by \cite[Lemma~2.9]{nw:gabriel-quillen}.
\end{ex}

The next statement gives an explanation why this condition is particularly important.

\begin{prop}[{\cite[Proposition~2.18 \& Corollary~2.19]{nw:gabriel-quillen}}]\label[prop]{prop:special-recollement}
    Let $\cU$ be an extension-closed subcategory of an exact category $\cC$.
    The following are equivalent:
    \begin{enumerate}
        \item $\cU$ is left special in $\cC$;
        \item the restriction functor $j^* \colon \psh_\lex(\cC) \to \psh_\lex(\cU)$ preserves all colimits.
    \end{enumerate}
    If this is the case, the left adjoint $j_\natural \colon \psh_\lex(\cU) \to \psh_\lex(\cC)$ to the restriction functor is fully faithful and features in a recollement
	\[\begin{tikzcd}
	 \psh^\cU_\lex(\cC)\ar[r, hookrightarrow] & \psh_\lex(\cC)\ar[l, bend right=60]\ar[l, bend right=25, phantom, "\dashv" rotate=-90]\ar[l, bend left=25, phantom, "\dashv"' rotate =-90]\ar[l, bend left=60]\ar[r, "j^*"] & \psh_\lex(\cU)\ar[l, bend right=60, "j_\natural"'] \ar[l, bend right=35, phantom, "\dashv" rotate=-90]\ar[l, bend left=25, phantom, "\dashv" rotate=-90]\ar[l, bend left=60, "j_*"]
	\end{tikzcd}\]
	of Grothendieck prestable categories, where $\psh_\lex^\cU(\cC)$ denotes the full subcategory of left exact presheaves vanishing on $\cU$.
\end{prop}

For the convenience of the reader, we include the following observation about recollements of Grothendieck prestable categories.
See \cite[Appendix~A.2]{the-nine:2} and \cite{bg:recollements} for related discussions.

\begin{lem}\label[lem]{lem:recollement-cofibres}
    Let $\cY$ and $\cZ$ be Grothendieck prestable categories, and suppose that $p \colon \cY \to \cZ$ is a colimit-preserving functor which admits both a left adjoint $p^\Ldec$ and a right adjoint $p^\Rdec$.
    
    Then $p^\Ldec$ is fully faithful if and only if $p^\Rdec$ is fully faithful.
    If this is the case, there exists a recollement
    \[\begin{tikzcd}
	   \cX\ar[r, hookrightarrow, "i"] & \cY\ar[l, bend right=70, "i^\Ldec"']\ar[l, bend right=40, phantom, "\dashv" rotate=-90]\ar[l, bend left=25, phantom, "\dashv"' rotate =-90]\ar[l, bend left=60, "i^\Rdec"]\ar[r, "p"] & \cZ\ar[l, bend right=70, "p^\Ldec"'] \ar[l, bend right=40, phantom, "\dashv" rotate=-90]\ar[l, bend left=25, phantom, "\dashv" rotate=-90]\ar[l, bend left=60, "p^\Rdec"]
	\end{tikzcd}\]
    and both $\cX \xrightarrow{i} \cY \xrightarrow{p} \cZ$ and $\cZ \xrightarrow{p^\Ldec} \cY \xrightarrow{p^\Ldec} \cX$ are bifibre sequences in the category of Grothendieck prestable categories and colimit-preserving functors.
\end{lem}
\begin{proof}
    If either $p^\Ldec$ or $p^\Rdec$ is fully faithful, $p$ is a localisation functor, and it is well-known that any adjoint to a localisation functor is fully faithful.
    
    Assume now that $p$ is a localisation and let $i \colon \cX \to \cY$ be the kernel of $p$.
    For every $y \in \cY$, the object $\cofib(p^\Ldec p(y) \to y)$ satisfies $p(\cofib(p^\Ldec p(y) \to y)) \simeq \cofib(pp^\Ldec p(y) \to p(y)) \simeq 0$, and therefore lies in $\cX$.
    Moreover, we obtain for every $x \in \cX$ a fibre sequence
    \[ \Hom_\cB(\cofib(p^\Ldec p(y) \to y),i(x)) \to \Hom_\cB(y,i(x)) \to \Hom_\cB(p^\Ldec p(y),i(x)). \]
    By adjunction, $\Hom_\cB(p^\Ldec p(y),i(x)) \simeq \Hom_\cC(p(y),pi(x)) \simeq *$, so $\cofib(p^\Ldec p(y) \to y)$ is a left adjoint object to $y$ with respect to $i$.
    The existence of $i^\Ldec$ follows.
    The existence of the right adjoint $i^\Rdec$ is shown similarly by considering $\fib(y \to p^\Rdec p(y))$.

    Since $i$ admits both adjoints, it preserves all limits and colimits.
    It follows that $\cX$ is Grothendieck prestable.

    To see that $p$ is the cofibre of $i$, observe that we have an induced adjunction
    \[ - \circ p \colon \Fun^\Ldec(\cZ,\cD) \rightleftarrows \Fun^\Ldec(\cY,\cD) \cocolon - \circ p^\Ldec \]
    for every Grothendieck prestable category $\cD$.
    Since $p^\Ldec$ is fully faithful, this adjunction identifies $\Fun^\Ldec(\cZ,\cD)$ with the full subcategory of $\Fun^\Ldec(\cY,\cD)$ spanned by those colimit-preserving functors $F \colon \cY \to \cD$ which invert the counit transformation $p^\Ldec p \Rightarrow \id_\cB$.
    Now observe that such a functor $F$ induces for every $y \in \cY$ a cofibre sequence $F(p^\Ldec p(y)) \to F(y) \to F(i i^\Ldec(y))$.
    Since $\cD$ is prestable, the map $F(p^\Ldec p(y)) \to F(y)$ is an equivalence if and only if $F(ii^\Ldec(y)) \simeq 0$.
    As $i^\Ldec$ is essentially surjective, this means that $F$ inverts the counit transformation if and only if $F$ vanishes on $\cX$.

    Consider now the sequence $\cZ \xrightarrow{p^\Ldec} \cY \xrightarrow{i^\Ldec} \cX$.
    Clearly, the composite of this sequence is trivial.
    Since $\cY$ is prestable, the cofibre sequence $p^\Ldec p(y) \to y \to ii^\Ldec(y)$ is also a fibre sequence.
    Hence any $y \in \cY$ with $i^\Ldec(y) \simeq 0$ lies in the essential image of $p^\Ldec$, showing that $p^\Ldec$ is the kernel of $i^\Ldec$.
    As before, we have an adjunction
    \[ - \circ i \colon \Fun^\Ldec(\cY,\cD) \rightleftarrows \Fun^\Ldec(\cX,\cD) \cocolon - \circ i^\Ldec \]
    which identifies $\Fun^\Ldec(\cX,\cD)$ with the full subcategory spanned by those colimit-preserving functors $F \colon \cY \to \cD$ which invert the unit transformation $\id_\cX \Rightarrow i i^\Ldec$.
    Considering the induced cofibre sequence $F(p^\Ldec p(y)) \to F(y) \to F(i i^\Ldec(y))$, which is also a fibre sequence by prestability of $\cD$, we see that $F$ inverts the unit transformation if and only if $F$ vanishes on the essential image of $p^\Ldec$.
\end{proof}

We will need to consider both idempotent completions and weak idempotent completions of exact categories, so let us collect some basic facts about these.

\begin{defn}
    An additive category $\cA$ is \emph{weakly idempotent complete} if for every retraction diagram $x \xrightarrow{i} y \xrightarrow{r} x$ in $\cA$ there exists an object $y'\in \cA$ and a morphism $j \colon y'\to y$ such that $i + j \colon x \oplus y' \to y$ is an equivalence.
\end{defn}

\begin{ex}
    Every stable category is weakly idempotent complete.
\end{ex}

We call an exact category (weakly) idempotent complete if its underlying additive category is (weakly) idempotent complete.

\begin{prop}[{\cite[Lemma~4.10]{sw:exact-stable}, \cite[Corollary~2.26]{nw:gabriel-quillen}}]
    Let $\cC$ be an exact category.
    \begin{enumerate}
        \item There exists a weakly idempotent complete exact category $\widem{\cC}$ which contains $\cC$ as an extension-closed subcategory such that the restriction functor
        \[ \Funex{}(\widem{\cC},\cD) \to \Funex{}(\cC,\cD) \]
        is an equivalence for every weakly idempotent complete exact category $\cD$.
        \item The idempotent completion $\idem{\cC}$ of $\cC$ carries an exact structure such that $\cC$ is an extension-closed subcategory and the restriction functor
        \[ \Funex{}(\idem{\cC},\cD) \to \Funex{}(\cC,\cD) \]
        is an equivalence for every idempotent complete exact category $\cD$.
    \end{enumerate}
    In either case, the exact sequences in the completion are precisely the retracts of exact sequences in $\cC$.
\end{prop}

In fact, the weak idempotent completion can be described relatively explicitly in terms of the idempotent completion.
Let us first recall some basic facts about the Grothendieck group $\K_0$.

\begin{rem}\label[rem]{rem:equality-in-k0}
    Let $\cC$ be an exact category.
    Recall that two elements $[x] - [y]$ and $[x'] - [y']$ in $\K_0(\cC)$ are equal if and only if there exist $a,b,s \in \cC$ and exact sequences $a \rightarrowtail x \oplus y'\oplus s \twoheadrightarrow b$ and $a \rightarrowtail x'\oplus y \oplus s \twoheadrightarrow b$.
    
    The proof of this statement in \cite[Lemma~3.4]{kw:shortening} applies verbatim to arbitrary exact categories.
\end{rem}

\begin{defn}\label[defn]{def:dense}
    An extension-closed subcategory $\cC_0 \subseteq \cC$ of an exact category is \emph{dense} if for every object $x \in \cC$ there exists some $x' \in \cC$ such that $x \oplus x' \in \cC_0$.
\end{defn}

\begin{ex}
    The inclusion of an exact category into its (weak) idempotent completion is dense.
\end{ex}

In most of the literature on algebraic K-theory, dense subcategories go by the name ``cofinal''.
The term ``dense'' has been adopted from \cite[Section~1.3]{the-nine:2} to resolve the obvious naming conflict with cofinal functors.
The next lemma asserts that \cref{def:dense} agrees with the definition of denseness in {\it loc.\ cit.}

\begin{lem}
    Let $\cC$ be a weakly idempotent complete exact category and let $\cC_0$ be an extension-closed subcategory of $\cC$.
    Then $\cC_0$ is dense in $\cC$ if and only if every object in $\cC$ is a retract of an object in $\cC_0$.
\end{lem}
\begin{proof}
    If $\cC_0$ is dense in $\cC$, every object in $\cC$ is in particular a retract of an object in $\cC_0$.
    If $x \in \cC$ is a retract of $y \in \cC_0$, weak idempotent completeness implies that there exists an object $x'\in \cC$ such that $x \oplus x' \simeq y$, showing the other direction.
\end{proof}

\begin{rem}\label[rem]{rem:k0-of-idem}
    \ \begin{enumerate}
        \item It is well-known that dense inclusions induce injections on $\K_0$.
        \item Every element in $\K_0(\idem{\cC})$ can be represented as a difference $[\widehat{x}] - [x]$ with $x \in \cC$: this can be accomplished by picking an arbitrary representative $[\widehat{y}] - [\widehat{z}]$, choosing $z' \in \idem{\cC}$ with $\widehat{z} \oplus z' \in \cC$, and observing that
        \[ [\widehat{y}] - [\widehat{z}] = [\widehat{y}] - [\widehat{z}] + [z'] - [z'] = [\widehat{y} \oplus \widehat{z}] - [\widehat{z} \oplus z']. \]
        \item If $\cC_0$ is dense in $\cC$, then two elements $[x] - [y]$ and $[x'] - [y']$ in $\K_0(\cC)$ are equal if and only if there exist $a,b \in \cC_0$ and $s \in \cC$ together with exact sequences $a \rightarrowtail x \oplus y'\oplus s \twoheadrightarrow b$ and $a \rightarrowtail x'\oplus y \oplus s \twoheadrightarrow b$ in $\cC$: starting with $a,b,s \in \cC$ and appropriate exact sequences witnessing the equality, choose $a',b'\in \cC$ with $a \oplus a' \in \cC_0$ and $b \oplus b' \in \cC_0$ and replace the object $a$ by $a \oplus a'$, the object $b$ by $b \oplus b'$ and the object $s$ by $s \oplus a'\oplus b'$.

    \end{enumerate}
\end{rem}

One characterisation of the weak idempotent completion is now the following.

\begin{lem}\label[lem]{lem:widem-explicitly}
    Let $\cC$ be an exact category.
    Then the following full subcategories of the idempotent completion $\idem{\cC}$ coincide:
    \begin{enumerate}
        \item\label{lem:widem-explicitly-1} $\{ \widehat{x} \in \idem{\cC} \mid [\widehat{x}] \in \K_0(\cC) \subseteq \K_0(\idem{\cC}) \}$;
        \item\label{lem:widem-explicitly-2} $\{ \widehat{x} \in \idem{\cC} \mid \exists\,x\in \cC \colon \widehat{x} \oplus x \in \cC \}$.
    \end{enumerate}
    Call this subcategory $\cC_0$.
    Then $\cC_0$ is closed under extensions and the inclusion $\cC \to \cC_0$ exhibits $\cC_0$ as a weak idempotent completion of $\cC$.
\end{lem}
\begin{proof}
    Let us first show that the two explicit subcategories of $\idem{\cC}$ agree.
    If $\widehat{x} \in \idem{\cC}$ and $x \in \cC$ are objects with $\widehat{x} \oplus x \in \cC$, then $[\widehat{x}] = [\widehat{x} \oplus x] - [x] \in \K_0(\cC)$.
    Conversely, if $[\widehat{x}] = [x] - [y] \in \K_0(\idem{\cC})$ with $x,y \in \cC$, then there exist by \cref{rem:k0-of-idem} objects $a,b \in \cC$ and $\widehat{s} \in \idem{\cC}$ together with exact sequences $a \rightarrow \widehat{x} \oplus y \oplus \widehat{s} \twoheadrightarrow b$ and $a \rightarrow x \oplus \widehat{s} \twoheadrightarrow b$.
    Then $x \oplus \widehat{s} \in \cC$, so $x \oplus y \oplus \widehat{s}$ is an object in $\cC$ with $\widehat{x} \oplus x \oplus y \oplus \widehat{s} \in \cC$.
    
    It is obvious that $\cC_0$ is an extension-closed subcategory of $\idem{\cC}$, and the first description of $\cC_0$ implies directly that $\cC_0$ is weakly idempotent complete.

    It now suffices to show that $\cC_0$ is a weak idempotent completion of $\cC$.
    Given a weakly idempotent complete exact category $\cD$, consider the commutative square
    \[\begin{tikzcd}
        \Funex{}(\cC_0,\cD)\ar[r]\ar[d] & \Funex{}(\cC,\cD)\ar[d] \\
        \Funex{}(\cC_0,\idem{\cD})\ar[r] & \Funex{}(\cC,\idem{\cD})
    \end{tikzcd}\]
    Since $\idem{(\cC_0)} \simeq \idem{\cC}$, the universal property of the idempotent completion implies that the lower horizontal arrow is an equivalence.
    Suppose that $f \colon \cC_0 \to \idem{\cD}$ is an exact functor which maps $\cC$ to $\cD$.
    Using the second description of $\cC_0$, the weak idempotent completeness of $\cD$ implies that $f(\cC_0) \subseteq \cD$ as well, which shows that the above square is a pullback.
\end{proof}

\section{Quotients of exact categories}\label{sec:quotients}

Since $\exact$ is a full subcategory of the category of Waldhausen categories, it is natural to ask to which extent the inclusion functor preserves cofibres.
To distinguish the two types of cofibres, we introduce the following notation.

\begin{defn}
 Let $\cC$ be an exact category and let $\cU$ be an extension-closed subcategory.
 \begin{enumerate}
     \item Denote by $\quotw{\cC}{\cU}$ the \emph{Waldhausen quotient} of $\cC$ by $\cU$, ie the cofibre of the inclusion functor $\cU \to \cC$ in $\wald$;
     \item Denote by $\quotex{\cC}{\cU}$ the \emph{exact quotient} of $\cC$ by $\cU$, ie the cofibre of the inclusion functor $\cU \to \cC$ in $\exact$.
 \end{enumerate}
\end{defn}

Note that the universal property of $\quotw{\cC}{\cU}$ implies that there is an exact comparison functor $\quotw{\cC}{\cU} \to \quotex{\cC}{\cU}$.
The following statement will only be used in \cref{sec:applications}.

\begin{prop}\label[prop]{prop:special-quotient}
 Let $\cC$ be an exact category and let $\cU$ be an extension-closed subcategory.
 If $\cU$ is left special, then the comparison map $\quotw{\cC}{\cU} \to \quotex{\cC}{\cU}$ is an equivalence.
\end{prop}
\begin{proof}
 By comparing universal properties, we have $\psh_\lex(\quotw{\cC}{\cU}) \simeq \psh_\lex^\cU(\cC)$.
 Combining \cref{prop:special-recollement} and \cref{lem:recollement-cofibres}, we obtain a cofibre sequence
 \[ \psh_\lex(\cU) \xrightarrow{j_\natural} \psh_\lex(\cC) \to \psh_\lex^\cU(\cC) \]
 of Grothendieck prestable categories. 
 Since $\psh_\lex(\cU)$ is generated by $\cU$ under colimits, the universal property of $\psh_\lex$ from \cref{prop:psh-lex-univ-prop} implies that $\psh_\lex(\quotex{\cC}{\cU})$
 is also a cofibre of $j_\natural$ in Grothendieck prestable categories.
 So the canonical functor $\psh_\lex(\quotw{\cC}{\cU}) \to \psh_\lex(\quotex{\cC}{\cU})$ is an equivalence.
 Since the Yoneda embedding is fully faithful, it follows that $\quotw{\cC}{\cU} \to \quotex{\cC}{\cU}$ is fully faithful.

 It is now enough to show that the quotient functor $p \colon \cC \to \quotex{\cC}{\cU}$ is essentially surjective.
 Note that $\quotex{\cC}{\cU}$ is necessarily generated by the essential image of $p$ under extensions, so we have to show that the essential image of $p$ is closed under extensions.
 Consider an exact sequence $p(x) \to y \to p(z)$ in $\quotex{\cC}{\cU}$.
 Since the inclusion functor $\psh_\lex^\cU(\cC) \to \psh_\lex(\cC)$ preserves cofibres, we obtain a cofibre sequence $p^*\yo(p(x)) \to p^*\yo(y) \to p^*\yo(p(z))$ in $\psh_\lex(\cC)$.
 Pulling back the egressive map $p^*\yo(y) \to p^*\yo(p(z))$ along the unit morphism $\yo(z) \to p^*\yo(p(z))$, we obtain an exact sequence $p^*\yo(p(x)) \to Y \to \yo(z)$.
 As $\cC$ is left special in $\psh_\lex(\cC)$, we find a pushout
 \[\begin{tikzcd}
     \yo(x')\ar[r, "\yo(i)"]\ar[d] & \yo(y')\ar[d] \\
     p^*\yo(p(x))\ar[r] & Y
 \end{tikzcd}\]
 in $\psh_\lex(\cC)$, where $i \colon x' \to y'$ is an ingressive morphism in $\cC$.
 Applying the left adjoint localisation functor $p_\natural \colon \psh_\lex(\cC) \to \psh_\lex(\quotex{\cC}{\cU})$ to this pushout, we obtain a pushout
 \[\begin{tikzcd}
     \yo(p(x'))\ar[r]\ar[d] & \yo(p(y'))\ar[d] \\
     \yo(p(x))\ar[r] & p_\natural Y
 \end{tikzcd}\]
 Consequently, $p_\natural Y \simeq \yo(p(x) \cup_{p(x')} p(y'))$.
 By construction, $p_\natural Y \simeq y$, so the essential image of $p$ is closed under extensions.
\end{proof}

\begin{rem}
    By considering opposite categories, \cref{prop:special-quotient} also shows that the coWaldhausen quotient of an exact category by a right special subcategory coincides with the exact quotient.
\end{rem}

Later in this section, we will see examples of extension-closed subcategories for which the Waldhausen cofibre differs from the exact cofibre.
Exhibiting such examples will become easier once we have obtained more explicit descriptions of some Waldhausen quotients.
The main point of the upcoming discussion is that the Waldhausen quotient is a Dwyer--Kan localisation of the underlying category if we take the quotient by a filtering subcategory.

We begin with an elementary observation which allows us to decide when a Dwyer--Kan localisation of a Waldhausen category inherits the structure of a Waldhausen category.

\begin{lem}\label[lem]{lem:localisation-Waldhausen}
    Let $\cC$ be a Waldhausen category and let $S \subseteq \cC$ be a collection of morphisms in $\cC$.
    Suppose that the Dwyer--Kan localisation $\ell \colon \cC \to \cC[S^{-1}]$ has the following properties:
    \begin{enumerate}
        \item $\ell$ preserves pushouts along ingressive morphisms;
        \item every morphism $\ell(x) \to y'$ in $\cC[S^{-1}]$ factors as a composition
        \[ \ell(x) \xrightarrow{\ell(f)} \ell(y) \xrightarrow{\sim} y' \]
        of a morphism in the image of $\ell$ followed by an equivalence. 
    \end{enumerate}
    Then $\cC[S^{-1}]$ is a Waldhausen category with respect to the subcategory
    \[ \cC[S^{-1}]_\ing := \{ f \in \cC[S^{-1}] \mid f \text{ is equivalent to the image of a morphism in } \cC_\ing \}. \]
    With respect to this Waldhausen structure, $\ell$ is exact, and a functor $F \colon \cC[S^{-1}] \to \cD$ to another Waldhausen category $\cD$ is exact if and only if $F \circ \ell$ is exact.
\end{lem}
\begin{proof}
    Note that every Dwyer--Kan localisation preserves zero objects.
    By definition, $\cC[S^{-1}]_\ing$ contains the groupoid core of $\cC[S^{-1}]$.
    
    Let $z' \leftarrow x' \rightarrowtail y'$ be a span in $\cC[S^{-1}]$ with the right leg ingressive.
    By definition, the right leg is equivalent to the image of an ingressive morphism in $\cC$, and the assumptions allow us to lift the entire span to $\cC$.
    Since $\cC$ admits pushouts along ingressives and $\ell$ preserves these, it follows that the original span has a pushout in $\cC[S^{-1}]$, and that the induced map $z' \to z' \cup_{x'} y'$ is ingressive.
    Hence $\cC[S^{-1}]$ is a Waldhausen category and $\ell$ is an exact functor.

    It follows from the definition of the Waldhausen structure on $\cC[S^{-1}]$ that a functor $F \colon \cC[S^{-1}] \to \cD$ is exact if and only if $F \circ \ell$ is exact.
\end{proof}

We now turn our attention to filtering subcategories of exact categories.
It will be convenient to rephrase \cref{def:filtering-intro} slightly.

\begin{defn}
  Let $S$ be a collection of morphisms in a category $\cC$.
  \begin{enumerate}
   \item A presheaf $X \colon \cC^\op \to \Ab$ is \emph{weakly $S$-effaceable} if for each $x \in \cC$ and each $\alpha \in X(x)$ there exists some $s \colon x' \to x$ in $S$ such that $s^*\alpha = 0 \in X(x')$.
   \item A presheaf $X \colon \cC^\op \to \Sp_{\geq 0}$ is \emph{weakly $S$-effaceable} if $\pi_kX$ is weakly $S$-effaceable for all $k \geq 0$.
  \end{enumerate}
\end{defn}

\begin{defn}\label[defn]{def:filtering}
 Let $\cC$ be an exact category and let $\cU \subseteq \cC$ be an extension-closed subcategory.
  \begin{enumerate}
   \item Denote by $\cC_\ing^\cU$ the subcategory of ingressive morphisms whose cofibre lies in $\cU$ and denote by $\cC_\egr^\cU$ the subcategory of egressives whose fibre lies in $\cU$.
   \item The subcategory $\cU \subseteq \cC$ is \emph{right filtering} if the represented functor
   \[ \hom_\cC(-,u) \colon \cC^\op \to \Sp_{\geq 0} \]
   is weakly $\cC_\ing^\cU$-effaceable for every $u \in \cU$.
   \item The subcategory $\cU \subseteq \cC$ is \emph{left filtering} if the corepresented functor
   \[ \hom_\cC(u,-) \colon \cC \to \Sp_{\geq 0} \]
   is weakly $\cC_\egr^\cU$-effaceable for every $u \in \cU$.
 \end{enumerate}
\end{defn}

Note that $\cU \subseteq \cC$ is left filtering if and only if $\cU^\op \subseteq \cC^\op$ is right filtering, and that this definition is equivalent to the condition formulated in \cref{def:filtering-intro}.

\begin{ex}\label[ex]{ex:filtering-1-cats}
 If $\cC$ is an ordinary exact category and $\cU \subseteq \cC$ is an extension-closed subcategory, the right filtering condition reduces to the following requirement: for every morphism $f \colon x \to u$ with $x \in \cC$ and $u \in \cU$, there exists an admissible epimorphism $p \colon x \twoheadrightarrow v$ and a morphism $f' \colon v \to u$ such that $f = f'p$.
 This provides the following examples:
 \begin{enumerate}
  \item\label{ex:filtering-1-cats-1} If $\cA$ is an abelian category and $\cU \subseteq \cA$ is a Serre subcategory, the canonical epi-mono-factorisation of morphisms shows that $\cU$ is both left and right filtering in $\cA$.
  \item\label{ex:filtering-1-cats-2} If $\cA$ is an ordinary additive category and $\cU \subseteq \cA$ is a full additive subcategory, the left (or right) filtering condition is a non-symmetric version of the notion of a Karoubi-filtering subcategory introduced in \cite[Section~5]{pw:khomology}---see \cite[Section~5]{kasprowski:fdc} for a simplified set of axioms.
  \item\label{ex:filtering-1-cats-3} If an extension-closed subcategory of an ordinary exact category is left filtering in the sense of \cite[Definition~1.3]{schlichting:delooping}, then it is also left filtering in the sense of \cref{def:filtering}.
 \end{enumerate}
\end{ex}

\begin{ex}\label[ex]{ex:right-exact-filtering}
 If $\cC$ is a stable category and $\cU \subseteq \cC$ is a full stable subcategory, then $\cU$ is left filtering in $\cC$.
\end{ex}

Further examples of filtering subcategories can be obtained from the following lemma.

\begin{lem}\label[lem]{lem:adjoint-filtering}
 Let $\cC$ be an exact category and let $L \colon \cC \to \cD$ be a Bousfield localisation of the underlying category.
 Assume one of the following:
 \begin{enumerate}
  \item\label{lem:adjoint-filtering-1} Every morphism in $\cC$ is egressive (ie $\cC$ is the opposite of a prestable category).
  \item\label{lem:adjoint-filtering-2} $\cC$ is an ordinary category and the unit morphism $x \to Lx$ is egressive for every $x \in \cC$.
 \end{enumerate}
 Then $\kernel(L) \subseteq \cC$ is left filtering.
\end{lem}
\begin{proof}
 The kernel of $L$ is clearly an extension-closed subcategory.
 Let $a \in \kernel(L)$ and $x \in \cC$.
 In either case, the naturality square
 \[\begin{tikzcd}
  a\ar[r]\ar[d] & x\ar[d] \\
  0 \simeq La\ar[r] & Lx	
 \end{tikzcd}\]
 associated to an arbitrary morphism $a \to x$ shows that $\pi_0\Hom_\cC(a,-)$ is weakly $\cC_\egr^{\kernel(L)}$-effaceable.
 In case \eqref{lem:adjoint-filtering-2}, this finishes the proof.
 In case \eqref{lem:adjoint-filtering-1}, we have natural equivalences
 \[ \pi_n\hom_\cC(a,-) \simeq \pi_0\hom_\cC(a, \Omega^n -), \]
 which reduces weak $\cC_\egr^{\kernel(L)}$-effaceablity to the $\pi_0$-statement.
\end{proof}

\begin{ex}\label[ex]{ex:filtering-subcats}
 \cref{lem:adjoint-filtering} applies to show that the following are examples of left filtering subcategories:
 \begin{enumerate}
    \item\label{ex:filtering-subcats-1} Let $R$ be a noetherian ring and denote by $\Mod_R^\fg$ the abelian category of finitely generated $R$-modules.
 	The full subcategory $\Mod^{\fg,\tf}_R$ of finitely generated torsionfree modules is a Bousfield localisation of $\Mod_R^\fg$: the left adjoint $(-)^\tf$ sends a module $M$ to its maximal torsionfree quotient.
 	Since the unit map $M \to M^\tf$ is an epimorphism, the full subcategory $\Mod_R^{\fg,\tors}$ of torsion modules is left filtering.

    While this is a special case of \cref{ex:filtering-1-cats}~\eqref{ex:filtering-1-cats-1}, we will use this observation in \cref{ex:w-quotients} below to conclude that $\Mod^{\fg,\tf}_R$ is the Waldhausen quotient of $\Mod_R^\fg$ by $\Mod_R^{\fg,\tors}$.
    %
   \item\label{ex:filtering-subcats-2} Let $\cC$ be a stable category with a t-structure $(\cC_{\geq 0},\cC_{\leq 0})$.
   Then the coconnective part $\cC_{\leq 0}$ carries an exact structure in which ingressives are the $\pi_0$-monomorphisms and egressives are all morphisms (because $\cC_{\leq 0}^\op$ is prestable).
   The Bousfield localisation
   \[ \tau_{\leq -1} \colon \cC_{\leq 0} \rightleftarrows \cC_{\leq -1} \cocolon \inc \]
   shows that $\kernel(\tau_{\leq -1}) \simeq \cC^\heartsuit$ is left filtering in $\cC_{\leq 0}$.
 \end{enumerate}
\end{ex}

We have the following characterisation of filtering subcategories.

\begin{lem}\label[lem]{lem:filtering-filtered}
 Let $\cC$ be an exact category and let $\cU \subseteq \cC$ be an extension-closed subcategory.
 For $x \in \cC$, denote by $\cC_\egr^\cU(x) \subseteq \cC_{x/}$ the full subcategory on morphisms in $\cC_\egr^\cU$.
 
 Then the following are equivalent:
 \begin{enumerate}
  \item\label{lem:filtering-filtered-1} $\cU$ is left filtering in $\cC$.
  \item\label{lem:filtering-filtered-2} $\cC_\egr^\cU(x)$ is filtered for every $x \in \cC$ and we have
   \[ \colim_{y \twoheadrightarrow y' \in \cC_\egr^\cU(y)} \hom_\cC(u,y') \simeq 0 \]
   for all $u \in \cU$ and $y \in \cC$.
  \item\label{lem:filtering-filtered-3} $\cC_\egr^\cU(x)$ is filtered for every $x \in \cC$ and the functor
   \[ \cC^\op \to \Sp_{\geq 0},\quad x \mapsto \colim_{y \twoheadrightarrow y' \in \cC_\egr^\cU(y)} \hom_\cC(x,y') \]
   inverts all morphisms in $\cC_\egr^\cU$ for every $y \in \cC$.
   \item\label{lem:filtering-filtered-4} $\cC_\egr^\cU(x)$ is filtered for every $x \in \cC$ and the functor
   \[ \cC^\op \to \Spc,\quad x \mapsto \colim_{y \twoheadrightarrow y' \in \cC_\egr^\cU(y)} \Hom_\cC(x,y') \]
   inverts all morphisms in $\cC_\egr^\cU$ for every $y \in \cC$.
 \end{enumerate}
\end{lem}
\begin{proof}
 Conditions~\eqref{lem:filtering-filtered-3} and \eqref{lem:filtering-filtered-4} are equivalent since $\Omega^\infty \colon \Sp_{\geq 0} \to \Spc$ is conservative and preserves filtered colimits.

 We show that \eqref{lem:filtering-filtered-2} and \eqref{lem:filtering-filtered-3} are equivalent.
 Since filtered colimits preserve fibre sequences of connective spectra, we have for every cofibre sequence $u \rightarrowtail x \twoheadrightarrow x'$ a fibre sequence
 \[ \colim_{y \twoheadrightarrow y' \in \cC_\egr^\cU(y)} \hom_\cC(x',y') \to \colim_{y \twoheadrightarrow y' \in \cC_\egr^\cU(y)} \hom_\cC(x,y') \to \colim_{y \twoheadrightarrow y' \in \cC_\egr^\cU(y)} \hom_\cC(u,y'). \]
 Assertion~\eqref{lem:filtering-filtered-2} says that the rightmost term vanishes, so assertion~\eqref{lem:filtering-filtered-3} follows.
 
 Conversely, apply assertion~\eqref{lem:filtering-filtered-3} to the egressive morphism $p \colon u \twoheadrightarrow 0$ to see that
 \[ 0 \simeq \colim_{y \twoheadrightarrow y' \in \cC_\egr^\cU(y)} \hom_\cC(0,y') \xrightarrow{p^*} \colim_{y \twoheadrightarrow y' \in \cC_\egr^\cU(y)} \hom_\cC(u,y') \]
 is an equivalence.

 Next, we show $\eqref{lem:filtering-filtered-2} \Longrightarrow \eqref{lem:filtering-filtered-1}$.
 Let $u \in \cU$, $y \in \cC$ and $\alpha \in \pi_k\hom_\cC(u,y)$. 
 Then
 \[ \pi_k \Big(\colim_{y \twoheadrightarrow y' \in \cC_\egr^\cU(y)} \hom_\cC(u,y') \Big) \simeq \colim_{y \twoheadrightarrow y' \in \cC_\egr^\cU(y)} \pi_k\hom_\cC(u,y') \simeq 0 \]
 because the indexing category is filtered.
 Hence the map
 \[ \pi_k\hom_\cC(u,y) \to \pi_k\Big(\colim_{y \twoheadrightarrow y' \in \cC_\egr^\cU(y)} \hom_\cC(u,y')\Big) \simeq 0 \]
 necessarily trivialises $\alpha$.
 Using another time that the indexing category is filtered, there must exist some egressive $p \colon y \twoheadrightarrow y'$ with fibre in $\cU$ such that $p_*\alpha = 0$.
 
 For the converse implication $\eqref{lem:filtering-filtered-1} \Longrightarrow \eqref{lem:filtering-filtered-2}$, we use the local characterisation of filteredness \cite[\href{https://kerodon.net/tag/02PS}{Theorem~02PS}]{kerodon}.
 Let $p_1 \colon x \twoheadrightarrow x_1$ and $p_2 \colon x \twoheadrightarrow x_2$ be egressives with fibre in $\cU$.
 Using the filtering assumption, the induced map $\fib(p_1) \to x_2$ factors over an ingressive morphism $u \rightarrowtail x_2$, and we obtain a commutative square
 \[\begin{tikzcd}
     x\ar[r, two heads, "p_1"]\ar[d, "p_2"'] & x_1\ar[d] \\
     x_2\ar[r, two heads, "q"] & y
 \end{tikzcd}\]
 in which $q$ also lies in $\cC_\egr^\cU$. This defines morphisms $p_1 \to qp_2 \leftarrow p_2$ in $\cC_\egr^\cU(x)$.
 
 Let now $\phi \in \Hom_{\cC_\egr^\cU(x)}(p_1,p_2)$ be an arbitrary basepoint and consider an element $\alpha \in \pi_k\Hom_{\cC_\egr^\cU(x)}(p_1,p_2)$.
 We have to show that there exists a morphism $q \colon p_2 \to p_3$ such that $q_*\alpha$ is trivial.
 Since $\cC_\egr^\cU(x)$ is a full subcategory of $\cC_{x/}$, we have an equivalence
 \[ \Hom_{\cC_\egr^\cU(x)}(p_1,p_2) \simeq \{ p_2 \} \mathop{\times}\limits_{\Hom_\cC(x,x_2)} \Hom_\cC(x_1,x_2), \]
 the pullback being taken along the restriction map $p_1^*$.
 Denote by $f$ the image of $\phi$ in $\Hom_\cC(x_1,x_2)$.
 The given element $\alpha$ lifts to an element
 \[ \alpha' \in \pi_{k}\Big( \{fp_1\} \times_{\Hom_\cC(x,x_2)} \Hom_\cC(x_1,x_2), \overline{f} \Big), \]
 where $\overline{f}$ denotes the induced basepoint of the fibre.
 Since $\cC$ is additive, there is a natural isomorphism
 \begin{align*}
     \pi_{k}\Big( \{fp_1\} \times_{\Hom_\cC(x,x_2)} \Hom_\cC(x_1,x_2), \overline{f} \Big)
     &\cong \pi_{k}(\Hom_\cC(\fib(p_1),x_2), 0) \\
     &\cong \pi_{k}\hom_\cC(\fib(p_1),x_2),
 \end{align*}
 so $\alpha'$ becomes trivial after postcomposition with a morphism $q \colon x_2 \twoheadrightarrow x_3$ in $\cC_\egr^\cU$.
 A quick diagram chase shows that the morphism $q \colon p_2 \to qp_2$ in $\cC_\egr^\cU$ trivialises $\alpha$ as required.
 
 Knowing that $\cC_\egr^\cU(x)$ is filtered for all $x$, the left filtering condition immediately implies that
 \[ \colim_{y \twoheadrightarrow y' \in \cC_\egr^\cU(y)} \Hom_\cC(u,y') \simeq 0 \]
 for all $u \in \cU$ and $y \in \cC$.
\end{proof}

We can now give a concise description of Waldhausen quotients by left filtering subcategories.

\begin{theorem}\label[theorem]{thm:W-quotient}
 Let $\cC$ be an exact category and let $\cU \subseteq \cC$ be a left filtering subcategory. Then the following holds:
 \begin{enumerate}
  \item\label{it:W-quotient1} The Dwyer--Kan localisation $\ell \colon \cC \to \cC[(\cC_\egr^\cU)^{-1}]$ is a cofibre of the inclusion $\cU \subseteq \cC$ in $\wald$ if we equip the target with the Waldhausen structure in which a morphism is ingressive if and only if it is equivalent to the image of an ingressive morphism in $\cC$.
  \item\label{it:W-quotient2} A morphism $f \in \cC$ satisfies $\ell(f) \simeq 0$ if and only if $f$ factors through some object of $\cU$.
  In particular, an object $x \in \cC$ satisfies $\ell(x) \simeq 0$ if and only if $x$ is a retract of an object in $\cU$.
  \item\label{it:W-quotient3} A morphism $f \colon x \to y$ in $\cC$ becomes invertible in $\cC[(\cC_\egr^\cU)^{-1}]$ if and only if there exists a morphism $g \colon y \to y'$ with the following properties:
  \begin{enumerate}
   \item $gf \in \cC_\egr^\cU$;
   \item $g$ is an egressive morphism in $\idem{\cC}$ whose fibre lies in $\idem{\cU}$.
  \end{enumerate}
  
  \item\label{it:W-quotient4} If $\cC$ is an ordinary category, then so is $\cC[(\cC_\egr^\cU)^{-1}]$. \end{enumerate}
\end{theorem}
\begin{proof}
	Since it seems more natural to work with the Yoneda embedding, let us formulate the proof in the case that $\cU$ is a right filtering subcategory. The subcategory $\cC_\ing^\cU$ of ingressives with cofibre in $\cU$ is closed under direct sums, so the Dwyer--Kan localisation $\ell \colon \cC \to \cD := \cC[(\cC_\ing^\cU)^{-1}]$ is also the localisation among additive categories.
	Therefore, $\ell$ induces a Bousfield localisation
	\[ \ell_! \colon \psh_\Sigma(\cC) \rightleftarrows \psh_\Sigma(\cD) \cocolon \ell^* \]
	of prestable categories.
	By considering the induced commutative square
	\[\begin{tikzcd}
	 \cC\ar[r, "\ell"]\ar[d, "\yo"'] & \cD\ar[d, "\yo"] \\
	 \psh_\Sigma(\cC)\ar[r, "\ell_!"] & \psh_\Sigma(\cD)
    \end{tikzcd}\]
	and using that $\yo$ is fully faithful, we observe that $\ell$ preserves finite limits if $\ell_!$ does so.
		
  For $x \in \cC$, consider the functor
  \[ F \colon \cC \to \Spc,\quad y \mapsto \colim_{x' \rightarrowtail x \in \cC_\ing^\cU(x)^\op} \Hom_\cC(x',y). \]
  The dual of \cref{lem:filtering-filtered} shows that $F$ inverts all ingressives with cofibre in $\cU$.
  Hence we have an equivalence
  \begin{equation}\label{eq:W-quotient-homs}
   F \simeq \ell^*\yo_{\cC[(\cC_\ing^\cU)^{-1}]}(\ell(x)) \simeq \ell^*\ell_!\yo_\cC(x),
  \end{equation}
  as can be seen by applying Step~1 in the proof of \cite[Theorem~6.9]{hls:localisation} to the subcategory $(\cC^\cU_\ing)^\op \subseteq \cC^\op$ (note that the subcategory need not be closed under pushouts for Step~1), or by appealing to \cite[Theorem~7.2.8]{cisinski:higher-cats}.
  
  Let $X \in \psh_\Sigma(\cC)$. For $x \in X$, consider the canonical morphism
  \begin{equation}\label{eq:W-quotient}
   \colim_{x' \rightarrowtail x \in \cC_\ing^\cU(x)^\op} X(x') \to \colim_{x' \rightarrowtail x \in \cC_\ing^\cU(x)^\op} \ell^*\ell_!X(x') \simeq \ell^*\ell_!X(x),
  \end{equation}
  where the equivalence arises from the fact that $\cC_\ing^\cU(x)^\op$ is filtered by \cref{lem:filtering-filtered}.
  Call $X$ \emph{good} if the map \eqref{eq:W-quotient} is an equivalence for all $x \in X$.
  The preceding argument shows that $\yo_\cC(y)$ is good for all $y \in \cC$.
  To see that every additive presheaf is good, it thus suffices to show that the collection of good presheaves is closed under (sifted) colimits.
  If $X \colon I \to \psh_\Sigma(\cC)$ is a diagram of good presheaves, we observe that
  \[ \colim_{x' \to x \in \cC_\ing^\cU(x)^\op} \colim_{i \in I} X_i(x') \simeq \colim_{i \in I} \colim_{x' \to x \in \cC_\ing^\cU(x)^\op} X_i(x')\]
  since colimits commute with each other, and analogously for $\ell^*\ell_!X$ since both $\ell^*$ and $\ell_!$ preserve colimits.
  Therefore, all additive presheaves on $\cC$ are good.
  This tells us that $\ell_!X$ can be computed as a filtered colimit of additive presheaves, and the component $X(x) \to \ell^*\ell_!X(x)$ of the unit map becomes identified with the structure map $X(x) \to \colim_{x' \rightarrowtail x \in \cC_\ing^\cU(x)^\op} X(x')$ of the colimit.
  In particular, we conclude that $\ell_!$ preserves finite limits, and therefore so does $\ell$.
  
  Assertion~\eqref{it:W-quotient4} follows directly from this because \eqref{eq:W-quotient-homs} expresses mapping anima in the localisation as filtered colimits of discrete anima.
  
  Since $\ell$ preserves finite limits, the colimit formula for mapping anima in $\cD$ allows us to apply \cref{lem:localisation-Waldhausen} to refine $\ell$ to an exact functor of coWaldhausen categories.
  Observing that every morphism of coWaldhausen categories $\cC \to \cD$ which vanishes on $\cU$ has to invert $\cC_\ing^\cU$, it is straightforward to verify assertion~\eqref{it:W-quotient1}.

  The explicit formula for $\ell_!$ also shows that $\ell_!X \simeq 0$ for an additive presheaf if and only if $X$ is weakly $\cC_\ing^\cU$-effaceable.

  For the non-trivial direction of assertion~\eqref{it:W-quotient2}, assume that $\ell(f) \simeq 0$.
 Then the unit map $u(x) \colon \yo(y)(x) \to \ell^*\ell_!\yo(y)(x)$ sends $f$ to $0$, so $f$ lifts to an element $\phi \in \pi_1(\cofib(u)(x))$.	
  Since $\cofib(u)$ is weakly $\cC_\ing^\cU$-effaceable, there exists a morphism $i \colon x' \rightarrowtail x$ in $\cC_\ing^\cU$ such that $i^*\phi = 0 \in \pi_1(\cofib(u)(x'))$.
  By naturality, this implies $f \circ i \simeq 0$, which entails that $f$ factors through $\cofib(i) \in \cU$.

  We are left with showing \eqref{it:W-quotient3}.
  If $\yo(f) \colon \yo(x) \to \yo(y)$, considered as a morphism in $\psh_\Sigma(\cC)$, has the property that $\ell_!\yo(f)$ is an equivalence, then $\cofib(\yo(f)) \in \ker(\ell_!)$.
  Hence $\cofib(\yo(f))$ is weakly $\cC_\ing^\cU$-effaceable.
  In particular, there exists some $i \colon x' \rightarrowtail y$ in $\cC^\cU_\ing$ such that the composition
  \[ \yo(x') \xrightarrow{\yo(i)} \yo(y) \to \cofib(f) \]
  is nullhomotopic. 
  Hence $f$ fits into a commutative diagram
  \[\begin{tikzcd}
	  x'\ar[dr, rightarrowtail, "i"]\ar[d, "g"'] & \\
	  x\ar[r, "f"] & y
  \end{tikzcd}\]
  Now take the pushout
  \[\begin{tikzcd}
	   x'\ar[r, tail, "i"]\ar[d, "g"'] & y\ar[d, "h"] \\
	   x\ar[r, tail, "j"] & z
  \end{tikzcd}\]
  in $\cC$.
  The morphism $f$ induces a retraction $r \colon z \to y$ to $h$, so there exists an object $\widehat{u}$ in $\idem{\cC}$ such that $h$ is equivalent to the summand inclusion $y \rightarrowtail y \oplus \widehat{u}$ in $\idem{\cC}$.
  Quillen's obscure axiom (\cref{lem:obscure-axiom}) implies that $g$ is ingressive in $\idem{\cC}$ with cofibre $\widehat{u}$.

  Let $p \colon \yo(z) \twoheadrightarrow \widehat{u}$ be the cofibre of $\yo(h)$, and let $s \colon \widehat{u} \to \yo(z)$ be a section of $p$.
  Since $\ell(j)$ and $\ell(f)$ are equivalences, the same is true for $\ell(h)$.
  The argument in the preceding paragraph implies that there exists a morphism $y' \to y$ such that the composite $y' \to z$ lies in $\cC_\ing^\cU$.
  Denote by $q \colon z \twoheadrightarrow u$ the cofibre of this morphism.
  Then we obtain an induced morphism $e \colon \yo(u) \to \widehat{u}$ such that $e\yo(q) \simeq p$.
  It follows that $\yo(q)s \colon \widehat{u} \to u$ and $e \colon u \to \widehat{u}$ are morphisms satisfying $e\yo(q)s \simeq ps \simeq \id_{\widehat{u}}$.
  So $\widehat{u}$ is a retract of $\yo(u)$, and thus lies in $\idem{\cU}$, showing assertion~\eqref{it:W-quotient3}.
 %
\end{proof}

\begin{rem}\label[rem]{rem:calculus-of-fractions}
    The proof of \cref{thm:W-quotient} in particular shows that morphisms in the localisation $\overline{\cC} := \cC[(\cC_\egr^\cU)^{-1}]$ can be described by a calculus of fractions.
    Concretely, we have shown that for $x,y \in \cC$ we have
    \[ \Hom_{\overline{\cC}}(\ell(x),\ell(y)) \simeq \colim_{y \twoheadrightarrow y' \in \cC_\egr^\cU(y)} \Hom_\cC(x,y'), \]
    and the indexing category $\cC_\egr^\cU(x)$ is filtered by \cref{lem:filtering-filtered}.
    Since $\pi_0 \colon \Spc \to \Set$ preserves colimits, we obtain
    \[ \pi_0\Hom_{\overline{\cC}}(\ell(x),\ell(y)) \simeq \colim_{y \twoheadrightarrow y' \in \cC_\egr^\cU(y)} \pi_0\Hom_\cC(x,y'). \]
    By the explicit description of filtered colimits in $\Set$, we conclude that every morphism in $\overline{\cC}$ is represented by a cospan $x \xrightarrow{f} y' \overset{q}{\twoheadleftarrow} y$ which corresponds to the morphism $\ell(q)^{-1}\ell(f)$.

    Combining the filtered colimit description with \cref{thm:W-quotient}.\eqref{it:W-quotient3}, it follows that every equivalence in $\overline{\cC}$ is represented by a cospan with both legs in $\cC_\egr^\cU$.
\end{rem}

\begin{cor}\label[cor]{cor:waldhausen-quotients}
    Let $\cC$ be an exact category and let $L \colon \cC \to \cD$ be a Bousfield localisation of the underlying category.
    Assume one of the following:
    \begin{enumerate}
     \item Every morphism in $\cC$ is egressive (ie $\cC^\op$ is prestable).
     \item $\cC$ is an ordinary category and the unit morphism $x \to Lx$ is egressive for every $x \in \cC$.
    \end{enumerate}
    Then $\cD$ inherits a Waldhausen structure making $L \colon \cC \to \cD$ the cofibre of the inclusion $\ker(L) \to \cC$ in $\wald$.
\end{cor}
\begin{proof}
    By \cref{lem:adjoint-filtering}, $\ker(L)$ is left filtering in $\cC$.
    Being a Bousfield localisation, $L$ is a Dwyer--Kan localisation at its collection of unit maps, which are egressives with fibre in $\ker(L)$ by assumption.
    Since $L$ preserves colimits and vanishes on $\ker(L)$, it inverts $\cC_\egr^{\ker(L)}$, so $L$ is a Dwyer--Kan localisation at $\cC_\egr^{\ker(L)}$.
    The corollary follows from \cref{thm:W-quotient}.
\end{proof}

\begin{ex}\label[ex]{ex:w-quotients}
 We can use \cref{thm:W-quotient} and \cref{cor:waldhausen-quotients} to identify some explicit Waldhausen quotients.
 \begin{enumerate}
  \item In the case of an ordinary additive category and a Karoubi-filtering full additive subcategory, we conclude that the usual quotient in ordinary additive categories is also the quotient in Waldhausen categories.
  \item For $R$ a noetherian ring, it follows from \cref{ex:filtering-subcats}~\eqref{ex:filtering-subcats-1} that the Waldhausen quotient of $\Mod_R^\fg$ by $\Mod_R^{\fg,\tors}$ is given by $\Mod_R^\tf$ equipped with the Waldhausen structure in which every monomorphism is ingressive.
  Note that the Waldhausen category $\Mod_R^\tf$ is usually not exact.
  \item Let $\cC$ be a stable category with a t-structure.
  By \cref{ex:filtering-subcats}~\eqref{ex:filtering-subcats-2}, the Waldhausen quotient of $\cC_{\leq 0}$ by the heart $\cC^\heartsuit$ is given by $\cC_{\leq -1}$ equipped with the Waldhausen structure $\cC_{\leq -1}^\mx$ in which every morphism is ingressive.

  Again, the Waldhausen quotient is evidently not an exact category.
  If the t-structure is bounded below, the exact quotient is actually trivial, as can be seen by the following argument.  
  Let $F \colon \cC_{\leq 0} \to \cD$ be an exact functor to another exact category such that $\cC^\heartsuit \subseteq \kernel(F)$. Then $F$ inverts the Postnikov section $x \to \tau_{\leq -1} x$ for all $x \in \cC_{\leq 0}$.
  Let $x \in \cC_{\leq -k}$ for $k > 0$. Then both $x \to 0$ and $\tau_{\leq -(k+1)}x \to 0$ are $\pi_0$-monomorphisms, so $F$ inverts the Postnikov section $x \to \tau_{\leq -(k+1)}x$ if and only if it inverts the suspended morphism
  \[ \Sigma x \to \Sigma \tau_{\leq -(k+1)}x \simeq \tau_{\leq -k} \Sigma x. \]
  After finitely many suspensions, this Postnikov section becomes $(-1)$-truncation, and therefore gets inverted.
  This implies that the Postnikov section $x \to \tau_{\leq -k}x$ is inverted by $F$ for all $x \in \cC_{\leq 0}$ and $k \geq 0$.
  Since the t-structure is bounded below, $\tau_{\leq -k}x \simeq 0$ for sufficiently large $k$, so $F \simeq 0$.
 \end{enumerate}
\end{ex}

\section{The generic fibration theorem}\label{sec:generic-fibration}

The proof of \cref{thm:localisation-intro} rests as usual on Waldhausen's generic fibration theorem \cite[Proposition~1.5.5]{waldhausen:algKspaces}, see also \cite[Theorem~8.11]{barwick:algK}.
In recalling the generic fibration theorem, we seize the opportunity to observe that there are two versions of the relative $S_\bullet$-construction, both of which model the cofibre of the induced map on K-theory.
For the sake of simplicity, we focus on inclusions of full Waldhausen subcategories.
Our discussion, in particular the use of Rezk's equifibration criterion, closely follows \cite[Section~2]{hls:localisation}.

Let us begin by recalling the relative $S_\bullet$-construction.
Denote by 
\[ \sigma^+ \colon \Delta \to \Delta \quad\text{and}\quad \sigma^- \colon \Delta \to \Delta \]
the endofunctors which add a maximal respectively a minimal object to each finite linear order.
Precomposition with these functors induces \emph{d\'ecalage functors}
\[ \dec^+ := - \circ (\sigma^+)^\op \colon \s\cX \to \s\cX,\quad \dec^- := - \circ (\sigma^-)^\op \colon \s\cX \to \s\cX \]
on the category of simplicial objects in any category $\cX$.
The $(p+1)$-th coface maps $d^{p+1} \colon [p] \to [p+1]$ and the $0$-th coface maps $d^0 \colon [p] \to [1+p]$ induce natural transformations
\[ \lambda \colon \dec^+ \Rightarrow \id \quad\text{and}\quad \phi \colon \dec^- \Rightarrow \id. \]
The inclusions $[0] \to [p] \sqcup [0] \cong [p+1]$ define a natural transformation $\con_{[0]} \Rightarrow \sigma^+$ from the constant functor on $[0]$, so we obtain a natural transformation
\[ \epsilon^+ \colon \dec^+ \Rightarrow \con_{(-)_0}. \]
Similarly, the inclusions $[0] \to [0] \sqcup [p] \cong [1+p]$  induce a natural transformation
\[ \epsilon^- \colon \dec^- \Rightarrow \con_{(-)_0}. \]
Analogously, the morphisms $[p+1] \cong [p] \sqcup [0] \to [0] \sqcup [0] \cong [1]$ and $[1+p] \cong [0] \sqcup [p] \to [0] \sqcup [0] \cong [1]$ induce natural transformations
\[ \con_{(-)_1} \Rightarrow \dec^+ \quad\text{and}\quad \con_{(-)_1} \Rightarrow \dec^-. \]

\begin{defn}
 Let $\cC$ be a Waldhausen category.
 A \emph{full Waldhausen subcategory} of $\cC$ is a full subcategory $\cU \subseteq \cC$ such that
 \[ \cU_\ing := \{ j \in \cU \cap \cC_\ing \mid \cofib(j) \in \cU \} \]
 determines a Waldhausen structure on $\cU$.
\end{defn}

\begin{defn}
 Let $\cC$ be a Waldhausen category and let $\cU \subseteq \cC$ be a full Waldhausen subcategory.
 Define the \emph{relative $S_\bullet$-constructions} of $\cU \subseteq \cC$ as the pullbacks
 \[\begin{tikzcd}
  S_\bullet^+(\cU \subseteq \cC)\ar[r]\ar[d] & \dec^+(S_\bullet\cC)\ar[d, "\lambda"] \\
  S_\bullet\cU\ar[r] & S_\bullet\cC
 \end{tikzcd}
 \quad\text{and}\quad
 \begin{tikzcd}
  S_\bullet^-(\cU \subseteq \cC)\ar[r]\ar[d] & \dec^-(S_\bullet\cC)\ar[d, "\phi"] \\
  S_\bullet\cU\ar[r] & S_\bullet\cC
 \end{tikzcd}\]
 The natural transformations $\cC \simeq \con_{S_1\cC} \Rightarrow \dec^{\pm}(S_\bullet\cC)$ induce natural transformations
 \[ \cC \Rightarrow S_\bullet^{\pm}(\cU \subseteq \cC). \]
\end{defn}

\begin{rem}
 The simplicial category $S_\bullet^-(\cU \subseteq \cC)$ is Waldhausen's relative $S_\bullet$-construction \cite[Definition~1.5.4]{waldhausen:algKspaces}.
 On the other hand, we may identify $S_q^+(\cU \subseteq \cC)$ with the full subcategory of $S_{q+1}\cC$ spanned by objects of the form
 \[\begin{tikzcd}[row sep=1em, column sep=1.5em]
  u_1\ar[r, rightarrowtail] & u_2\ar[r, rightarrowtail]\ar[d, twoheadrightarrow] & u_3\ar[r, rightarrowtail]\ar[d, twoheadrightarrow]&  \ldots\ar[r, rightarrowtail] & u_q\ar[r, rightarrowtail]\ar[d, twoheadrightarrow] & x\ar[d, twoheadrightarrow] \\
  & u_{1,2}\ar[r, rightarrowtail] & u_{1,3}\ar[r, rightarrowtail]\ar[d, twoheadrightarrow]& \ldots\ar[r, rightarrowtail] & u_{1,q}\ar[r, rightarrowtail]\ar[d, twoheadrightarrow] & x_1\ar[d, twoheadrightarrow] \\
  & & u_{2,3}\ar[r, rightarrowtail] & \ldots\ar[r, rightarrowtail] & \ldots\ar[r, rightarrowtail] & x_2\ar[d, twoheadrightarrow] \\
  & & & & & \ldots\ar[d, twoheadrightarrow] \\
  & & & & & x_q
 \end{tikzcd}\]
 such that forgetting the rightmost column (ie applying $d_{q+1}$) yields an object in $S_q\cU$.
\end{rem}

The proof of the generic fibration theorem uses Rezk's equifibration criterion; see \cref{app:rezk} for further details.

\begin{prop}[Rezk's equifibration criterion]\label[prop]{thm:rezk-equifib}
 Let $I$ be a small category and consider a pullback square
 \[\begin{tikzcd}
  W\ar[r]\ar[d, "\sigma"'] & X\ar[d, "\tau"] \\
  Y\ar[r] & Z
 \end{tikzcd}\]
 in $\Fun(I,\Spc)$.
 Suppose that $\tau$ is \emph{equifibred}, meaning that the commutative square
 \begin{equation}\label{eq:rezk-equifib}\begin{tikzcd}
  X(i)\ar[r, "X(\alpha)"]\ar[d, "\tau(i)"'] & X(j)\ar[d, "\tau(j)"] \\
  Z(i)\ar[r, "Z(\alpha)"] & Z(j)
 \end{tikzcd}\end{equation}
 is a pullback for every morphism $\alpha \colon i \to j$ in $I$.
 Then
  \[\begin{tikzcd}
  \colim_I W\ar[r]\ar[d, "\colim_I \sigma"'] & \colim_I X\ar[d, "\colim_I \tau"] \\
  \colim_I Y\ar[r] & \colim_IZ
 \end{tikzcd}\]
 is a pullback in $\Spc$.
\end{prop}

\begin{rem}
    We will apply \cref{thm:rezk-equifib} in the case that $I = \Delta^\op$.
    In that case, the pasting lemmas for pullback squares imply that $\tau$ is equifibred if \eqref{eq:rezk-equifib} is a pullback for every face map.
\end{rem}

Recall from the \hyperlink{target:conventions}{Conventions section} that we use the term ``additive invariant'' in Barwick's sense \cite{barwick:algK}.


\begin{theorem}[Generic fibration theorem]\label[theorem]{thm:gen-fib}
 Let $\cC$ be a Waldhausen category and let $\cU \subseteq \cC$ be a full Waldhausen subcategory.
 Let $H \colon \wald \to \Spc$ be an additive invariant.
 Then the squares
 \[\begin{tikzcd}[column sep=1em]
  H(\cU)\ar[r]\ar[d] & H(\cC)\ar[d] \\
  * \simeq \abs{H(S_\bullet^+(\cU \subseteq \cU))}\ar[r] & \abs{H(S_\bullet^+(\cU \subseteq \cC))}
 \end{tikzcd}\]
 and
 \[\begin{tikzcd}[column sep=1em]
  H(\cU)\ar[r]\ar[d] & H(\cC)\ar[d] \\
  * \simeq \abs{H(S_\bullet^-(\cU \subseteq \cU))}\ar[r] & \abs{H(S_\bullet^-(\cU \subseteq \cC))}
 \end{tikzcd}\]
 are pullbacks of anima.
\end{theorem}
\begin{proof}
 We only provide the proof for $S_\bullet^+(\cU \subseteq \cC)$, which is slightly more non-standard.
 The proof for $S_\bullet^-$ is completely analogous.
 
 For the lower left corner, note that $S_\bullet^+(\cU \subseteq \cU)$ is $\dec^+(S_\bullet\cU)$, which is a split simplicial object.
 Hence $\abs{H(S_\bullet^-(\cU \subseteq \cU))} \simeq H(S_0\cU) \simeq *$ by \cite[Lemma~6.1.3.17]{HTT}.

 Additivity implies that the functors
 \begin{align*}
  \cU \times S_q^+(\cU \subseteq \cC) &\to S_{q+1}^+(\cU \subseteq \cC) \\
  (u,u_1 \rightarrowtail \ldots \rightarrowtail u_q \rightarrowtail x) &\mapsto (u \rightarrowtail u \sqcup u_1 \rightarrowtail \ldots \rightarrowtail u \sqcup u_q \rightarrowtail u \sqcup x)
 \end{align*}
 induce equivalences
  \[ H(\cU) \times H(S_q^+(\cU \subseteq \cC)) \to H(S_{q+1}^+(\cU \subseteq \cC)). \]
 By induction, we conclude that the functor
 \[ \cU^q \times \cC \to S_q^+(\cU \subseteq \cC),\quad ((u_i)_i,x) \mapsto (u_1 \rightarrowtail u_1 \sqcup u_2 \rightarrowtail \ldots \rightarrowtail u_1 \sqcup \ldots \sqcup u_q \sqcup x) \]
 induces an equivalence upon application of $H$.
 In particular, we have
 \[ H(S_q^+(\cU \subseteq \cU)) \simeq H(\cU)^{q+1}. \]
 Before realisation, the degree $q$ part of the square in question has been identified with the square
 \[\begin{tikzcd}
   * \times H(\cU)\ar[r]\ar[d] & * \times H(\cC)\ar[d] \\
   H(\cU)^q \times H(\cU)\ar[r] & H(\cU)^q \times H(\cC)
 \end{tikzcd}\]
 which is a product of pullback squares. Hence it is a pullback.
 We want to show that the transformation
 \[ H(S_\bullet^+(\cU \subseteq \cU)) \to H(S_\bullet^+(\cU \subseteq \cC)) \]
 is equifibred.
 The preceding argument also shows that the commutativity square for the face map $d_i$ is identified with a square of the form
   \[\begin{tikzcd}
   H(\cU)^{q+1} \times H(\cU)\ar[r, "\delta_i"]\ar[d] & H(\cU)^q \times H(\cU)\ar[d] \\
   H(\cU)^{q+1} \times H(\cC)\ar[r, "\delta_i"] & H(\cU)^q \times H(\cC)
  \end{tikzcd}\]
 From the explicit identification, we obtain the following description of the maps $\delta_i$:
  \begin{enumerate}
   \item for $i=0$, the map $\delta_i$ is projection away from the first factor;
   \item for $1 \leq i \leq q$, the map $\delta_i$ is induced by the fold map $H(\cU) \times H(\cU) \to H(\cU)$ on the $i$-th and $(i-1)$-th factor;
   \item for $i = q+1$, the map $\delta_i$ is induced by the map $\delta_{q+1}' \colon H(\cU) \times H(\cC) \to H(\cC) \times H(\cC) \to H(\cC)$ involving the fold map on $H(\cC)$ and the last copy of $H(\cU)$).
  \end{enumerate}
 In the first two cases, the induced map on vertical fibres is obviously an equivalence.
  For the last case, consider the commutative cube
   \[\begin{tikzcd}[row sep=1em]
 & H(\cU) \times H(\cU)\ar[dl, "\simeq"']\ar[rr, "\delta_{q+1}'"]\ar[dd] & & H(\cU)\ar[dl, "\id"]\ar[dd] & \\
  H(\cU) \times H(\cU)\ar[rr, crossing over, "\pr" xshift=1em]\ar[dd] & & H(\cU) \\
  & H(\cU)\times H(\cC)\ar[rr, "\delta_{q+1}'" xshift=-1.5em]\ar[dl, "\simeq"'] & & H(\cC)\ar[dl, "\id"] \\
  H(\cU) \times H(\cC)\ar[rr, "\pr" xshift=1em] & &  H(\cC)\ar[uu, leftarrow, crossing over]
 \end{tikzcd}\]
 in which the diagonal arrows on the left are given by the appropriate version of the shear map.
 Since $H$ is group-like, these maps are equivalences as indicated, and it follows that the square in question is also a pullback in this case.

  The theorem follows by applying Rezk's equifibration criterion. 
\end{proof}

We will specifically be interested in applying \cref{thm:gen-fib} to the additive invariant $\abs{\iota S_\bullet(-)}$, in which case we obtain a fibre sequence
\[ \K(\cU) \to \K(\cC) \to \Omega\abs{\iota S_\bullet S_\bullet^{\pm}(\cU \subseteq \cC)} \]
for each full Waldhausen subcategory $\cU \subseteq \cC$.
The next statement gives a criterion to identify the relative term with the K-theory of another Waldhausen category, and summarises well-known arguments (see eg \cite{staffeldt:fundamental,hls:localisation}).

\begin{prop}\label[prop]{prop:identify-cofibre}
 Let $\cC$ be a Waldhausen category, let $\cU \subseteq \cC$ be a full Waldhausen subcategory, and let $P \colon \cC \to \cD$ be an exact functor.
 \begin{enumerate}
     \item\label{prop:identify-cofibre-1} Suppose that the collection of morphisms
        \[ \cC_\ing^\cU := \{ i \colon x \rightarrowtail y \mid i\text{ is ingressive with cofibre in $\cU$} \} \]
    forms a subcategory of $\cC$.
    If $P$ inverts $\cC_\ing^\cU$ and induces an equivalence
    \[ \abs{(S_p\cC)_\ing^{S_p\cU}} \xrightarrow{\sim} \iota S_p\cD \]
    for all $p \geq 1$, then
    \[ \K(\cU)\to \K(\cC) \to \K(\cD) \]
    is a cofibre sequence of connective spectra.
    \item\label{prop:identify-cofibre-2} Suppose $\cC$ is exact and that $\cU$ is an extension-closed subcategory.
    Set
    \[ \cC_\egr^\cU := \{ p \colon y \twoheadrightarrow z \mid p\text{ is egressive with fibre in }\cU \}. \]
    If $P$ inverts $\cC_\egr^\cU$ and induces an equivalence
    \[ \abs{(S_p\cC)_\egr^{S_p\cU}} \xrightarrow{\sim} \iota S_p\cD \]
    for all $p \geq 1$, then
    \[ \K(\cU)\to \K(\cC) \to \K(\cD) \]
    is a cofibre sequence of connective spectra.
 \end{enumerate}
\end{prop}

Before giving the proof, recall that the category of small categories $\catinf$ embeds fully faithfully into the category of simplicial anima via the \emph{Rezk nerve}
\[ \N^\rezk \colon \catinf \to \s\Spc, \quad \cC \mapsto \big[ [n] \mapsto \iota\Fun([n],\cC) \big]. \]
The essential image of the Rezk nerve is given by the full subcategory of complete Segal anima, and the Rezk nerve admits a left adjoint $\mathrm{ac} \colon \s\Spc \to \catinf$.
By inspection, there exists a commutative diagram
\[\begin{tikzcd}
    & \Spc\ar[dl, "\inc"']\ar[dr, "\con"] & \\
    \catinf\ar[rr, "\N^\rezk"] & & \s\Spc 
\end{tikzcd}\]
Passing to left adjoints, we obtain a natural equivalence $\colim_{\Delta^\op} \simeq \abs{\mathrm{ac}(-)}$.
Denoting $\colim_{\Delta^\op}$ also by $\abs{-}$ as before, we obtain a natural equivalence
\[ \abs{\N^\rezk(\cC)} \simeq \abs{\cC} \]
for every small category $\cC$.
See \cite[Section~2.1]{the-nine:2} and \cite{hs:rezk-nerve} for further details.

\begin{proof}[Proof of \cref{prop:identify-cofibre}]
    In both cases, the assumptions imply that the induced map $\K_0(\cC) \to \K_0(\cD)$ is surjective.
    Consequently, it suffices to check that we have fibre sequences on the level of the underlying infinite loop spaces.

    Denote by $\Fun^\simeq([p],\cD)$ the full subcategory of $\Fun([p],\cD)$ on those functors which invert all morphisms in $[p]$, and let $\Fun^\simeq(\bullet,\cD)$ be the associated simplicial category.

    Forgetting about all choices of subquotients and applying the functor $P$ induces a transformation
    \[ S_\bullet^-(\cU \subseteq \cC) \to \Fun^\simeq(\bullet,\cD) \]
    through exact functors because $P$ is assumed to invert $\cC_\ing^\cU$.
    Since taking constant functors defines an equivalence $\con_\cD \xrightarrow{\sim} \Fun^\simeq(\bullet,\cD)$, applying $\iota S_\bullet$ yields a transformation
     \begin{equation}\label{eq:gen-cofib-comparison}
      \iota S_\bullet S_\bullet^-(\cU \subseteq \cC) \to \iota S_\bullet \con_\cD
    \end{equation}
    Fixing $p \in \bbN$, the simplicial category $S_pS_q^-(\cU \subseteq \cC)$ is equivalent to the full subcategory of $\Fun([p-1] \times [q],\cC)$ on objects given by commutative diagrams
 \[\begin{tikzcd}[column sep=1.5em, row sep=1em]
  x_0^1\ar[r, rightarrowtail]\ar[d,rightarrowtail] &  x_0^2\ar[r, rightarrowtail]\ar[d,rightarrowtail] &   \ldots\ar[r, rightarrowtail] &   x_0^p\ar[d,rightarrowtail] \\
  x_1^2\ar[r, rightarrowtail]\ar[d,rightarrowtail] &  x_1^2\ar[r, rightarrowtail]\ar[d,rightarrowtail] &   \ldots\ar[r, rightarrowtail] &   x_1^p\ar[d,rightarrowtail] \\
  \vdots\ar[d, rightarrowtail] & \vdots\ar[d, rightarrowtail] & & \vdots\ar[d, rightarrowtail] \\
  x_q^1\ar[r, rightarrowtail] &   x_q^2\ar[r, rightarrowtail] &   \ldots\ar[r, rightarrowtail] &   x_q^p
 \end{tikzcd}\]
 such that all vertical subquotients lie in $\cU$ and the horizontal maps define ingressive morphisms in $S_{q+1}\cC$.
 The latter condition is equivalent to requiring the vertical transformations to define ingressive morphisms in $S_p\cC$, which yields a natural identification
 \[ S_pS_\bullet^-(\cU \subseteq \cC) \simeq S_\bullet^-(S_p\cU \subseteq S_p\cC). \]
 Fixing the first simplicial direction in \eqref{eq:gen-cofib-comparison}, the degree $p$ part of this transformation is therefore identified with the map
 \[ \iota S_\bullet^-(S_p\cU \subseteq S_p\cC) \to \iota S_p\cD \]
 which forgets about choices of subquotients and applies $S_pP$.
  
 Since the $0$-th face map of the $S_\bullet$-construction does not feature in the relative $S_\bullet$-construction, forgetting about choices of subquotients defines an equivalence
  \[ \iota S_\bullet^-(S_p\cU \subseteq S_p\cC) \xrightarrow{\sim} \N^\rezk((S_p\cC)_\ing^{S_p\cU}). \]
 Hence, the degree $p$ component of \eqref{eq:gen-cofib-comparison} becomes the map
 \[ \abs{(S_p\cC)_\ing^{S_p\cU}} \to \iota S_p\cD \]
 induced by $S_pP$ after realisation.
 Since all of these maps are equivalences by assumption, it follows that \eqref{eq:gen-cofib-comparison} is an equivalence.
 Assertion~\eqref{prop:identify-cofibre-1} now follows from \cref{thm:gen-fib}.

 For the second assertion, we proceed similarly.
 Remembering only the rightmost column of objects in $S_\bullet^+(\cU \subseteq \cC)$ and applying $P$ induces a transformation
 \[ S_\bullet^+(\cU \subseteq \cC) \to \Fun^\simeq(\bullet,\cD) \]
 because $P$ inverts $\cC_\egr^\cU$.
 Therefore, we obtain a map
  \begin{equation}\label{eq:gen-cofib-comparison-2}
  \iota S_\bullet S_\bullet^+(\cU \subseteq \cC) \to \iota S_\bullet \con_\cD.
  \end{equation}
 For fixed $p$, we have an identification
 \[ S_pS_\bullet^+(\cU \subseteq \cC) \simeq S_\bullet^+(S_p\cU \subseteq S_p\cC) \]
 so that \eqref{eq:gen-cofib-comparison-2} becomes
 \[ \iota S_\bullet^+(S_p\cU \subseteq S_p\cC) \to \iota S_p\cD \]
 in degree $p$.
 Exactness of $\cC$ implies that restricting to the last column induces an equivalence between $S_q^+(S_p\cU \subseteq S_p\cC)$ and the full subcategory of $\Fun([q],S_p\cC)$ on functors $[q] \to (S_p\cC)_\egr^{S_p\cU}$.
 Since the last face maps do not feature in $S_\bullet^+$, these equivalences assemble to an equivalence
 \[ \iota S_\bullet^+(S_p\cU \subseteq S_p\cC) \xrightarrow{\sim} \N^\rezk((S_p\cC)_\egr^{S_p\cU}). \]
 Therefore, the degree $p$ component of \eqref{eq:gen-cofib-comparison-2} realises to
 \[ \abs{(S_p\cC)_\egr^{S_p\cU}} \to \iota S_p\cD. \]
 Since these maps are assumed to be equivalences, assertion~\eqref{prop:identify-cofibre-2} also follows from \cref{thm:gen-fib}.
\end{proof}

\section{The localisation theorem}\label{sec:localisation}

This section is dedicated to the proof of \cref{thm:localisation-intro}, which occupies \cref{sec:the-proof}.
After finishing the proof, we will explain in \cref{sec:no-idempotent-completeness} how the theorem can be generalised to strongly filtering subcategories which are not idempotent complete.

For a Waldhausen category $\cC$, we will denote by $F_p\cC$ the category of filtered objects $x_1 \rightarrowtail \ldots \rightarrowtail x_p$ in $\cC$.
This category carries an induced Waldhausen structure in which the ingressive morphisms are the Reedy cofibrant transformations.
With this definition, the forgetful functor $S_p\cC \to F_p\cC$ is an equivalence of Waldhausen categories.

\subsection{The localisation theorem for a Bousfield localisation}\label{sec:the-easy-proof}
    Before proving the general case of \cref{thm:localisation-intro}, let us point out that the theorem admits a very short proof if the left filtering subcategory $\cU$ arises as the kernel of a Bousfield localisation $L \colon \cC \to \cD$ as in \cref{lem:adjoint-filtering}.
    See \cref{lem:adjoint-strongly-filtering} below for a proof that such a subcategory $\cU$ is in fact strongly filtering, so we are proving a special case of the theorem.
    In this case, we claim that we have a cofibre sequence
    \[ \K(\cU) \to \K(\cC) \to \K(\cD). \]
    Since the unit morphisms are egressives with fibres in $\cU$, \cref{thm:W-quotient} exhibits $\cD$ as the Waldhausen quotient of $\cC$ by $\cU$.
    
    By \cref{prop:identify-cofibre}, it suffices to show that
    \[ \abs{(S_p\cC)_\egr^{S_p\cU}} \to \iota S_p(\cC[(\cC_\egr^\cU)^{-1}]) \]
    is an equivalence for all $p$.
    Through the identification $S_p\cC \simeq F_p\cC$, we want to show that the functor
    \[ \lambda \colon (F_p\cC)_\egr^{F_p\cU} \to \iota F_p\cD \]
    induced by $L$ becomes an equivalence after realisation for all $p$.
    Note that the right adjoint $R \colon \cD \to \cC$ preserves ingressives, so it induces a functor
    \[ \rho \colon \iota F_p\cD \to (F_p\cC)_\egr^{F_p\cU} \]
    in the opposite direction.
    Evidently, $\lambda\rho \simeq \id$, and the unit transformation provides a natural transformation $\id \Rightarrow \rho\lambda$.
    The theorem follows.

\subsection{Proof of the main theorem}\label{sec:the-proof}
We now turn to the proof of the general case.
That is, we are considering an idempotent complete, strongly left filtering subcategory $\cU$ in an exact category $\cC$.
The filtering assumption implies through \cref{thm:W-quotient} that the Waldhausen quotient $\quotw{\cC}{\cU}$ is given by the Dwyer--Kan localisation $\cC[(\cC_\egr^\cU)^{-1}]$.
By \cref{prop:identify-cofibre}, it therefore suffices to show that
\[ \abs{(S_p\cC)_\egr^{S_p\cU}} \to \iota S_p(\cC[(\cC_\egr^\cU)^{-1}]) \]
is an equivalence for all $p$.

The next lemma, which is a minimal refinement of the $2$-Segal condition for $\iota S_\bullet\cC$, allows us to reduce arguments in $S_p\cC$ to the case $p=2$.

\begin{lem}\label[lem]{lem:2-segal}
 Let $p \geq 2$. The canonical functor
 \[ S_p\cC \to S_{\{0,1,2\}}\cC \mathop{\times}\limits_{S_{\{0,2\}}\cC} S_{\{0,2,\ldots,p\}}\cC \]
 is an equivalence of exact categories.
\end{lem}
\begin{proof}
 The assertion is trivial for $p=2$, so assume $p \geq 3$.
 Since the restriction along $[p] \to \Ar[p]$, $i \mapsto (0,i)$ induces an equivalence $S_p\cC \xrightarrow{\sim} F_p\cC$ with the category of filtered objects, the functor in question is identified with the canonical equivalence
 \[ F_p\cC \to F_{2}\cC \mathop{\times}\limits_{F_{\{2\}}\cC} F_{\{2,\ldots,p\}}\cC. \]
 We additionally claim that a morphism in $S_p\cC$ is ingressive if and only if its image is ingressive.
 So let
 \[ f = \{f_{i,j}\}_{0 \leq i < j \leq p} \colon \{x_{i,j}\}_{0 \leq i < j \leq p} \to \{y_{i,j}\}_{0 \leq i < j \leq p} \]
 be a morphism in $S_p\cC$ such that $f_{0,j}$ and $f_{j-1,j}$ are ingressive for all $j$.
 Then $f_{j-2,j}$ fits into a commutative diagram
\[\begin{tikzcd}
 x_{j-2,j-1}\ar[r, tail]\ar[d, tail, "f_{j-2,j-1}"] & x_{j-2,j}\ar[r, two heads]\ar[d, "f_{j-2,j}"] & x_{j-1,j}\ar[d, tail, "f_{j-1,j}"] \\
 y_{j-2,j-1}\ar[r, tail] & y_{j-2,j}\ar[r, two heads] & y_{j-1,j}
\end{tikzcd}\]
with both rows exact.
By \cref{cor:five-lemma}, it follows that $f_{j-2,j}$ is ingressive.
Proceeding inductively, it follows that all $f_{i,j}$ are ingressive.
\end{proof}

Next, we establish a number of auxiliary lemmas.

\begin{lem}\label[lem]{lem:cofibres-of-ingressives}
    If $\cU$ is an idempotent complete, left filtering subcategory of an exact category $\cC$, then $\cU$ is closed under cofibres of ingressives in the weak idempotent completion of $\cC$: if $i \colon u \to v$ is ingressive in $\widem{\cC}$ with $u,v \in \cU$, then $\cofib(i) \in \cU$ (and $i$ is ingressive in $\cU$).
\end{lem}
\begin{proof}
    Since $\cofib(i) \in \widem{\cC}$, there exists an object $x \in \cC$ such that $\cofib(i) \oplus x \in \cC$.
    Writing the inclusion $x \rightarrowtail \cofib(i) \oplus x$ as the composite $x \rightarrowtail v \oplus x \twoheadrightarrow \cofib(i) \oplus x$, it is clear that this morphism becomes invertible in $\quotw{\cC}{\cU}$.
    Then the projection $\cofib(i) \oplus x \to x$ also becomes invertible in $\quotw{\cC}{\cU}$, so
    \cref{thm:W-quotient} implies that there exists a morphism $p \colon x \to y$ in $\cC$ which is egressive in $\idem{\cC}$ with fibre in $\idem{\cU}$ and such that $0 + p \colon \cofib(i) \oplus x \to y$ lies in $\cC_\egr^\cU$.
    Evidently, $\fib(0+p) \simeq \cofib(i) \oplus \fib(p)$, so $\cofib(i)$ is a retract of an object in $\cU$.
    By idempotent completeness of $\cU$, we have $\cofib(i) \in \cU$.
\end{proof}

\begin{cor}\label[cor]{cor:egressives-cancellation}
    Let $\cU$ be an idempotent complete, left filtering subcategory of an exact category $\cC$.
    Let $p \colon x \to y$ and $q \colon y \to z$ be morphisms in $\cC$ such that both $p$ and $qp$ lie in $\cC_\egr^\cU$.
    Then $q$ also lies in $\cC_\egr^\cU$.
\end{cor}
\begin{proof}
    \cref{cor:ingressives-widem} implies that $q$ is egressive in the weak idempotent completion of $\cC$.
    It follows from \cref{lem:nine-lemma} that the fibre of $\fib(q)$ participates in an exact sequence
    \[ \fib(p) \rightarrowtail \fib(qp) \twoheadrightarrow \fib(q) \]
    in $\widem{\cC}$.
    \cref{lem:cofibres-of-ingressives} shows that $\fib(q)$ also lies in $\cU$.
\end{proof}

\begin{cor}\label[cor]{lem:automatic-egressive}
    Let $\cU$ be an idempotent complete, left filtering subcategory of an exact category $\cC$.
    Consider a commutative square
    \[\begin{tikzcd}
        x_0\ar[r, tail, "i"]\ar[d, two heads, "p_0"] & x_1\ar[d, two heads, "p_1"] \\
        y_0\ar[r, tail, "j"] & y_1
    \end{tikzcd}\]
    in $\cC$ in which the horizontal arrows are ingressive and the vertical arrows are egressive with fibre in $\cU$.
    Then the induced map $p_2 \colon x_2 \to y_2$ on horizontal cofibres is also egressive with fibre in $\cU$.
\end{cor}
\begin{proof}
    Set $x := y_1 \times_{y_2} x_2$.
    Then $p_1$ admits a factorisation
    \[ x_1 \twoheadrightarrow x \to y_1, \]
    where the first map is egressive by the dual of \cref{lem:exact-reedy}.
    Note that the fibre of $x_1 \twoheadrightarrow x$ is equivalent to $\fib(p_0)$.
    Therefore, the lemma follows from \cref{cor:egressives-cancellation}.
\end{proof}

The following consequence of \cref{lem:automatic-egressive} will be used later.

\begin{cor}\label[cor]{cor:lift-projections-of-egressives}
    Let $\cU$ be an idempotent complete, strongly left filtering subcategory of an exact category $\cC$.
    Let $x \in F_p\cC$, let $1 \leq k \leq p$, and suppose that $q \colon x_k \twoheadrightarrow y_k$ is a morphism in $\cC_\egr^\cU$.
    Then there exist an object $y' \in F_p\cC$, a morphism $q' \colon x \twoheadrightarrow y$ in $(F_p\cC)_\egr^{F_p\cU}$ and a morphism $g \colon y_k \to y_k'$ in $\cC_\egr^\cU$ such that $q'_k \simeq gq$.
\end{cor}
\begin{proof}
    For $p=2$, there are two cases to consider.
    If $k=1$, we can obtain $q'$ by taking the pushout of $q$ along $x_1 \rightarrowtail x_2$.
    If $k=2$, we find a commutative square
    \[\begin{tikzcd}
        x_1\ar[r, tail, "i"]\ar[d, two heads, "q'_1"'] & x_2\ar[d, two heads, "gq"] \\
        y_1\ar[r, tail, "j"] & y_2
    \end{tikzcd}\]
    because $\cU$ is strongly filtering.
    Note that $g \in \cC_\egr^\cU$ by \cref{cor:egressives-cancellation}.
    Then \cref{lem:automatic-egressive} shows that this defines the required morphism.

    For $p > 2$, the corollary follows from the case $p=2$ using \cref{lem:2-segal}.
\end{proof}

The next statement will allow us to reduce the proof to considering $S_1\cU \subseteq S_1\cC$, and also clarifies the role of the assumption on $\cU$ to be strongly filtering.

\begin{prop}\label[prop]{lem:sp-filtering}
    Let $\cC$ be an exact category and let $\cU$ be an idempotent complete, extension-closed subcategory.
    Then the following are equivalent:
    \begin{enumerate}
        \item $\cU$ is strongly left filtering in $\cC$;
        \item $S_p\cU$ is left filtering in $S_p\cC$ for all $p \geq 0$.
    \end{enumerate}
\end{prop}
\begin{proof}
    Suppose that both $\cU \subseteq \cC$ and $S_2\cU \subseteq S_2\cC$ are left filtering.
    Then we may consider for any ingressive morphism $v \rightarrowtail y$ with $v \in \cU$ and any object $x \rightarrowtail y \twoheadrightarrow z$ in $S_2\cC$ the morphism in $S_2\cC$ represented by the diagram
    \[\begin{tikzcd}
     0\ar[r ,tail]\ar[d] & v\ar[d, tail]\ar[r, two heads, "\id"] & v\ar[d] \\
     x\ar[r, tail] & y\ar[r, two heads] & z
    \end{tikzcd}\]
    A morphism in $(S_2\cC)_\egr^{S_2\cU}$ which annihilates the given morphism witnesses that $\cU$ is strongly left filtering in $\cC$.

    For the converse direction, we begin with the case $p=2$.
    The class of objects $u \in S_2\cC$ such that $\Hom_{S_2\cC}(u,-)$ is weakly $(S_2\cC)_\egr^{S_2\cU}$-effaceable is closed under extensions, so it suffices to consider objects of the form $u_f := u \xrightarrow{\id} u \twoheadrightarrow 0$ and $u_b := 0 \rightarrowtail u \xrightarrow{\id} u$ with $u \in \cU$.
    
    Since $\Hom_{S_2\cC}(u_f,x) \simeq \Hom_\cC(u,x_0)$ for every $x = (x_0 \rightarrowtail x_1 \twoheadrightarrow x_2) \in S_2\cC$, and every egressive $x_0 \twoheadrightarrow x_0'$ with fibre in $\cU$ extends to an egressive in $S_2\cC$ with fibre in $S_2\cU$ by taking pushouts, it is straightforward to check that $\Hom_{S_2\cC}(u_f,-)$ is weakly $(S_2\cC)_\egr^{S_2\cU}$-effaceable.

    Consider now an element $\alpha \in \pi_k\Hom_{S_2\cC}(u_b,x)$.
    Denote by $\alpha_l \in \pi_k\Hom_\cC(d_lu,d_lx)$ the induced element under the face map $d_l \colon S_2\cC \to \cC$.
    Suppose we have fixed a morphism $p_1 \colon x_1 \twoheadrightarrow y_1$ in $\cC_\egr^\cU$ such that $(p_1)_*\alpha_1 = 0 \in \pi_k\Hom_\cC(u_1,y_1)$; such a morphism always exists because $\cU$ is left filtering.
    Since $\cU$ is strongly left filtering, we obtain a commutative square
    \[\begin{tikzcd}
        x_0\ar[r, tail]\ar[d, two heads, "p_0"'] & x_1\ar[d, two heads, "gp_1"] \\
        x_0'\ar[r, tail] & x_1'
    \end{tikzcd}\]
    such that $g \colon y_1 \to x_1'$ and $gp_1$ lies in $\cC_\egr^\cU$.
    \cref{lem:automatic-egressive} implies that the induced morphism $x_2 \to \cofib(x_0' \rightarrow x_1')$ lies in $\cC_\egr^\cU$, so this defines an egressive morphism $q \colon x \twoheadrightarrow x'$ in $S_2\cC$ with fibre in $S_2\cU$.
    Observing that
    \[ \Hom_{S_2\cC}(u_b,x') \simeq \Hom_\cC(0,x_0') \mathop{\times}\limits_{\Hom_\cC(0,x_1')} \Hom_\cC(u,x_1') \simeq \Hom_\cC(u,x_1'), \]
    we conclude that $q_*\alpha = 0 \in \pi_k\Hom_{S_2\cC}(u_b,x')$.
    
    Let now $p > 2$ and consider $u \in S_p\cU$, $x \in S_p\cC$ and $\alpha \in \pi_k\Hom_{S_p\cC}(u,x)$.
    By \cref{lem:2-segal}, the restriction functor
    \[ (r_p,d_p) \colon S_p\cC \to S_{\{0,1,2\}}(\cC) \times_{S_{\{0,2\}}\cC} S_{\{0,2,\ldots,p}(\cC) \]
    is an equivalence.
    In particular, we have
    \[ \Hom_{S_p\cC}(u,x) \simeq \Hom_{S_{\{0,1,2\}}(\cC)}(r_p(u),r_p(x)) \mathop{\times}\limits_{\Hom_\cC(u_2,x_2)} \Hom_{S_{\{0,2,\ldots,p\}}(\cC)}(d_p(u),d_p(x)). \]
    By induction, there exists a morphism $q \colon d_p(x) \twoheadrightarrow z$ in $(S_{\{0,2,\ldots,p\}}(\cC))_\egr^{S_{\{0,2,\ldots,p\}}(\cU)}$ such that $q_*d_p(\alpha) = 0$.
    Let $q_2 \colon x_2 \twoheadrightarrow z_2$ be the corresponding component of $q$.
    The first part of the argument shows that there exists a morphism $f_2 \colon z_2 \twoheadrightarrow y_2$ in $\cC_\egr^\cU$ such that $f_2q_2 \colon x_2 \twoheadrightarrow y_2$ is the component of a morphism $f \colon r_p(x) \twoheadrightarrow y$ in $(S_{\{0,1,2\}}(\cC))_\egr^{S_{\{0,1,2\}}(\cU)}$ such that $f_*r_p(\alpha) = 0$.
    By taking iterated pushouts, we can extend $f_2$ to a morphism $g \colon z \twoheadrightarrow z'$ in $(S_{\{0,2,\ldots,p\}}(\cC))_\egr^{S_{\{0,2,\ldots,p\}}(\cU)}$.
    Note that $(gq)_*d_p(\alpha) = 0$ as well.

    It follows that the image of $\alpha$ in the group
    \[ \pi_k\big( \Hom_{S_{\{0,1,2\}}(\cC)}(r_p(u),y) \times_{\Hom_\cC(u_2,y_2)} \Hom_{S_{\{0,2,\ldots,p\}}(\cC)}(d_p(u),z') \big) \]
    lifts to an element $\beta \in \pi_{k+1}\Hom_\cC(u_2,y_2)$.
    Since $\cU$ is filtering, there exists a morphism $\phi \colon y_2 \twoheadrightarrow y_2'$ in $\cC_\egr^\cU$ which annihilates $\beta$.
    Invoking the first part of the proof again, we find a morphism $f' \colon y \twoheadrightarrow y'$ in $(S_{\{0,1,2\}}\cC)_\egr^{S_{\{0,1,2\}}\cU}$ whose component at $2$ factors over $\phi$, and taking iterated pushouts extends $f'_2$ to a morphism $g' \colon z' \twoheadrightarrow z''$ in $(S_{\{0,2,\ldots,p\}}\cC)_\egr^{S_{\{0,2,\ldots,p\}}\cU}$.
    It then follows that $(f'f,g'gq)_*\alpha = 0$, which is exactly what we needed to show.    
\end{proof}

Momentarily, we will have use for the following characterisation of cofinal functors between filtered categories.
See \cite[Proposition~3.2.2]{ks:cats-and-sheaves} for a variant of this statement for ordinary categories.

\begin{lem}\label[lem]{lem:filtered-cofinal}
    Let $I$ and $J$ be filtered categories and let $f \colon I \to J$ be a functor.
    Then the following are equivalent:
    \begin{enumerate}
     \item\label{lem:filtered-cofinal-1} the functor $f$ is cofinal;
     \item\label{lem:filtered-cofinal-2} the slice $j/f := I \times_J J_{j/}$ is filtered for all $j \in J$;
     \item\label{lem:filtered-cofinal-3} the following conditions hold:
      \begin{enumerate}
          \item\label{lem:filtered-cofinal-3a} for every $j \in J$ there exist $i \in I$ and a morphism $j \to f(i)$;
          \item\label{lem:filtered-cofinal-3b} for every morphism $f(i) \to j$ with $i \in I$ and $j \in J$, there exist $i' \in I$ and a morphism $j \to f(i')$ in $J$ such that that the composite morphism $f(i) \to f(i')$ lies in the image of $f$.
      \end{enumerate}
    \end{enumerate}
\end{lem}
\begin{proof}
    \eqref{lem:filtered-cofinal-2} implies \eqref{lem:filtered-cofinal-1} by Quillen's Theorem A.
    
    To show that \eqref{lem:filtered-cofinal-3} implies \eqref{lem:filtered-cofinal-2}, we want to use the local characterisation of filteredness \cite[\href{https://kerodon.net/tag/02PS}{Theorem~02PS}]{kerodon}.
    Let $(i,j \overset{\beta}{\to} f(i))$ and $(i',j \overset{\beta'}{\to} f(i'))$ be objects in $j/f$.
    Observe that the mapping space in the slice $j/f$ is identified with
    \begin{align*}
        \Hom_{j/f}((i,\beta),(i',\beta'))
        &\simeq \Hom_I(i,i') \mathop{\times}\limits_{\Hom_J(f(i),f(i'))} \Hom_{J_{j/}}(\beta,\beta') \\
        &\simeq \Hom_I(i,i') \mathop{\times}\limits_{\Hom_{J}(j,f(i'))} \{ \beta' \},
    \end{align*}
    where the last term denotes the fibre of the map induced by $f$ and precomposition with $\beta$.
    
    Since $I$ is filtered, there exist an object $i_+ \in I$ and morphisms $\alpha_+ \colon i \to i_+$ as well as $\alpha'_+ \colon i' \to i_+$.
    Then $f(\alpha_+)\beta$ and $f(\alpha'_+)\beta'$ define a map $S^0 \to \Hom_J(j,f(i_+))$.
    Since $J$ is filtered, there exist $j_+ \in J$ and $\beta_+ \colon f(i_+) \to j_+$ such that $\beta_+f(\alpha_+)\beta \simeq \beta_+ f(\alpha'_+)\beta'$.
    By assumption \eqref{lem:filtered-cofinal-3b}, there exist $i_\infty \in I$ and a morphism $\beta_\infty \colon j_+ \to f(i_\infty)$ such that $\beta_\infty\beta_+ \simeq f(\alpha_\infty)$ for some morphism $\alpha_\infty \colon i_+ \to i_\infty$.
    Hence
    \[ f(\alpha_\infty \alpha_+)\beta \simeq \beta_\infty\beta_+f(\alpha_+)\beta \simeq \beta_\infty \beta_+f(\alpha'_+)\beta'\simeq f(\alpha_\infty\alpha'_+) \beta', \]
    which witnesses that $\alpha_\infty\alpha_+ \colon i \to i_\infty$ and $\alpha_\infty\alpha'_+ \colon i' \to i_\infty$ refine to morphisms with the same target in $j/f$.
    
    The argument for maps $\xi \colon S^k \to \Hom_{j/f}((i,\beta),(i',\beta'))$ with $k \geq 0$ is similar.
    Observe that we have a fibre sequence
    \[ \Omega_{\beta'} \Hom_{J}(j,f(i')) \to \Hom_{j/f}((i,\beta),(i',\beta')) \overset{r}{\to} \Hom_I(i,i') \]
    which is natural with respect to morphisms in $j/f$.
    The filteredness of $I$ implies that there exists a morphism $\alpha_+ \colon i' \to i_+$ in $I$ such that $\alpha_+ \circ r(\xi)$ is nullhomotopic.
    Note that $\alpha_+$ also defines a morphism $\widetilde\alpha_+ \colon (i',\beta') \to (i_+,f(\alpha_+)\beta')$, so $(\widetilde\alpha_+)_*\xi$ lifts to a map $\widehat{\xi} \colon S^k \to \Omega_{f(\alpha_+)\beta'} \Hom_{J}(j,f(i_+))$.
    Since $J$ is filtered, there exists a morphism $\beta_+ \colon f(i_+) \to j_+$ such that $(\beta_+)_*\widehat{\xi}$ is nullhomotopic.
    Applying condition~\eqref{lem:filtered-cofinal-3b}, we find a morphism $\beta_\infty \colon j_+ \to f(i_\infty)$ such that $\beta_\infty\beta_+ \simeq f(\alpha_\infty)$ for some $\alpha_\infty \colon i_+ \to i_\infty$ in $I$.
    Again, $\alpha_\infty$ lifts to a morphism $\widetilde\alpha_\infty \colon (i_+,f(\alpha_+)\beta') \to (i',f(\alpha_\infty\alpha_+)\beta')$ in $j/f$, and we conclude that $(\widetilde\alpha_\infty\widetilde\alpha_+)_*\xi$ is nullhomotopic.
    Hence $j/f$ is filtered.

    To show that \eqref{lem:filtered-cofinal-1} implies \eqref{lem:filtered-cofinal-3}, note that the cofinality of $f$ implies in particular that
    \[ \colim_{i \in I} \Hom_J(j,f(i)) \to \colim_{j' \in J} \Hom_J(j,j') \]
    is an equivalence for all $j \in J$.
    The right hand term is non-empty, so \eqref{lem:filtered-cofinal-3a} holds.
    Similarly, the map
    \[ \colim_{i' \in I} \Hom_J(f(i),f(i')) \to \colim_{j \in J} \Hom_J(f(i),j) \]
    is an equivalence for all $i \in I$.
    By \eqref{lem:filtered-cofinal-3a}, every morphism $f(i) \to j$ represents the same class as some morphism $f(i) \to f(i')$.
    The colimit of the Yoneda embedding is the terminal presheaf, so both sides are contractible.
    Applying this observation to the left hand term, we find that every map $f(i) \to f(i')$ in $J$ represents the same element as $\id_{f(i)}$ in $\pi_0$.
    By the filteredness of $I$, this means that every map $f(i) \to f(i')$ becomes a map of the form $f(\alpha)$ after composition with some map $f(\beta) \colon f(i') \to f(i'')$.
    Hence \eqref{lem:filtered-cofinal-3b} holds as well.
\end{proof}

\begin{lem}\label[lem]{lem:sp-quotient}
Let $\cU$ be an idempotent complete, strongly left filtering subcategory of the exact category $\cC$.
Then the functor $\quotw{S_p\cC}{S_p\cU} \to S_p(\quotw{\cC}{\cU})$ is an equivalence of categories.
\end{lem}
\begin{proof}
 Recall that $S_p\cC$ is equivalent to the category $F_p\cC$ of filtered objects, so it suffices to show both statements for $F_p$.

 \cref{lem:sp-filtering} shows that $F_p\cU$ is left filtering in $F_p\cC$, so \cref{thm:W-quotient} describes $\quotw{S_p\cC}{S_p\cU}$ as a Dwyer--Kan localisation of $F_p\cC$.
 For essential surjectivity, it is therefore enough to show that $F_p\cC \to F_p(\quotw{\cC}{\cU})$ is essentially surjective.

 Using \cref{rem:calculus-of-fractions}, an ingressive morphism in $\quotw{\cC}{\cU}$ can be represented by a chain of morphisms
 \[ x \overset{r}{\twoheadrightarrow} x' \overset{r'}{\twoheadleftarrow} x'' \overset{j}{\rightarrowtail} y'' \overset{s'}{\twoheadrightarrow} y' \overset{s}{\twoheadleftarrow} y \]
 such that $j$ is ingressive and all other morphisms are egressives with fibre in $\cU$.
 Since $\cU$ is strongly left filtering in $\cC$, there exist a morphism $g \colon y' \to \widetilde{y}$ and a commutative square 
 \[\begin{tikzcd}
     x''\ar[r, tail, "j"]\ar[d, two heads, "p"] & y''\ar[d, two heads, "gs'"] \\
     \widetilde{x}\ar[r, tail, "\widetilde{j}"] & \widetilde{y} 
 \end{tikzcd}\]
 whose vertical morphisms lie in $\cC_\egr^\cU$.
 Since $\cU$ is closed under cofibres of ingressives in $\widem{\cC}$ by \cref{lem:cofibres-of-ingressives}, it follows that $g \in \cC_\egr^\cU$ as well.
 So we may represent the same morphism by the chain
 \[ x \overset{r}{\twoheadrightarrow} x' \overset{r'}{\twoheadleftarrow} x'' \overset{p}{\twoheadrightarrow} \widetilde{x} \overset{\widetilde{j}}{\rightarrowtail} \widetilde{y} \overset{gs}{\twoheadleftarrow} y. \]
 Taking the pullback of $\widetilde{j}$ along $gs$, we find that every ingressive morphism $x \rightarrowtail y$ in $\quotw{\cC}{\cU}$ can be lifted to an ingressive morphism in $\cC$ with target $y$.

 By induction, this implies directly that every object in $F_p(\quotw{\cC}{\cU})$ can be lifted to an object in $F_p\cC$.

 Next, we show that the comparison functor is fully faithful.
 Let $x,y \in S_p\cC$.
 Using the formula for anima of natural transformations (see \cite[Proposition~2.3]{glasman:hodge} or \cite[Proposition~5.1]{ghn:lax-colimits}), the functor $\quotw{F_p\cC}{F_p\cU} \to F_p(\quotw{\cC}{\cU})$ induces the comparison map
 \[ \colim_{y \to y' \in (F_p\cC)_\egr^{F_p\cU}(y)} \lim_{(i,j) \in \Tw[p]} \Hom_\cC(x_i,y_j') \to \lim_{(i,j) \in \Tw[p]} \colim_{y_j \to y_j' \in \cC_\egr^{\cU}(y_j)} \Hom_\cC(x_i,y_j'). \]
 Here, $\Tw[p]$ denotes the twisted arrow category on the poset $[p]$.
 Since the colimit is filtered, the domain is identified with
 \[ \lim_{(i,j) \in \Tw[p]} \colim_{y \to y' \in (F_p\cC)_\egr^{F_p\cU}(y)} \Hom_\cC(x_i,y_j'), \]
 and it suffices to show that the functor
 \[ (F_p\cC)_\egr^{F_p\cU}(y) \to \cC_\egr^{\cU}(y_j) \]
 induced by projection onto the $j$-th component is cofinal.
 We will check that condition~\eqref{lem:filtered-cofinal-3} of \cref{lem:filtered-cofinal} is satisfied.
 Condition~\eqref{lem:filtered-cofinal-3a} is immediate from \cref{cor:lift-projections-of-egressives}.
 For condition~\eqref{lem:filtered-cofinal-3b}, let $y \twoheadrightarrow z$ in $(F_p\cC)_\egr^{F_p\cU}$ be given and let $z_j \to y_j'$ be a morphism in $\cC$ such that the composite $y_j \to y_j'$ lies in $\cC_\egr^\cU$.
 By \cref{thm:W-quotient}, there exists a morphism $y_j' \to y_j''$ such that $z_j \to y_j''$ is a morphism in $\cC_ \egr^\cU$.
 Then \cref{cor:lift-projections-of-egressives} applies to $z_j \to y_j''$ to show that condition~\eqref{lem:filtered-cofinal-3b} is also satisfied.
\end{proof}

We are now prepared to finish the proof of \cref{thm:localisation-intro}.
We have a commutative diagram
\[\begin{tikzcd}
    S_p\cC\ar[dr]\ar[d] & \\
    \quotw{S_p\cC}{S_p\cU}\ar[r, "\sim"] & S_p(\quotw{\cC}{\cU})
\end{tikzcd}\]
in which the bottom horizontal arrow is an equivalence by \cref{lem:sp-quotient}.
Consequently, it suffices to show that
\[ \abs{(S_p\cC)_\egr^{S_p\cU}} \to \iota \quotw{S_p\cC}{S_p\cU} \]
is an equivalence.
Note that $S_p\cU$ is also idempotent complete.
As $S_p\cU$ is left filtering in $S_p\cC$ by \cref{lem:sp-filtering}, we are left with proving the following lemma.

\begin{lem}\label[lem]{lem:filtering-subcat-groupoid-core-localisation}
    Let $\cD$ be an exact category and let $\cV$ be an idempotent complete, left filtering subcategory of $\cD$.
    Then the induced map
    \[ \abs{\cD_\egr^\cV} \to \iota \quotw{\cD}{\cV} \]
    is an equivalence.
\end{lem}
\begin{proof}
    We will apply Quillen's Theorem A to the functor $\ell \colon \cD_\egr^\cV \to \iota \quotw{\cD}{\cV}$.
    Once more, \cref{thm:W-quotient} allows us to identify $\cQ := \quotw{\cD}{\cV}$ with the Dwyer--Kan localisation of $\cD$ at $\cD_\egr^\cV$.
    Since localisations are essentially surjective, it suffices to consider objects $\ell(x) \in \quotw{\cD}{\cV}$.
    Let
    \[ f = (f_\cD,f_\cQ) \colon K \to \cD_\egr^\cV \times_{\iota \cQ} (\iota \cQ)_{\ell(x)/} \]
    be an arbitrary finite diagram.
    We regard $f_\cQ$ as a natural transformation $\tau$ from the constant diagram $\con_{\ell(x)}$ on $\ell(x)$ to the diagram $\ell \circ f_\cD$ in $\iota\cQ$.
    Since we are working in the groupoid core of $\cQ$, the transformation $\tau$ is invertible and induces a transformation $\tau^{-1}$ in the opposite direction.
    The transformation $\tau^{-1}$ corresponds to a point in
    \begin{align*}
        \lim_{K^\op} \Hom_{\quotw{\cD}{\cV}}(\ell \circ f_\cD,\ell(x))
        &\simeq \lim_{K^\op} \colim_{x \twoheadrightarrow x' \in \cD_\egr^\cV(x)} \Hom_\cD(f_\cD,x') \\
        &\simeq \colim_{x \twoheadrightarrow x' \in \cD_\egr^\cV(x)} \lim_{K^\op} \Hom_\cD(f_\cD,x'),
    \end{align*}
    where the second equivalence holds since $\cD_\egr^\cV(x)$ is filtered.
    For an appropriate choice of $x'$, the transformation $\tau^{-1}$ therefore induces a transformation $\gamma \colon f_\cD \Rightarrow \con_{x'}$ in $\cC$ such that $\ell \circ \gamma$ extends to a cone on $f_\cQ$.
    All components of $\gamma$ become invertible in $\cQ$.
    By \cref{thm:W-quotient} and \cref{cor:egressives-cancellation}, it follows from the finiteness of $K$ that there exists a morphism $q \colon x'\twoheadrightarrow x''$ in $\cD_\egr^\cV$ such that $q \circ \gamma$ defines a cone of $f_\cD$ in $\cD_\egr^\cV$.
    Similarly, the cone on $f_\cQ$ induced by $\ell \circ \gamma$ gives rise to another cone with cone point $\ell(x'')$ via composition with $\ell(q)$, which finally yields the required cone on $f$.
    Hence $\cD_\egr^\cV \times_{\iota\cQ} \iota\cQ_{\ell(x)/}$ is filtered, and the lemma follows.
\end{proof}

Applying \cref{lem:filtering-subcat-groupoid-core-localisation}, we find that $\abs{(S_p\cC)_\egr^{S_p\cU}} \to \iota S_p(\quotw{\cC}{\cU})$ is an equivalence for all $p$, which concludes the proof of \cref{thm:localisation-intro}.

\subsection{The generalisation to non-idempotent complete categories}\label{sec:no-idempotent-completeness}
Just as for the localisation theorem for stable categories, there is a slightly more general version of \cref{thm:localisation-intro} which does not require $\cU$ to be idempotent complete.
Its formulation requires the following notation.

\begin{defn}
 Let $\cC$ be an exact category and let $\cU \subseteq \cC$ be an extension-closed subcategory.
 For a subgroup $A \subseteq \K_0(\idem{\cC})$, define
 \[ \Idemrel{A}{\cU} \subseteq \idem{\cU} \]
 to be the full subcategory consisting of those objects $\widehat{u}$ with $[\widehat{u}] \in A$.
\end{defn}

\begin{rem}
 The full subcategory $\Idemrel{A}{\cU} \subseteq \idem{\cU}$ is closed under extensions, and thus acquires an exact structure.
\end{rem}

\begin{theorem}\label[theorem]{thm:localisation}
 Let $\cC$ be an exact category and let $\cU \subseteq \cC$ be a strongly left filtering subcategory of $\cC$.
Then the projection onto the Waldhausen quotient induces a cofibre sequence
 \[ K(\Idemrel{\K_0(\cC)}{\cU}) \to K(\cC) \to K(\quotw{\cC}{\cU}) \]
 of connective spectra such that the composite map
 \[ K(\cU) \to K(\Idemrel{\K_0(\cC)}{\cU}) \to K(\cC) \]
 is induced by the inclusion of $\cU$ into $\cC$.
 Moreover, the map $K(\cU) \to K(\Idemrel{\K_0(\cC)}{\cU})$ is $(-1)$-truncated.
\end{theorem}

We will reduce the proof of \cref{thm:localisation} to an application of \cref{thm:localisation-intro} using arguments going back to \cite[Section~5]{pw:khomology} and \cite[Section~2]{hs:localisation-hermitian-k}.
Recall \cref{rem:equality-in-k0,rem:k0-of-idem} about $\K_0$ of an exact category. 

\begin{lem}\label[lem]{lem:idemrel}
    Let $\cC$ be an exact category and let $\cU \subseteq \cC$ be an extension-closed subcategory.
    \begin{enumerate}
        \item The category $\Idemrel{\K_0(\cC)}{\cU}$ is contained in the weak idempotent completion of $\cC$.
        \item The subcategory $\cU$ is dense in $\Idemrel{\K_0(\cC)}{\cU}$. In particular, the induced map $\K(\cU) \to \K(\Idemrel{\K_0(\cC)}{\cU})$ is $(-1)$-truncated.
    \end{enumerate}
\end{lem}
\begin{proof}
    Note that $\Idemrel{\K_0(\cC)}{\cU}$ is contained in $\Idemrel{\K_0(\cC)}{\cC}$.
    By \cref{lem:widem-explicitly}, $\Idemrel{\K_0(\cC)}{\cC}$ is the weak idempotent completion of $\cC$.

    Since $\cU$ is dense in $\idem{\cU}$, there exists for every $\widehat{u} \in \Idemrel{\K_0(\cC)}{\cU}$ an object $u' \in \idem{\cU}$ with $\widehat{u} \oplus u' \in \cU$.
    Since $\K_0(\cC)$ is a subgroup, it follows that $[u'] = [\widehat{u} \oplus u'] - [\widehat{u}] \in \K_0(\cC)$ as well.
    Hence $\cU$ is dense in $\Idemrel{\K_0(\cC)}{\cU}$.
    
    The assertion about the induced map on K-theory follows from the cofinality theorem \cite[Theorem~10.11]{barwick:algK}.
\end{proof}

For the reader familiar with the case of stable categories, the following lemma explains why the subcategory $\Idemrel{\K_0(\cC)}{\cU}$ appears in this formulation of the localisation theorem.
Recall that stable categories are weakly idempotent complete.

\begin{cor}\label[cor]{cor:retracts-partial-completion}
    Let $\cC$ be a weakly idempotent complete exact category and let $\cU$ be an extension-closed subcategory.
    Then the following are equivalent:
    \begin{enumerate}
        \item\label{lem:retracts-partial-completion-1} $\cU$ is closed under retracts in $\cC$;
        \item\label{lem:retracts-partial-completion-2} the inclusion $\cU \subseteq \Idemrel{\K_0(\cC)}{\cU}$ is an equivalence.
    \end{enumerate}
\end{cor}
\begin{proof}
    Assume \eqref{lem:retracts-partial-completion-1} and let $\widehat{u} \in \Idemrel{\K_0(\cC)}{\cU}$.
    Since $\cC$ is weakly idempotent complete, \cref{lem:idemrel} implies that  $\widehat{u} \in \cC$.
    By assumption, $\widehat{u}$ is a retract of an object in $\cU$, so it follows that $\widehat{u} \in \cU$.

    Conversely, assume \eqref{lem:retracts-partial-completion-2} and let $x \in \cC$ be a retract of an object $u \in \cU$.
    Since $\cC$ is weakly idempotent complete, we have $u \simeq x \oplus y$ for some $y \in \cC$.
    As $x$ is in particular an object in $\idem{\cU}$ and satisfies $[x] = [u] - [y] \in \K_0(\cC)$, it follows that $x \in \Idemrel{\K_0(\cC)}{\cU} \simeq \cU$.
\end{proof}

\begin{defn}
 Let $\cC$ be an exact category and let $\cV \subseteq \idem{\cC}$ be an extension-closed subcategory.
 Define
 \[ \idem{\cC}_\cV \subseteq \idem{\cC} \]
 as the full subcategory on those objects $\widehat{x} \in \idem{\cC}$ such that there exists $v \in \cV$ with $\widehat{x} \oplus v \in \cC$.
\end{defn}

Since $\idem{\cC}_\cV$ is closed under extensions in $\idem{\cC}$, it carries an induced exact structure.
Note that $\cC \subseteq \idem{\cC}_\cV$ because $\cV$ contains a zero object.

\begin{lem}\label[lem]{lem:pass-to-completions}
    Let $\cC$ be an exact category and let $\cU$ be an extension-closed subcategory.
    \begin{enumerate}
        \item\label{lem:pass-to-completions-0} If $A \subseteq \K_0(\idem{\cC})$ is a subgroup with $\K_0(\cC) \subseteq A$, then $\Idemrel{A}{\cU}$ is an extension-closed subcategory of $\idem{\cC}_{\Idemrel{A}{\cU}}$.
        \item\label{lem:pass-to-completions-1} The category $\idem{\cC}_{\Idemrel{\K_0(\cC)}{\cU}}$ is contained in $\widem{\cC}$, and both maps
        \[ \K(\cC) \to \K(\idem{\cC}_{\Idemrel{\K_0(\cC)}{\cU}}) \to \K(\widem{\cC}) \]
        are equivalences.
        \item\label{lem:pass-to-completions-2} The square
        \[\begin{tikzcd}
         \K(\Idemrel{\K_0(\cC)}{\cU})\ar[r]\ar[d] & \K(\idem{\cU})\ar[d] \\
         \K(\idem{\cC}_{\Idemrel{\K_0(\cC)}{\cU}})\ar[r] &  \K(\idem{\cC}_{\idem{\cU}})
        \end{tikzcd}\]
        is bicartesian.
        \item\label{lem:pass-to-completions-3} The induced functor
        \[ \quotw{\cC}{\cU} \to \quotw{\idem{\cC}_{\idem{\cU}}}{\idem{\cU}} \]
        is an equivalence.
    \end{enumerate}
\end{lem}
\begin{proof}
    For \eqref{lem:pass-to-completions-0}, we only have to show that $\Idemrel{A}{\cU}$ is actually contained in $\idem{\cC}_{\Idemrel{A}{\cU}}$.
    Observe that any $\widehat{u} \in \idem{\cU}$ with $[\widehat{u}] \in A$ admits a complement $u' \in \idem{\cU}$ with $\widehat{u} \oplus u' \in \cU \subseteq \cC$. Then $[u'] = [\widehat{u} \oplus u'] - [\widehat{u}] \in A$ because $\K_0(\cC) \subseteq A$ and $A$ is a subgroup.

    For \eqref{lem:pass-to-completions-1}, let $\widehat{x} \in \idem{\cC}$ such that there exists $\widehat{u} \in \Idemrel{\K_0(\cC)}{\cU}$ with $\widehat{x} \oplus \widehat{u} \in \cC$.
    Then $\widehat{u} \in \idem{\cC}_{\Idemrel{\K_0(\cC)}{\cU}}$ by assertion~\eqref{lem:pass-to-completions-0}, so $\cC$ is dense in $\idem{\cC}_{\Idemrel{\K_0(\cC)}{\cU}}$.
    \cref{lem:idemrel} shows that $\widehat{u} \in \widem{\cC}$, so $\widehat{x}$ is also an object of $\widem{\cC}$.
    The cofinality theorem \cite[Theorem~10.11]{barwick:algK} implies that
    both $\K(\cC) \to \K(\idem{\cC}_{\Idemrel{\K_0(\cC)}{\cU}})$ and $\K(\cC) \to \K(\widem{\cC})$ are $(-1)$-truncated, so it suffices to check that $\K_0(\cC) \to \K_0(\widem{\cC})$ is surjective.
    For every $\widehat{x} \in \widem{\cC}$, there exists some $x \in \cC$ with $\widehat{x} \oplus x \in \cC$, so
    the relation $[\widehat{x}] = [\widehat{x} \oplus x] - [x] \in \K_0(\widem{\cC})$ shows surjectivity.

    For assertion~\eqref{lem:pass-to-completions-2}, we will show that the horizontal cofibres agree.
    Using the cofinality theorem \cite[Theorem~10.11]{barwick:algK}, this amounts to the assertion that the induced homomorphism
    \[ \gamma \colon A := \K_0(\idem{\cU})/\K_0(\Idemrel{\K_0(\cC)}{\cU}) \to \K_0(\idem{\cC}_{\idem{\cU}})/\K_0(\idem{\cC}_{\Idemrel{\K_0(\cC)}{\cU}}) =:B \]
    is a bijection.

    By definition, any object $\widehat{x} \in \idem{\cC}_{\Idemrel{\K_0(\cC)}{\cU}}$ admits a complement $\widehat{u} \in \Idemrel{\K_0(\cC)}{\cU}$ with $\widehat{x} \oplus \widehat{u} \in \cC$.
    Hence $[\widehat{x}] = - [\widehat{u}] \in B$, so $\gamma$ is surjective.

    We are left with proving injectivity.
    By \cref{rem:k0-of-idem}, an arbitrary element in $\K_0(\idem{\cU})$ has a representative $[\widehat{u}] - [v]$ with $v \in \cU$, so suppose that $[\widehat{u}] - [v] = [x] - [y] \in \K_0(\idem{\cC}_{\idem{\cU}})$ with $x,y \in \cC$.
    Then there exist $a,b \in \cC$ and $\widehat{s} \in \idem{\cC}$ together with exact sequences $a \rightarrowtail \widehat{u} \oplus y \oplus \widehat{s} \twoheadrightarrow b$ and $a \rightarrowtail x \oplus v \oplus \widehat{s} \twoheadrightarrow b$.
    Since $\cC$ is extension-closed, it follows that $x \oplus v \oplus \widehat{s} \in \cC$.
    As $y \oplus \widehat{s} \oplus x \oplus v$ is also an object in $\cC$, it follows that $\widehat{u} \in \widem{\cC}$.
    Since $\K_0(\cC) \cong \K_0(\widem{\cC})$, this shows $\widehat{u} \in \Idemrel{\K_0(\cC)}{\cU}$ as required.

    For \eqref{lem:pass-to-completions-3}, it suffices to show that the restriction functor
    \[ \Fun^\exct_{\idem{\cU}}(\idem{\cC}_{\idem{\cU}},\cD) \to \Fun^\exct_\cU(\cC,\cD) \]
    is an equivalence for every small Waldhausen category $\cD$, where either side denotes the category of exact functors vanishing on the respective subcategory.
    This restriction functor features in a commutative square
    \[\begin{tikzcd}
        \Fun^\exct_{\idem{\cU}}(\idem{\cC}_{\idem{\cU}},\cD)\ar[r]\ar[d, "\yo_*"'] &  \Fun^\exct_\cU(\cC,\cD)\ar[d, "\yo_*"] \\
        \Fun^\exct_{\idem{\cU}}(\idem{\cC}_{\idem{\cU}},\psh_\lex(\cD))\ar[r] &  \Fun^\exct_\cU(\cC,\psh_\lex(\cD))
    \end{tikzcd}\]
    where we have used that a functor from an exact category to a pointed and cocomplete category is exact (with respect to the maximal Waldhausen structure on the target) if and only if it preserves finite coproducts and sends exact sequences to cofibre sequences, see \cite[Theorem~2.7]{nw:gabriel-quillen}.
    
    Any exact functor $\idem{\cC}_{\idem{\cU}} \to \psh_\lex(\cD)$ which vanishes on $\idem{\cU}$ inverts the projection morphisms $\widehat{x} \oplus \widehat{u} \to \widehat{x}$ with $\widehat{x} \in \idem{\cC}$ and $\widehat{u} \in \idem{\cU}$.
    Thus, if such a functor maps $\cC$ into $\cD$, it necessarily also maps $\idem{\cC}_{\idem{\cU}}$ to $\cD$,
    showing that the above square is a pullback.
    
    Observe that an exact functor $\idem{\cC}_{\idem{\cU}} \to \psh_\lex(\cD)$ vanishes on $\idem{\cU}$ if and only if it vanishes on $\cU$.
    Since $\psh_\lex(\cC) \to \psh_\lex(\idem{\cC}_{\idem{\cU}})$ is an equivalence, it follows from \cref{prop:psh-lex-univ-prop} that the lower horizontal arrow is an equivalence, finishing the proof of \eqref{lem:pass-to-completions-3}.
\end{proof}

\begin{proof}[Proof of \cref{thm:localisation}]
    By \cref{lem:pass-to-completions}, it suffices to show that
    \[ \K(\idem{\cU}) \to \K(\idem{\cC}_{\idem{\cU}}) \to \K(\quotw{\idem{\cC}_{\idem{\cU}}}{\idem{\cU}}) \]
    is a cofibre sequence.
    Due to \cref{thm:localisation-intro}, this amounts to proving that $\idem{\cU}$ is strongly left filtering in $\idem{\cC}_{\idem{\cU}}$.
    In the following, we abbreviate
    \[ \widehat{\cC} := \idem{\cC}_{\idem{\cU}}. \]
    Let us first show that $\idem{\cU}$ is left filtering in $\widehat{\cC}$.
    Since every object in $\idem{\cU}$ is a retract of an object in $\cU$, it suffices to check that $\Hom_{\widehat{\cC}}(u,-)$ is weakly $\widehat{\cC}_\egr^{\idem{\cU}}$-effaceable for all $u \in \cU$.
    Let $\widehat{x} \in \widehat{\cC}$ and choose $\widehat{u} \in \idem{\cU}$ such that $\widehat{x} \oplus \widehat{u} \in \cC$.
    Consider an element $\alpha \in \pi_k\Hom_{\widehat{\cC}}(u,\widehat{x})$.
    Since $\cU$ is left filtering in $\cC$, there exists a morphism $p \colon \widehat{x} \oplus \widehat{u} \twoheadrightarrow x'$ in $\cC_\egr^\cU$ such that the image of $\alpha$ under
    \[ \pi_k\Hom_{\widehat{\cC}}(u,\widehat{x}) \to \pi_k\Hom_{\widehat{\cC}}(u,\widehat{x} \oplus \widehat{u}) \xrightarrow{p_*} \pi_k\Hom_{\widehat{\cC}}(u,x') \]
    is trivial.
    Applying the filtering assumption also to the composite $\widehat{u} \rightarrowtail \widehat{x} \oplus \widehat{u} \overset{p}{\twoheadrightarrow} x'$, we find a morphism $p' \colon x' \twoheadrightarrow x''$ in $\cC_\egr^\cU$ such that the morphism $\widehat{x} \oplus \widehat{u} \xrightarrow{p' \circ p} x''$ factors through the projection to $\widehat{x}$.
    Then the induced morphism $q \colon \widehat{x} \to x''$ is an egressive morphism in $\idem{\cC}$ whose fibre is a retract of an object of $\cU$ by \cref{cor:egressives-cancellation}.
    By construction, the induced map $q_* \colon \pi_k\Hom_{\widehat{\cC}}(u,\widehat{x}) \to \pi_k\Hom_{\widehat{\cC}}(u,x'')$ annihilates $\alpha$.
    Hence $\idem{\cU}$ is left filtering in $\widehat{\cC}$.

    Consider now an exact sequence $\widehat{x} \rightarrowtail \widehat{y} \twoheadrightarrow \widehat{z}$ in $\widehat{\cC}$ and an egressive morphism $\widehat{y} \twoheadrightarrow \widetilde{y}$ with fibre $\widehat{v} \in \idem{\cU}$.
    Choose $\widehat{u}, \widehat{w} \in \idem{\cU}$ such that $\widehat{x} \oplus \widehat{u} \in \cC$ and $\widehat{z} \oplus \widehat{w} \in \cC$.
    Then
    \[ \widehat{x} \oplus \widehat{u} \overset{i}{\rightarrowtail} \widehat{y} \oplus \widehat{u} \oplus \widehat{w}  \twoheadrightarrow \widehat{z} \oplus \widehat{w} \]
    is an exact sequence in $\cC$.
    Pick $\overline{v} \in \idem{\cU}$ such that $\widehat{v} \oplus \widehat{u} \oplus \widehat{w} \oplus \overline{v} \in \cU$.
    Since $\cU$ is left filtering in $\cC$, the induced morphism
    \[ \widehat{v} \oplus \widehat{u} \oplus \widehat{w} \oplus \overline{v} \to \widehat{y} \oplus \widehat{u} \oplus \widehat{w} \]
    becomes nullhomotopic upon composition with some morphism $q \colon \widehat{y} \oplus \widehat{u} \oplus \widehat{w} \twoheadrightarrow \overline{y}$ in $\cC_\egr^\cU$ and which factors over $\widetilde{y}$.
    
    As $\cU$ is strongly filtering, there exist a morphism $g \colon \overline{y} \to y''$ in $\cC$ with $gq \in \cC_\egr^\cU$, a morphism $p \colon \widehat{x} \oplus \widehat{u} \twoheadrightarrow x''$ in $\cC_\egr^\cU$ and an ingressive morphism $j' \colon x'' \rightarrowtail y''$ in $\cC$ which fit into the top commutative square of the diagram
    \[\begin{tikzcd}
     \widehat{x} \oplus \widehat{u}\ar[r, tail, "i"]\ar[d, two heads, "p"'] & \widehat{y} \oplus \widehat{u} \oplus \widehat{w}\ar[d, two heads, "gq"] \\
     x''\ar[r, tail, "j'"]\ar[d, two heads, "p'"'] & y''\ar[d, two heads, "q'"] \\
     x'\ar[r, tail, "j"] & y'
    \end{tikzcd}\]
    The morphism $p' \colon x'' \twoheadrightarrow x'$ in $\widehat{\cC}_\egr^{\idem{\cU}}$ is then chosen such that the component $\widehat{u} \to x''$ of $p$ becomes nullhomotopic after composition with $p'$.
    The bottom commutative square is obtained from $p'$ by taking a pushout.
    In particular, $q'$ is egressive with fibre in $\idem{\cU}$.
    By construction, the composite $p'p$ factors over the projection $\widehat{x} \oplus \widehat{u} \twoheadrightarrow \widehat{x}$, which is egressive with fibre in $\idem{\cU}$.
    By \cref{cor:egressives-cancellation}, the induced morphism $\overline{p} \colon \widehat{x} \to x'$ is also egressive in $\widehat{\cC}$ with fibre in $\idem{\cU}$.
    Therefore, we obtain a commutative diagram
    \[\begin{tikzcd}[row sep=1.5em]
     \widehat{x}\ar[r, tail]\ar[d, tail] & \widehat{y}\ar[d, tail] \\
     \widehat{x} \oplus \widehat{u}\ar[r, tail, "i"]\ar[d, two heads, "\pr"']\ar[dd, bend left=30, "p'p"] & \widehat{y} \oplus \widehat{u} \oplus \widehat{w}\ar[dd, two heads, "q'gq"] \\
     \widehat{x}\ar[d, two heads, "\overline{p}"'] &  \\
     x'\ar[r, tail, "j"] & y'
    \end{tikzcd}\]
    Recalling that $q'gq$ factors over the egressive morphism $\widehat{y} \oplus \widehat{u} \oplus \widehat{w} \twoheadrightarrow \widetilde{y}$, the outer square produces the required factorisation.
\end{proof}

\section{Applications}\label{sec:applications}
In this section, we will see that \cref{thm:localisation-intro,thm:localisation} subsume a number of well-known localisation results in algebraic K-theory.
The first immediate consequence is the localisation theorem for stable categories, see eg \cite[Theorem~6.1]{hls:localisation}.

\begin{prop}\label[prop]{prop:localisation-stable}
 Let $\cC$ be a stable category and let $\cU \subseteq \cC$ be a full stable subcategory which is closed under retracts in $\cC$. Then
 \[ \K(\cU) \to \K(\cC) \to \K(\cC/\cU) \]
 is a cofibre sequence of connective spectra, where $\cC/\cU$ denotes the Verdier quotient of $\cC$ by $\cU$.
\end{prop}
\begin{proof}
    Since stable categories are weakly idempotent complete, \cref{cor:retracts-partial-completion} shows that $\cU \simeq \Idemrel{\K_0(\cC)}{\cU}$.
    Any stable subcategory is strongly left filtering, so we obtain from \cref{thm:localisation} a cofibre sequence
    \[ \K(\cU) \to \K(\cC) \to \K(\quotw{\cC}{\cU}). \]
    By \cref{thm:W-quotient}, the Waldhausen quotient is the Dwyer--Kan localisation of $\cC$ at the subcategory of morphisms whose fibre lies in $\cU$, which is a description of the Verdier quotient $\cC/\cU$, see \cite[Theorem~I.3.3]{ns:on-tc}.
\end{proof}

For the sake of completeness, let us point out that the localisation theorem for stable categories easily generalises to a localisation theorem for inclusions of left or right special subcategories.

\begin{prop}\label[prop]{prop:localisation-special}
    Let $\cC$ be an exact category and let $\cU$ be an extension-closed subcategory which is left special or right special.
    \begin{enumerate}
        \item If $\cU$ is idempotent complete, then
        \[ \K(\cU) \to \K(\cC) \to \K(\quotex{\cC}{\cU}) \]
        is a cofibre sequence of connective spectra.
        \item Without additional assumptions on $\cU$, there exists a cofibre sequence
        \[ \K(\Idemrel{\K_0(\cC)}{\cU}) \to \K(\cC) \to \K(\quotex{\cC}{\cU}). \]
    \end{enumerate}
\end{prop}
\begin{proof}
    As explained in \cite[Corollary~5.4]{nw:gabriel-quillen}, the first assertion follows from the localisation theorem for stable categories and the Gillet--Waldhausen theorem \cite[Theorem~1.7]{sw:exact-stable}.

   If $\cU$ is not idempotent complete, \cref{lem:pass-to-completions}.\eqref{lem:pass-to-completions-2} yields a bicartesian square
   \begin{equation}\label{eq:localisation-special}\begin{tikzcd}
         \K(\Idemrel{\K_0(\cC)}{\cU})\ar[r]\ar[d] & \K(\idem{\cU})\ar[d] \\
         \K(\idem{\cC}_{\Idemrel{\K_0(\cC)}{\cU}})\ar[r] &  \K(\idem{\cC}_{\idem{\cU}})
    \end{tikzcd}\end{equation}
    We claim that $\idem{\cU}$ is left special in $\idem{\cC}_{\idem{\cU}}$.
    For an egressive morphism $p \colon \widehat{x} \twoheadrightarrow \widehat{w}$ in $\idem{\cC}_{\idem{\cU}}$ with $\widehat{w} \in \idem{\cU}$, pick $w' \in \idem{\cU}$ and $\widehat{u} \in \idem{\cU}$ such that $\widehat{w} \oplus w' \in \cU$ and $\fib(p) \oplus \widehat{u} \in \cC$.
    Then the morphism $p' \colon \widehat{x} \oplus w' \oplus \widehat{u} \twoheadrightarrow \widehat{w} \oplus w'$ induced by $p$, $\id_{w'}$ and $\widehat{u} \to 0$ is egressive in $\cC$.
    Since $\cU$ is left special, there exists a commutative diagram
    \[\begin{tikzcd}
        v\ar[d]\ar[dr, two heads, "q'"] & \\
        \widehat{x} \oplus w' \oplus \widehat{u}\ar[r, two heads, "p'"]& \widehat{w} \oplus w'
    \end{tikzcd}\]
    with $q$ an egressive morphism in $\cU$.
    Pulling back the inclusion $\widehat{w} \rightarrowtail \widehat{w} \oplus w'$ along the egressives $p'$ and $q$ then yields a commutative diagram
    \[\begin{tikzcd}
        \widehat{v}\ar[d]\ar[dr, two heads, "q"] & \\
        \widehat{x} \oplus \widehat{u}\ar[r, two heads]& \widehat{w}
    \end{tikzcd}\]
    with $q$ an egressive morphism in $\idem{\cU}$.
    The lower horizontal map factors over $p \colon \widehat{x} \twoheadrightarrow \widehat{w}$ by construction, so this yields the required factorisation.

    The first assertion combined with \cref{lem:pass-to-completions}.\eqref{lem:pass-to-completions-3} identifies the vertical cofibres in \eqref{eq:localisation-special} with $\K(\quotex{\cC}{\cU})$.
    Moreover, we have $\K(\cC) \simeq \K(\idem{\cC}_{\Idemrel{\K_0(\cC)}{\cU}})$ by \cref{lem:pass-to-completions}.\eqref{lem:pass-to-completions-1}, which yields the desired cofibre sequence.
\end{proof}

Next, we are able to recover Schlichting's localisation theorem for s-filtering subcategories.

\begin{prop}[{Schlichting, \cite[Theorem~2.1]{schlichting:delooping}}]
    Let $\cC$ be an ordinary exact category and let $\cU$ be an idempotent complete, left special and left filtering subcategory.
    Denoting by $\cQ$ the quotient exact category of $\cC$ by $\cU$ (in the $(2,1)$-category of ordinary exact categories), there exists a cofibre sequence
    \[ \K(\cU) \to \K(\cC) \to \K(\cQ) \]
    of connective spectra.
\end{prop}
\begin{proof}
    As observed by Schlichting, there are two possible proofs for this statement.
    Note that \cref{thm:W-quotient} and \cref{prop:special-quotient} show that $\quotex{\cC}{\cU} \simeq \cQ$.

    Since $\cU$ is left special, \cref{prop:localisation-special} shows that $\K(\cU) \to \K(\cC) \to \K(\quotex{\cC}{\cU})$ is a cofibre sequence of connective spectra.
    
    Alternatively, one can combine \cref{lem:sp-filtering} with \cite[Lemma~1.17~(4)]{schlichting:delooping} to see that any idempotent complete, left special and left filtering subcategory of an ordinary exact category is strongly left filtering.
    Thus, we can appeal directly to \cref{thm:localisation}.
\end{proof}

Even though the proof does not require any of our main results, let us also observe that Sarazola's localisation theorem from \cite{sarazola:cotorsion} is a consequence of the localisation theorem for special subcategories.

\begin{prop}[{Sarazola, \cite[Theorem~6.1]{sarazola:cotorsion}}]
    Let $\cC$ be a weakly idempotent complete, ordinary exact category with enough injective objects and let $\cU \subseteq \cC$ be a full subcategory which contains all injectives and has the property that if any two objects in an exact sequence $x \rightarrowtail y \twoheadrightarrow z$ lie in $\cU$, then so does the third.

    Define $w_\cU$ as the collection of morphisms which are a composition of an ingressive morphism with cofibre in $\cU$ followed by an egressive morphism with injective fibre.

    Then the Dwyer--Kan localisation $\cC[w_\cU^{-1}]$ is a pointed and finitely cocomplete category, and
    \[ \K(\cU) \to \K(\cC) \to \K(\cC[w_\cU^{-1}]) \]
    is a cofibre sequence of connective spectra.
\end{prop}
\begin{proof}
 Sarazola shows in \cite[Theorem~4.1]{sarazola:cotorsion} that the ingressives together with the subcategory $w_\cU$ define the structure of a category of cofibrant objects on $\cC$.
 It follows from \cite[Proposition~7.5.6]{cisinski:higher-cats} that $\cC[w_\cU^{-1}]$ is pointed and finitely cocomplete, and the localisation functor $\cC \to \cC[w_\cU^{-1}]$ preserves pushouts along ingressives (in fact, it is the initial functor with this property, see \cite[Proposition~7.5.11]{cisinski:higher-cats}).

 Any exact functor with target an exact category inverts $w_\cU$, so the exact cofibre $\pi \colon \cC \to \quotex{\cC}{\cU}$ factors over $\cC[w_\cU^{-1}]$.
 Every cofibre sequence in $\cC[w_\cU^{-1}]$ is equivalent to the image of an exact sequence in $\cC$, so the comparison functor $\gamma \colon \cC[w_\cU^{-1}] \to \quotex{\cC}{\cU}$ sends cofibre sequences to exact sequences.
 In particular, the morphism $\gamma(x) \to 0$ is ingressive in $\quotex{\cC}{\cU}$ for all $x \in \cC[w_\cU^{-1}]$.
 Since $\quotex{\cC}{\cU}$ is generated by the essential image of $\pi$ under extensions, \cref{cor:five-lemma} implies that the morphism $y \to 0$ is ingressive for every object $y \in \quotex{\cC}{\cU}$.
 By \cref{cor:prestable-exact}, the category $\quotex{\cC}{\cU}$ is prestable and equipped with the maximal Waldhausen structure.
 
 The composite functor $\cC[w_\cU^{-1}] \to \quotex{\cC}{\cU} \to \derb(\quotex{\cC}{\cU})$ induces an exact functor $\stab(\cC[w_\cU^{-1}]) \to \derb(\quotex{\cC}{\cU})$.
 Since the composite $\cC \to \cC[w_\cU^{-1}] \to \stab(\cC[w_\cU^{-1}])$ vanishes on $\cU$, we also obtain an induced exact functor $\quotex{\cC}{\cU} \to \stab(\cC[w_\cU^{-1}])$ which in turn induces an exact functor $\derb(\quotex{\cC}{\cU}) \to \stab(\cC[w_\cU^{-1}])$ because the target is stable.
 Unwinding universal properties, one finds that these are mutually inverse equivalences $\derb(\quotex{\cC}{\cU}) \simeq \stab(\cC[w_\cU^{-1}])$.
 Moreover, note that $\derb(\quotex{\cC}{\cU}) \simeq \stab(\quotex{\cC}{\cU})$.

 Since $\cC$ has enough injectives and $\cU$ is closed under cofibres of ingressives, it is immediate that $\cU$ is a right special subcategory of $\cC$. 
 Now \cref{prop:localisation-special} yields a cofibre sequence
 \[ \K(\Idemrel{\K_0(\cC)}{\cU}) \to \K(\cC) \to \K(\quotex{\cC}{\cU}). \]
 Since $\cD \to \stab(\cD)$ induces an equivalence in K-theory for every pointed, finitely cocomplete category, the preceding discussion implies $\K(\quotex{\cC}{\cU}) \simeq \K(\cC[w_\cU^{-1}])$.

 We are left with showing that $\cU = \Idemrel{\K_0(\cC)}{\cU}$.
 To this end, we claim that $w_\cU$ satisfies the two-out-of-six property.
 Assuming this for a moment, if $x \in \cC$ is a retract of an object $u \in \cU$, the two-out-of-six property applied to $0 \to x \to u \to x$ implies that $0 \to x$ is in $w_\cU$, and therefore $x \in \cU$ by \cite[Lemma~5.1]{sarazola:cotorsion}.
 Since $\cC$ is weakly idempotent complete, \cref{cor:retracts-partial-completion} implies that $\cU = \Idemrel{\K_0(\cC)}{\cU}$.

 Consequently, it suffices to show the two-out-of-six property for $w_\cU$.
 Since morphisms in $w_\cU$ are closed under extensions \cite[Proposition~5.4]{sarazola:cotorsion} and an object $x \in \cC$ lies in $\cU$ if and only if $0 \to x$ is a morphism in $w_\cU$ \cite[Lemma~5.1]{sarazola:cotorsion}, it follows that an ingressive morphism is an equivalence if and only if its cofibre lies in $\cU$.
 Now let $x_1 \xrightarrow{f_1} x_2 \xrightarrow{f_2} x_3 \xrightarrow{f_3} x_4$ be morphisms such that $f_2f_1$ and $f_3f_2$ lie in $w_\cU$.
 By \cite[Proposition~3.3]{sarazola:cotorsion}, there exist factorisations $f_1 \colon x_1 \overset{i}{\rightarrowtail} x_2' \xrightarrow{g} x_2$
 and $f_3 \colon x_3 \overset{j}{\rightarrowtail} x_4' \xrightarrow{h} x_4$ with $i$ and $j$ ingressives and $g$ and $h$ in $w_\cU$.
 Then $x_1 \overset{i}{\rightarrowtail} x_2' \xrightarrow{f_2g} x_3 \overset{j}{\rightarrowtail} x_4'$ is a sequence of composable morphisms such that $f_2gi \simeq f_2f_1$ is in $w_\cU$, and also $jf_2g$ is in $w_\cU$ because $hjf_2 \simeq f_3f_2$ is in $w_\cU$ and $w_\cU$ satisfies the two-out-of-three property \cite[Proposition~5.6]{sarazola:cotorsion}.
 We conclude from the commutative diagram
 \[\begin{tikzcd}
     x_1\ar[r, tail, "i"]\ar[d, "f_2f_1"', "\sim"] & x_2'\ar[d, "\sim"', "jf_2g"] \\
     x_3\ar[r, tail, "j"] & x_4'
 \end{tikzcd}\]
 that the zero morphism $0 \colon \cofib(i) \to \cofib(j)$ is in $w_\cU$.
 Taking the pushout along the identity on $\cofib(i)$, it follows that the summand inclusion $\cofib(i) \to \cofib(i) \oplus \cofib(j)$ is also in $w_\cU$.
 Consequently, $\cofib(j) \in \cU$, which means that $j$ lies in $w_\cU$.
 By the two-out-of-three-property, we find that $f_1$, $f_2$ and $f_3$ all lie in $w_\cU$.
\end{proof}

The examples of left filtering subcategories which arise as kernels of Bousfield localisations are also strongly left filtering.

\begin{lem}\label[lem]{lem:adjoint-strongly-filtering}
    Let $\cC$ be an exact category and let $L \colon \cC \to \cD$ be a Bousfield localisation of the underlying category.
    Assume one of the following:
    \begin{enumerate}
     \item\label{lem:adjoint-strongly-filtering-1} Every morphism in $\cC$ is egressive.
     \item\label{lem:adjoint-strongly-filtering-2} $\cC$ is an ordinary category and the unit morphism $x \to Lx$ is egressive for    every $x \in \cC$.
    \end{enumerate}
    Then $\cU := \ker(L)$ is strongly left filtering and $\cU \simeq \Idemrel{\K_0(\cC)}{\cU}$.
\end{lem}
\begin{proof}
    We already know from \cref{lem:adjoint-filtering} that $\cU$ is filtering, so consider an ingressive morphism $i \colon x \rightarrowtail y$ and a morphism $q \colon y \twoheadrightarrow y'$ in $\cC_\egr^\cU$.
    The unit transformation, which has only egressive components, gives rise to a commutative diagram
    \[\begin{tikzcd}
        x\ar[r, tail, "i"]\ar[d, two heads] & y\ar[r, two heads, "q"]\ar[d, two heads] & y'\ar[d, two heads] \\
        L(x)\ar[r, tail, "L(i)"] & L(y)\ar[r, "L(q)"] & L(y')
    \end{tikzcd}\]
    Since $q$ lies in $\cC_\egr^{\cU}$, the morphism $L(q)$ is an equivalence.
    Hence the outer square provides the required factorisation.

    Suppose now that $\widehat{u} \in \idem{\cU}$ satisfies $[\widehat{u}] \in K_0(\cC)$, and recall from \cref{lem:idemrel} that $\widehat{u} \in \widem{\cC}$.
    Then there exist $x,y,a,b \in \cC$ and $\widehat{s} \in \widem{\cC}$ together with exact sequences $a \rightarrowtail \widehat{u} \oplus y \oplus \widehat{s} \twoheadrightarrow b$ and $a \rightarrowtail x \oplus \widehat{s} \twoheadrightarrow b$.
    In particular, $x \oplus \widehat{s}$ and $\widehat{u} \oplus y \oplus \widehat{s}$ lie in $\cC$.
    The obvious inclusion map induces a commutative square
    \[\begin{tikzcd}
        y \oplus \widehat{s} \oplus x\ar[r, tail]\ar[d, two heads] & \widehat{u} \oplus y \oplus \widehat{s} \oplus x\ar[d, two heads] \\
        L(y \oplus \widehat{s} \oplus x)\ar[r] & L(\widehat{u} \oplus y \oplus \widehat{s} \oplus x)
    \end{tikzcd}\]
    in which the vertical maps are the unit maps of the adjunction.
    Since $L(\widehat{u} \oplus y \oplus \widehat{s} \oplus x) \simeq \widem{L}(\widehat{u}) \oplus L(y \oplus \widehat{s} \oplus x)$, the lower horizontal map is an equivalence.
    By assumption, the fibres of the vertical maps lie in $\cU$.
    Since $\cU$, being the kernel of $L$, is closed under cofibres of ingressives, it follows that $\widehat{u}$ lies in $\cU$.
\end{proof}

From this, we can deduce Barwick's theorem of the heart.

\begin{prop}[{Theorem of the heart, \cite[Theorem~6.1]{barwick:heart}}]\label[prop]{prop:heart}
 Let $\cC$ be a small stable category with a t-structure.
 \begin{enumerate}
  \item\label{prop:heart-1} If the t-structure is bounded below, then both
   \[ \K(\cC_{\geq 0}) \to \K(\cC) \quad\text{and}\quad \K(\cC^\heartsuit) \to \K(\cC_{\leq 0}) \]
   are equivalences.
  \item\label{prop:heart-2} If the t-structure is bounded, then $\K(\cC^\heartsuit) \to \K(\cC)$ is an equivalence.
 \end{enumerate}
\end{prop}
\begin{proof}
    If the t-structure is bounded below, the inclusion functor $\cC_{\geq 0} \to \cC$ exhibits $\cC$ as a Spanier--Whitehead stabilisation of $\cC_{\geq 0}$.
    Since K-theory is invariant under passage to the Spanier--Whitehead stabilisation, the induced map $\K(\cC_{\geq 0}) \to \K(\cC)$ is an equivalence.

    \cref{lem:adjoint-strongly-filtering} shows that $\cC^\heartsuit$ is a strongly left filtering subcategory of $\cC_{\leq 0}$, and \cref{ex:w-quotients} identifies the Waldhausen quotient as $\cC^\mx_{\leq -1}$, where the superscript indicates the Waldhausen structure in which every morphism is ingressive.
    Since the t-structure is bounded below, the Spanier--Whitehead stabilisation of $\cC_{\leq -1}^\mx$ is trivial, so \cref{thm:localisation-intro} implies that $\K(\cC^\heartsuit) \to \K(\cC_{\leq 0})$ is an equivalence.
    This proves \eqref{prop:heart-1}.

    Since K-theory, as a functor on exact categories, is invariant under passage to opposite categories (see \cite[Corollary~5.16.1]{barwick:heart} for a highly coherent version of this statement), the map $\K(\cC_{\leq 0}) \to \K(\cC)$ is identified with the map $\K((\cC^\op)_{\geq 0}) \to \K(\cC^\op)$.
    If the t-structure on $\cC$ is bounded above, \eqref{prop:heart-1} implies that this map is an equivalence, which establishes \eqref{prop:heart-2}.
\end{proof}

\begin{rem}
    The proof of \cref{prop:heart} illustrates that \cref{thm:localisation-intro,thm:localisation} are genuinely statements about \emph{connective} algebraic K-theory.
    Note that all three Waldhausen categories $\cC^\heartsuit$, $\cC_{\leq 0}$ and $\cC_{\leq -1}^\mx$ have reasonably defined non-connective K-theory (being either exact or pointed and finitely cocomplete).
    However, if \cref{thm:localisation} had an analogue for non-connective K-theory, this would immediately imply the non-connective version of the theorem of the heart, which is known to fail by \cite[Theorem~B]{rsw:heart}. 
\end{rem}

\begin{rem}
    I feel obliged to point out the following gap in the published proof of \cite[Theorem~6.1]{barwick:heart}: using the notation of \cref{prop:heart}, one finds at the bottom of page 2180 an objectwise construction of factorisations in a certain labelled Waldhausen structure on $(\cC_{\geq 0})^\op$.
    Unfortunately, the given construction does not extend to morphisms, obstructing the use of the Generic Fibration Theorem II \cite[Theorem~6.4]{barwick:heart}.

    The argument for \cref{prop:heart} proceeds differently by using the localisation theorem for filtering subcategories.
    Note that the simplified proof of the localisation theorem described in \cref{sec:the-easy-proof} is sufficient for the proof of \cref{prop:heart}.
    In fact, upon specialising to the case $\cC = \ker(\Perf(\mathrm{ku}) \to \Perf(\mathrm{KU})$), this argument is virtually identical to the proof in \cite{bm:localisation-K-of-KU} establishing the fibre sequence $\K(\bbZ) \to \K(\mathrm{ku}) \to \K(\mathrm{KU})$.
\end{rem}

To conclude, we also reprove Quillen's localisation theorem for abelian categories.
A special case of the following argument can be found in \cite{efimov:devissage}, where Efimov shows that Quillen's d\'evissage theorem is a consequence of the localisation theorem.

\begin{prop}[{Quillen's localisation theorem, \cite[Theorem~5]{quillen:k-theory}}]\label[prop]{prop:quillen-localisation}
 Let $\cA$ be a small abelian category and let $\cU \subseteq \cA$ be a Serre subcategory. Then
 \[ \K(\cU) \to \K(\cA) \to \K(\cA/\cU) \]
 is a cofibre sequence of connective spectra, where $\cA/\cU$ denotes the abelian quotient category.
\end{prop}
\begin{proof}
    Using eg \cite[Section~7.4]{bckw:novikov}, we find that for every small abelian category $\cB$, the stable envelope $\derb(\cB)$ is equivalent to the Dwyer--Kan localisation of the ordinary category of bounded chain complexes at the subcategory of quasi-isomorphisms. 

    Since $\cU \to \cA \to \cA/\cU$ is a cofibre sequence of exact categories, the induced functor $p \colon \derb(\cA) \to \derb(\cA/\cU)$ is a Verdier projection, and the localisation theorem for stable categories yields a cofibre sequence
    \[ \K(\ker(p)) \to \K(\derb(\cA)) \xrightarrow{\K(p)} \K(\derb(\cA/\cU)). \]
    Either by the Gillet--Waldhausen theorem \cite[Theorem~1.11.7]{TT} or by the theorem of the heart \ref{prop:heart}, the middle and right hand terms are identified as $\K(\cA)$ and $\K(\cA/\cU)$, respectively.
    
    The kernel of $p$ is given by the full subcategory $\derb_\cU(\cA)$ of complexes whose homology lies in $\cU$.
    Hence the standard t-structure on $\derb(\cA)$ restricts to a bounded t-structure on $\derb_\cU(\cA)$ whose heart is equivalent to $\cU$.
    Applying \cref{prop:heart} again, the fibre in K-theory becomes identified with $\K(\cU)$. 
\end{proof}

\begin{rem}
    It is possible to derive \cref{prop:quillen-localisation} directly from the localisation theorem for filtering subcategories without recourse to derived categories and the theorem of the heart.

    Any Serre subcategory $\cU$ of an abelian category $\cA$ is strongly filtering, so \cref{thm:localisation-intro} yields a cofibre sequence
    \[ \K(\cU) \to \K(\cA) \to \K(\quotw{\cA}{\cU}) \]
    because $\cU$ is idempotent complete.
    As we have seen in \cref{ex:w-quotients}, the Waldhausen cofibre $\quotw{\cA}{\cU}$ need not coincide with the abelian quotient category $\cA/\cU$, so it remains to show that the comparison functor $\gamma \colon \quotw{\cA}{\cU} \to \cA/\cU$ induces an equivalence in K-theory.
    
    Setting $\cB := \quotw{\cA}{\cU}$, one can check by hand that the functor $S_p\cB \to S_p(\cA/\cU)$ is a localisation at the collection $(S_p\cB)_\ing^0$ of ingressives with trivial cofibre. 
    This will use the fact that all localisations involved in the argument support a calculus of fractions.
    Since $(S_p\cB)_\ing^0$ is closed under pushouts and satisfies two-out-of-six,
    \cite[Proposition~6.8]{hls:localisation} implies that
    \[ \abs{ (S_p\cB)_\ing^0 } \to \iota S_p(\cA/\cU) \]
    is an equivalence for all $p$.
    Therefore, \cref{prop:identify-cofibre} implies that
    \[ \K(\gamma) \colon \K(\quotw{\cA}{\cU}) \to \K(\cA/\cU) \]
    is an equivalence.
\end{rem}

\appendix

\section{Diagram lemmas in exact categories}\label{sec:diagram-lemmas}

This appendix collects a number of useful statements about exact categories that are used frequently throughout the main paper.
Let us first recall the classical definition of an exact category.

\begin{defn}[{\cite[Definition~3.1]{barwick:heart}}]\label[defn]{def:exact}
 An \emph{exact category} is a category $\cC$ together with a wide subcategory $\cC_\ing$ of \emph{ingressive} morphisms, denoted by feathered arrows $\rightarrowtail$, and a wide subcategory $\cC_\egr$ of \emph{egressive} morphisms, denoted by double-headed arrows $\twoheadrightarrow$, such that the following holds:
 \begin{enumerate}
  \item\label{def:exact-1} $\cC$ is additive;
  \item\label{def:exact-3} for every object $x \in \cC$, the unique map $0 \to x$ is ingressive and the unique map $x \to 0$ is egressive;
  \item\label{def:exact-4} every span $z \leftarrow x \rightarrowtail y$ in $\cC$ admits a pushout and the structure morphism $z \to z \sqcup_x y$ is ingressive;
  \item\label{def:exact-5} every cospan $x \twoheadrightarrow z \leftarrow y$ in $\cC$ admits a pullback and the structure morphism $x \times_z y \to y$ is egressive;
  \item\label{def:exact-6} for a commutative square
  \[\begin{tikzcd}
   w\ar[r, "i"]\ar[d,"p"'] & y\ar[d, "q"] \\
   x\ar[r, "j"] & z
  \end{tikzcd}\]
  in $\cC$ the following are equivalent:
  \begin{enumerate}
   \item $i$ is ingressive, $p$ is egressive, and the square is a pushout;
   \item $j$ is ingressive, $q$ is egressive, and the square is a pullback.
  \end{enumerate}
 \end{enumerate}
 Any square satisfying the equivalent conditions in \eqref{def:exact-6} is called an \emph{exact square}.
 If $x \simeq 0$, we call such a square also an \emph{exact sequence}, and denote it by abuse of notation as $w \rightarrowtail y \twoheadrightarrow z$.
 
 An \emph{exact functor} $f \colon (\cC,\cC_\ing,\cC_\egr) \to (\cD,\cD_\ing,\cD_\egr)$ of exact categories is a functor $f \colon \cC \to \cD$ which preserves zero objects, ingressive morphisms and egressive morphisms and sends exact squares to exact squares.
\end{defn}

By forgetting the choice of egressive morphisms, we obtain a forgetful functor $\exact \to \wald$.
A key fact is the following characterisation of admissible pushout squares.

\begin{lem}[{\cite[Proposition~A.1]{klemenc:stablehull} and \cite[Lemma~4.7]{barwick:heart}}]\label[lem]{lem:pushouts}
    Let $\cC$ be an exact category and consider a commutative square
    \[\begin{tikzcd}
        x\ar[r, rightarrowtail, "i"]\ar[d] & y\ar[d] \\
        x'\ar[r, rightarrowtail, "j"]& y'
    \end{tikzcd}\]
    in which $i$ and $j$ are ingressive as indicated.
    Then the following are equivalent:
    \begin{enumerate}
        \item the square is bicartesian;
        \item the square is a pushout;
        \item the induced map $\cofib(i) \to \cofib(j)$ on cofibres is an equivalence.
    \end{enumerate}
    Moreover, in this case there exists an exact sequence $x \rightarrowtail x' \oplus y \twoheadrightarrow y'$.
\end{lem}

\cref{lem:pushouts} entails that the forgetful functor $\exact \to \wald$ is fully faithful \cite[Corollary~4.8.1]{barwick:heart}, so there is no ambiguity about what an exact functor between exact categories is.
Moreover, it is easy to deduce the next statements.
For ordinary exact categories, these are all well-known, and the proofs provided eg in \cite{buehler} carry over in straightforward fashion.

\begin{cor}\label[lem]{lem:exact-reedy}
    Let $\cC$ be an exact category and consider a commutative diagram
    \[\begin{tikzcd}
        x\ar[r, tail]\ar[d, "f"'] & y\ar[r, two heads]\ar[d, "g"] & z\ar[d, "h"] \\
        x'\ar[r, tail] & y'\ar[r, two heads] & z'
    \end{tikzcd}\]
    Then the left square is Reedy cofibrant (ie $x' \sqcup_x y \to y'$ is ingressive) if and only if $h$ is ingressive.
\end{cor}
\begin{proof}
    This follows by considering the same factorisation as in \cite[Proposition~3.1]{buehler}.
\end{proof}

\begin{cor}\label[cor]{cor:five-lemma}
    Consider the situation of \cref{lem:exact-reedy}.
    If $f$ and $h$ are ingressive, then $g$ is also ingressive.
\end{cor}
\begin{proof}
    The morphism $g$ factors as $y \to y \cup_x x' \to y'$.
    The first map is ingressive since it is a pushout of $f$, and the second morphism is ingressive by \cref{lem:exact-reedy}.
\end{proof}

\begin{cor}\label[cor]{lem:nine-lemma}
    Let $\cC$ be an exact category and consider a commutative diagram
    \[\begin{tikzcd}
        x\ar[r, tail]\ar[d, tail, "f"'] & y\ar[r, two heads]\ar[d, tail, "g"] & z\ar[d, tail, "h"] \\
        x'\ar[r, tail] & y'\ar[r, two heads] & z'
    \end{tikzcd}\]
    Then the induced sequence of vertical cofibres
    \[ \cofib(f) \to \cofib(g) \to \cofib(h) \]
    is also exact.
\end{cor}
\begin{proof}
    By \cref{lem:exact-reedy}, the left square in the original diagram is Reedy cofibrant, which immediately implies that $\cofib(f) \to \cofib(g)$ is ingressive.
    The sequence in question is a cofibre sequence, so we are done.
\end{proof}

\begin{cor}\label[cor]{cor:prestable-exact}
    Let $\cC$ be a category.
    Then the following are equivalent:
    \begin{enumerate}
        \item\label{cor:prestable-exact-1} $\cC$ admits an exact structure in which the map $x \to 0$ is ingressive for all $x \in \cC$;
        \item\label{cor:prestable-exact-2} $(\cC,\cC)$ is an exact category;
        \item\label{cor:prestable-exact-3} $\cC$ is prestable.
    \end{enumerate}
\end{cor}
\begin{proof}
    Let $(\cC,\cC_\ing)$ be an exact structure as in \eqref{cor:prestable-exact-1}.
    Applying \cref{cor:five-lemma} to the commutative diagram
    \[\begin{tikzcd}
        0\ar[r, tail]\ar[d, tail] & x\ar[r, two heads]\ar[d, "f"] & x\ar[d, tail] \\
        y\ar[r, tail] & y\ar[r, two heads] & 0
    \end{tikzcd}\]
    associated to any morphism $f$ shows that $\cC_\ing = \cC$.
    It follows that \eqref{cor:prestable-exact-1} and \eqref{cor:prestable-exact-2} are equivalent.

    If every morphism in $\cC$ is ingressive in an exact structure, then $\cC$ admits a suspension functor.
    Since every exact sequence $x \rightarrowtail 0 \twoheadrightarrow \Sigma x$ is also a fibre sequence, it follows that $\Sigma$ is fully faithful.
    As $0 \to \Sigma z$ is always egressive, it follows that every morphism whose target is a suspension admits a fibre, and the resulting fibre sequence $x \to y \to \Sigma z$ is also a cofibre sequence by (the dual of) \cref{lem:pushouts}.
    This shows $\eqref{cor:prestable-exact-2} \Rightarrow \eqref{cor:prestable-exact-3}$.

    If $\cC$ is prestable, the embedding $\cC \to \stab(\cC)$ into its Spanier--Whitehead stabilisation is fully faithful, preserves finite colimits, and the essential image of $\cC$ is closed under extensions \cite[Proposition~C.1.2.2]{SAG}.
    Hence $(\cC,\cC)$ is an exact category.
\end{proof}

\begin{rem}
    The first part of the proof of \cref{cor:prestable-exact} shows that if $x \to 0$ is ingressive in an exact category, then every morphism with domain $x$ is ingressive.
\end{rem}

\begin{lem}[Quillen's obscure axiom]\label[lem]{lem:obscure-axiom}
    Let $\cC$ be an exact category.
    Suppose that $f \colon x \to y$ is a morphism which admits a cofibre and that there exists an ingressive morphism $j \colon y \rightarrowtail z$ such that $jf$ is ingressive.
    Then $f$ is also ingressive.
\end{lem}
\begin{proof}
 The proof of \cite[Proposition~2.16]{buehler} applies almost verbatim.
\end{proof}

\begin{cor}\label[cor]{cor:pushout-ingressive-two-out-of-three}
    Let $\cC$ be an exact category and consider a pushout square
    \[\begin{tikzcd}
        x\ar[r, "i"]\ar[d, tail, "f"] & y\ar[d, tail, "g"] \\
        x'\ar[r, "j", tail] & y'
    \end{tikzcd}\]
    in which $f$, $g$ and $j$ are ingressive as indicated.
    Then $i$ is also ingressive.
\end{cor}
\begin{proof}
    Since the co-Yoneda embedding sends pushouts to pullbacks, the composite morphism $y \xrightarrow{g} y' \to \cofib(j)$ is a cofibre of $i$.
    Quillen's obscure axiom implies that $i$ is ingressive because $g$ and $gi \simeq jf$ are ingressive.
\end{proof}

\begin{cor}\label[cor]{cor:ingressives-widem}
    Let $\cC$ be an exact category and let $i \colon x \to y$ and $g \colon y \to z$ be morphisms in $\cC$ such that $gi$ is ingressive.
    Then $i$ is ingressive in the weak idempotent completion of $\cC$.
\end{cor}
\begin{proof}
    The pushout
    \[\begin{tikzcd}
     x\ar[r, tail, "gi"]\ar[d, "i"'] & z\ar[d, "j"] \\
     y\ar[r, tail] & z'
    \end{tikzcd}\]
    admits a retraction $z' \to z$.
    Consequently, $j$ is ingressive in $\idem{\cC}$, and $z' \simeq z \oplus \cofib(j)$.
    So $\cofib(j)$ lies in $\widem{\cC}$.
    Now apply \cref{cor:pushout-ingressive-two-out-of-three}.
\end{proof}

\begin{rem}
 By dualising, one also obtains corresponding statements for egressive morphisms.
\end{rem}


\section{Rezk's equifibration criterion}\label{app:rezk}
The proof of the generic fibration theorem in \cref{sec:generic-fibration} relies on Rezk's equifibration criterion.
This was stated and proved in very convenient form as Lemma~3.3.14 in the preprint version of \cite{the-nine:2}, but was removed in the published version.
For the reader's convenience, I reproduce the argument from \textit{loc.\ cit.} here, with some minor improvements I learned from \cite[Section~4.6]{cisinski:higher-cats}.
This exhibits the equifibration criterion entirely as a consequence of the (un)straightening equivalence for left fibrations.
I claim no originality.

To be precise, we will use the following version of unstraightening, where $\LFib(I) \subseteq \catinf_{/I}$ denotes the full subcategory spanned by the left fibrations over $I$.

\begin{theorem}[{Straightening/Unstraightening, \cite[Theorem~3.2.0.1]{HTT}, \cite{hhr:straightening}}]\label[theorem]{thm:straight}
 For each small category $I$, there exists an equivalence
 \[ \Un^\cc \colon \Fun(I,\Spc) \xrightarrow{\sim} \LFib(I) \]
 such that for every functor $f \colon I \to J$, there exists a commutative square
 \[\begin{tikzcd}
  \Fun(J,\Spc)\ar[r, "\Un^\cc"]\ar[d, "f^*"'] & \LFib(J)\ar[d, "f^*"] \\
  \Fun(I,\Spc)\ar[r, "\Un^\cc"] & \LFib(I)
 \end{tikzcd}\]
 where the left vertical arrow is the restriction functor along $f$, and the right vertical functor is given by pullback along $f$.
 
 Moreover, $\Un^\cc$ is inverse to the forgetful equivalence $\Spc_{/*} \to \Fun(*,\Spc) \simeq \Spc$ for $I = *$.
\end{theorem}

We refer to $\Un^\cc$ as the \emph{(cocartesian) unstraightening functor}.
Its inverse, the \emph{(cocartesian) straightening functor}, will be denoted by
\[ \St^\cc \colon \LFib(I) \to \Fun(I,\Spc). \]
The following consequences of the unstraightening equivalence are immediate:
\begin{enumerate}
    \item if $X \colon I \to \Spc$ is a functor, then the evident functor
    \[ \Un^\cc(X) \to I \times \abs{\Un^\cc(X)} \]
    exhibits $\abs{\Un^\cc(X)}$ as a colimit of $X$ \cite[Corollary~3.3.4.6]{HTT};
    \item colimits are universal in anima: for every morphism $f \colon A \to B$ in $\Spc$, the induced pullback functor
    \[ f^* \colon \Spc_{/B} \to \Spc_{/A} \]
    preserves colimits \cite[Lemma~6.1.3.14]{HTT}.
\end{enumerate}

\begin{defn}
 Let $X, Y \colon I \to \Spc$. A natural transformation $\tau \colon X \Rightarrow Y$ is \emph{equifibred} if the commutative square
 \[\begin{tikzcd}
  X(i)\ar[r, "X(\alpha)"]\ar[d, "\tau(i)"'] & X(j)\ar[d, "\tau(j)"] \\
  Y(i)\ar[r, "Y(\alpha)"] & Y(j)
 \end{tikzcd}\]
 is a pullback for every morphism $\alpha \colon i \to j$ in $I$.
\end{defn}

\begin{defn}
 A left fibration $p \colon \cC \to I$ is \emph{locally constant} if the straightened functor $\St^\cc(p) \colon I \to \Spc$ inverts every morphism in $I$.
\end{defn}

\begin{lem}\label[lem]{lem:equifib-lc}
 Let $X, Y \colon I \to \Spc$ be functors and let $\tau \colon X \Rightarrow Y$ be a natural transformation.
 Then $\tau$ is equifibred if and only if the left fibration $\Un^\cc(\tau) \colon \Un^\cc(X) \to \Un^\cc(Y)$ is locally constant.
\end{lem}
\begin{proof} 
 Let $\alpha \colon i \to j$ be a morphism in $I$, which we regard as a functor $\alpha \colon [1] \to I$.
 By the naturality of the unstraightening equivalence, the pullback \[ \alpha^*\Un^\cc(\tau) \colon \alpha^*\Un^\cc(X) \to \alpha^*\Un^\cc(Y)\]
 straightens to the natural transformation $X \circ \alpha \Rightarrow Y \circ \alpha$ represented by the commutative square
  \[\begin{tikzcd}
  X(i)\ar[r, "X(\alpha)"]\ar[d, "\tau(i)"'] & X(j)\ar[d, "\tau(j)"] \\
  Y(i)\ar[r, "Y(\alpha)"] & Y(j)
 \end{tikzcd}\]
 in $\Spc$.
 Similarly, a functor $\eta \colon [1] \to \alpha^*\Un^\cc(Y)$ over $[1]$ straightens to a natural transformation $* \Rightarrow Y \circ \alpha$.
 
 Since (un)straightening is an equivalence, it preserves pullbacks.
 It follows that the left fibration
 \[ [1] \times_{\alpha^*\Un^\cc(Y)} \alpha^*\Un^\cc(X) \to [1] \]
 straightens to the canonical map
 \begin{equation}\label{eq:equifib-lc}
     \fib_y(\tau(i)) \to \fib_{y'}(\tau(j)).
 \end{equation}
 Therefore, $\Un^\cc(\tau)$ is locally constant if and only if the comparison map \eqref{eq:equifib-lc} is an equivalence for every choice of $\alpha$ and $\eta$.
 This in turn is the same as the assertion that $\tau$ is equifibred.
\end{proof}

\begin{lem}\label[lem]{lem:realise-constant-fib}
 Let $p \colon \cC \to \cD$ be a locally constant left fibration. Then the commutative square
 \[\begin{tikzcd}
  \cC\ar[r]\ar[d, "p"'] & \abs{\cC}\ar[d, "\abs{p}"] \\
  \cD\ar[r] & \abs{\cD}
 \end{tikzcd}\]
 is a pullback.
\end{lem}
\begin{proof}
 Straightening $p$, we obtain a functor $\St^\cc(p) \colon \cD \to \Spc$ which factors over $\abs{\cD}$ since $p$ is locally constant.
 Call the induced functor $\overline{p} \colon \abs{\cD} \to \Spc$.
 Unstraightening $\overline{p}$, we obtain from the naturality of unstraightening a pullback square
  \[\begin{tikzcd}
  \cC\ar[r]\ar[d, "p"'] & \Un^\cc(\overline{p})\ar[d, "\abs{p}"] \\
  \cD\ar[r] & \abs{\cD}
 \end{tikzcd}\]
 Since $\abs{\cD}$ is an anima and $\overline{p}$ is a left fibration, $\Un^\cc(\overline{p})$ is also an anima.
 
 The induced map $\abs{\cC} \to \Un^\cc(\overline{p})$ is equivalent to the canonical map
 \[ \colim_\cD \St^\cc(p) \to \colim_{\abs{\cD}} \overline{p}. \]
 Since $\cD \to \abs{\cD}$ is a localisation, it is cofinal.
 Hence this map is an equivalence, and the lemma follows.
 \end{proof}

\begin{lem}\label[lem]{lem:equifib-colim}
 Let $X^\rhd, Y^\rhd \colon I^\rhd \to \Spc$ be diagrams and let $\tau^\rhd \colon X^\rhd \Rightarrow Y^\rhd$ be a natural transformation.
 Assume that $\tau \colon X \Rightarrow Y$ is equifibred and that $Y^\rhd$ is a colimit diagram.
 
 Then $\tau^\rhd$ is equifibred if and only if $X^\rhd$ is a colimit diagram.
\end{lem}
\begin{proof}
 Suppose that $X^\rhd$ is a colimit diagram.
 \cref{lem:equifib-lc} implies that $\Un^\cc(\tau)$ is locally constant.
 Using the formula for colimits in $\Spc$, we obtain from \cref{lem:realise-constant-fib} a pullback square
 \[\begin{tikzcd}
  \Un^\cc(X)\ar[r]\ar[d, "\Un^\cc(\tau)"'] & X(\infty)\ar[d, "\tau(\infty)"] \\
  \Un^\cc(Y)\ar[r] & Y(\infty)
 \end{tikzcd}\]
 By pasting with the pullback diagram
  \[\begin{tikzcd}
  X(i)\ar[r]\ar[d, "\tau(i)"'] & \Un^\cc(X)\ar[d, "\Un^\cc(\tau)"] \\
  Y(i)\ar[r] & \Un^\cc(Y)
 \end{tikzcd}\]
 (which we obtain from the naturality of unstraightening), it follows that $\tau^\rhd$ is equifibred.
 
 Assume conversely that $\tau^\rhd$ is equifibred.
 Then each square
 \[\begin{tikzcd}
  X(i)\ar[r]\ar[d, "\tau(i)"'] & X(\infty)\ar[d] \\
  Y(i)\ar[r] & Y(\infty)
  \end{tikzcd}\]
  is a pullback, which means that $X \simeq Y \times_{{Y(\infty)}} {X(\infty)}$, where we write $X(\infty)$ and $Y(\infty)$ for the respective constant functor.
  Unstraightening, we conclude that
  \[ \Un^\cc(X) \to \Un^\cc(Y) \times_{Y(\infty)} X(\infty) \]
  is an equivalence.
  
  Moreover, we obtain from the universal property of realisation a commutative diagram
 \[\begin{tikzcd}
  \Un^\cc(Y \times_{Y(\infty)} X(\infty))\ar[r]\ar[d, "\Un^\cc(\tau)"'] & \abs{\Un^\cc(Y) \times_{Y(\infty)} X(\infty)}\ar[r]\ar[d] & X(\infty)\ar[d, "\tau(\infty)"] \\
  \Un^\cc(Y)\ar[r] & \abs{\Un^\cc(Y)}\ar[r, "\sim"] & Y(\infty)
 \end{tikzcd}\]
 By the universality of colimits, the canonical map
 \[ \abs{\Un^\cc(Y) \times_{Y(\infty)} X(\infty))} \to \abs{\Un^\cc(Y)} \times_{Y(\infty)} X(\infty) \]
 is an equivalence.
 It follows that the right square is a pullback, which implies that
 \[ \colim_I X \to X(\infty) \]
 is an equivalence.
\end{proof}

\begin{proof}[Proof of \cref{thm:rezk-equifib}]
 Since $\tau$ is equifibred, it is locally constant by \cref{lem:equifib-lc}.
 Then \cref{lem:realise-constant-fib} and the formula for colimits in $\Spc$ imply that
 \[\begin{tikzcd}
  \Un^\cc(X)\ar[r]\ar[d, "\Un^\cc(\tau)"'] & \colim_I X\ar[d] \\
  \Un^\cc(Z)\ar[r] & \colim_I Z
 \end{tikzcd}\]
 is a pullback.
 Pasting this pullback with the obvious pullback square involving the unstraightenings of $W$, $X$, $Y$ and $Z$, the outer square in the commutative diagram
 \[\begin{tikzcd}
  \Un^\cc(W)\ar[r]\ar[d, "\Un^\cc(\sigma)"'] & \colim_I Y \mathop{\times}\limits_{\colim_I Z} \colim_I X\ar[r]\ar[d] & \colim_I X\ar[d] \\
  \Un^\cc(Y)\ar[r] & \colim_I Y\ar[r] & \colim_I Z
 \end{tikzcd}\]
 is a pullback.
 Since the right square is a pullback by definition, the left square is also a pullback.
 Since $\tau$ is equifibred, the same is true for $\sigma$.
 It follows that the left square corresponds to an equifibred transformation, so \cref{lem:equifib-colim} implies that
 \[ \colim_I W \simeq \colim_I Y \mathop{\times}\limits_{\colim_I Z} \colim_I X \]
 as desired.
\end{proof}

\bibliographystyle{alpha}
\bibliography{LocalisingInvariants}

\newcommand{\etalchar}[1]{$^{#1}$}
\begin{thebibliography}{BCKW25}

\bibitem[Bar15]{barwick:heart}
C.~Barwick.
\newblock On exact {{\(\infty\)}}-categories and the theorem of the heart.
\newblock {\em Compos. Math.}, 151(11):2160--2186, 2015.

\bibitem[Bar16]{barwick:algK}
C.~Barwick.
\newblock On the algebraic {{\(K\)}}-theory of higher categories.
\newblock {\em J. Topol.}, 9(1):245--347, 2016.

\bibitem[BCKW25]{bckw:novikov}
U.~Bunke, D.-C. Cisinski, D.~Kasprowski, and C.~Winges.
\newblock Controlled objects in left-exact $\infty$-categories and the
  {N}ovikov conjecture.
\newblock {\em Bull. Soc. Math. Fr.}, 153(2):295--458, 2025.

\bibitem[BG]{bg:recollements}
C.~Barwick and S.~Glasman.
\newblock A note on stable recollements.
\newblock \href{https://arxiv.org/abs/1607.02064}{arXiv:1607.02064}.

\bibitem[BGT13]{bgt:alg-k}
A.~Blumberg, D.~Gepner, and G.~Tabuada.
\newblock A universal characterization of higher algebraic {{\(K\)}}-theory.
\newblock {\em Geom. Topol.}, 17(2):733--838, 2013.

\bibitem[BM08]{bm:localisation-K-of-KU}
A.~J. Blumberg and M.~A. Mandell.
\newblock The localization sequence for the algebraic {{\(K\)}}-theory of
  topological {{\(K\)}} -theory.
\newblock {\em Acta Math.}, 200(2):155--179, 2008.

\bibitem[B{\"u}h10]{buehler}
T.~B{\"u}hler.
\newblock Exact categories.
\newblock {\em Expo. Math.}, 28(1):1--69, 2010.

\bibitem[CDH{\etalchar{+}}25]{the-nine:2}
B.~Calm\`es, E.~Dotto, Y.~Harpaz, F.~Hebestreit, M.~Land, K.~Moi, D.~Nardin,
  T.~Nikolaus, and W.~Steimle.
\newblock Hermitian {K}-theory for stable {{\(\infty\)}}-categories. {II}:
  {Cobordism} categories and additivity.
\newblock {\em Acta Math.}, 235(2):149--400, 2025.

\bibitem[Cis19]{cisinski:higher-cats}
D.-C. Cisinski.
\newblock {\em Higher categories and homotopical algebra}, volume 180 of {\em
  Camb. Stud. Adv. Math.}
\newblock Cambridge: Cambridge University Press, 2019.

\bibitem[Efi]{efimov:devissage}
A.~Efimov.
\newblock Some remarks on {Q}uillen's {D}\'evissage theorem.
\newblock \href{https://arxiv.org/abs/2505.13260}{arXiv:2505.13260}.

\bibitem[GHN17]{ghn:lax-colimits}
D.~Gepner, R.~Haugseng, and T.~Nikolaus.
\newblock Lax colimits and free fibrations in {{\(\infty\)}}-categories.
\newblock {\em Doc. Math.}, 22:1225--1266, 2017.

\bibitem[Gla16]{glasman:hodge}
S.~Glasman.
\newblock A spectrum-level {Hodge} filtration on topological {Hochschild}
  homology.
\newblock {\em Sel. Math., New Ser.}, 22(3):1583--1612, 2016.

\bibitem[HHR25]{hhr:straightening}
F.~Hebestreit, G.~Heuts, and J.~Ruit.
\newblock A short proof of the straightening theorem.
\newblock {\em Trans. Am. Math. Soc., Ser. B}, 12:697--747, 2025.

\bibitem[HLS24]{hls:localisation}
F.~Hebestreit, A.~Lachmann, and W.~Steimle.
\newblock The localisation theorem for the {{\(K\)}}-theory of stable
  {{\(\infty\)}}-categories.
\newblock {\em Proc. R. Soc. Edinb., Sect. A, Math.}, 154(6):1749--1785, 2024.

\bibitem[HS04]{hs:localisation-hermitian-k}
J.~Hornbostel and M.~Schlichting.
\newblock Localization in {Hermitian} {{\(K\)}}-theory of rings.
\newblock {\em J. Lond. Math. Soc., II. Ser.}, 70(1):77--124, 2004.

\bibitem[HS25]{hs:rezk-nerve}
F.~Hebestreit and J.~Steinebrunner.
\newblock A short proof that {R}ezk's nerve is fully faithful.
\newblock {\em Int. Math. Res. Not.}, 2025(4):rnaf021, 2025.

\bibitem[Kas15]{kasprowski:fdc}
D.~Kasprowski.
\newblock On the {{\(K\)}}-theory of groups with finite decomposition
  complexity.
\newblock {\em Proc. Lond. Math. Soc. (3)}, 110(3):565--592, 2015.

\bibitem[Kle22]{klemenc:stablehull}
J.~Klemenc.
\newblock The stable hull of an exact {{\(\infty\)}}-category.
\newblock {\em Homology Homotopy Appl.}, 24(2):195--220, 2022.

\bibitem[KS06]{ks:cats-and-sheaves}
M.~Kashiwara and P.~Schapira.
\newblock {\em Categories and sheaves}, volume 332 of {\em Grundlehren Math.
  Wiss.}
\newblock Berlin: Springer, 2006.

\bibitem[KW20]{kw:shortening}
D.~Kasprowski and C.~Winges.
\newblock Shortening binary complexes and commutativity of {{\(K\)}}-theory
  with infinite products.
\newblock {\em Trans. Am. Math. Soc., Ser. B}, 7:1--23, 2020.

\bibitem[Lur]{SAG}
J.~Lurie.
\newblock Spectral {A}lgebraic {G}eometry.
\newblock Available on the author's homepage:
  \href{https://www.math.ias.edu/~lurie/papers/SAG-rootfile.pdf}{https://www.math.ias.edu/{\textasciitilde}lurie/papers/SAG-rootfile.pdf}.

\bibitem[Lur09]{HTT}
J.~Lurie.
\newblock {\em {Higher topos theory}}, volume 170.
\newblock Princeton, NJ: Princeton University Press, 2009.

\bibitem[Lur23]{kerodon}
J.~Lurie.
\newblock Kerodon.
\newblock \href{https://kerodon.net}{https://kerodon.net}, 2023.

\bibitem[NS18]{ns:on-tc}
T.~Nikolaus and P.~Scholze.
\newblock On topological cyclic homology.
\newblock {\em Acta Math.}, 221(2):203--409, 2018.

\bibitem[NW]{nw:gabriel-quillen}
M.~Nielsen and C.~Winges.
\newblock The presentable stable envelope of an exact category.
\newblock To appear in \textit{Algebr.~Geom.~Topol.},
  \href{https://arxiv.org/abs/2506.02598}{arXiv:2506.02598}.

\bibitem[PW89]{pw:khomology}
E.~K. Pedersen and C.~A. Weibel.
\newblock K-theory homology of spaces.
\newblock Algebraic topology, {Proc}. {Int}. {Conf}., {Arcata}/{Calif}. 1986,
  {Lect}. {Notes} {Math}. 1370, 346-361, 1989.

\bibitem[Qui73]{quillen:k-theory}
D.~Quillen.
\newblock Higher algebraic {{\(K\)}}-theory. {I}.
\newblock Algebr. {{\(K\)}}-{Theory} {I}, {Proc}. {Conf}. {Battelle} {Inst}.
  1972, {Lect}. {Notes} {Math}. 341, 85-147, 1973.

\bibitem[RSW]{rsw:heart}
M.~Ramzi, V.~Sosnilo, and C.~Winges.
\newblock Every spectrum is the {K}-theory of a stable $\infty$-category.
\newblock \href{https://arxiv.org/abs/2401.06510}{arXiv:2401.06510}.

\bibitem[Sar20]{sarazola:cotorsion}
M.~Sarazola.
\newblock Cotorsion pairs and a {{\(K\)}}-theory localization theorem.
\newblock {\em J. Pure Appl. Algebra}, 224(11):28, 2020.

\bibitem[Sch04]{schlichting:delooping}
M.~Schlichting.
\newblock Delooping the {{\(K\)}}-theory of exact categories.
\newblock {\em Topology}, 43(5):1089--1103, 2004.

\bibitem[Sta89]{staffeldt:fundamental}
R.~E. Staffeldt.
\newblock On fundamental theorems of algebraic {K}-theory.
\newblock {\em \(K\)-Theory}, 2(4):511--532, 1989.

\bibitem[SW26]{sw:exact-stable}
V.~Saunier and C.~Winges.
\newblock On exact categories and their stable envelopes.
\newblock {\em Math. Z.}, 312(1):28, 2026.
\newblock Id/No 28.

\bibitem[TT90]{TT}
R.~W. Thomason and T.~Trobaugh.
\newblock Higher algebraic {K}-theory of schemes and of derived categories.
\newblock The {Grothendieck} {Festschrift}, {Collect}. {Artic}. in {Honor} of
  the 60th {Birthday} of {A}. {Grothendieck}. {Vol}. {III}, {Prog}. {Math}. 88,
  247--435, 1990.

\bibitem[Wal85]{waldhausen:algKspaces}
F.~Waldhausen.
\newblock Algebraic {K}-theory of spaces.
\newblock Algebraic and geometric topology, {Proc}. {Conf}., {New}
  {Brunswick}/{USA} 1983, {Lect}. {Notes} {Math}. 1126, 318-419, 1985.

\end{thebibliography}

\end{document}